\documentclass[11pt]{article}
\usepackage{amsmath,amssymb,color}

\evensidemargin \oddsidemargin

\newtheorem{theorem}{Theorem}[section]

\newtheorem{corollary}[theorem]{Corollary}
\newtheorem{lemma}[theorem]{Lemma}
\newtheorem{proposition}[theorem]{Proposition}

\newtheorem{definition}{Definition}[section]

\numberwithin{equation}{section}

\def\al{\alpha}
\def\be{\beta}
\def\ga{\gamma}
\def\de{\delta}
\def\ep{\epsilon}

\def\e {\varepsilon}

\def\th{\theta}

\def\hro{\hat{\rho}}

\newcommand{\R}{\mathbb{R}}
\newcommand{\C}{\mathbb{C}}
\renewcommand{\Re}{\mathop{\mathrm{Re}}}
\renewcommand{\Im}{\mathop{\mathrm{Im}}}

\newcommand{\supp}{\mathop{\mathrm{supp}}}

\newcommand{\norm}[1]{\left\Vert #1 \right\Vert}

\newcommand{\bka}[1]{\left \langle #1 \right \rangle}
\newcommand{\bke}[1]{\left( #1 \right)}
\newcommand{\bkt}[1]{\left[ #1 \right]}
\newcommand{\bket}[1]{\left\{ #1 \right\}}

\newcommand{\wei}[1]{\langle #1 \rangle}
\newcommand{\myproof}{\noindent{\it Proof}.\quad }
\newcommand{\myendproof}{{\hfill $\square$}}
\newcommand{\myremark}{\noindent{\it Remark}.\quad }

\newcommand{\Pc}{P_c}

\newcommand{\bE}{\mathbf{E}}
\renewcommand{\L}{\mathcal{L}}

\newcommand{\loc}{_{\mathrm{loc}}}
\newcommand{\LLLL}{\lambda}
\newcommand{\La}{\Lambda}
\newcommand{\Ga}{\Gamma}
\newcommand{\Si}{\Sigma}
\newcommand{\Th}{\Theta}

\newcommand{\wt}{\widetilde}
\newcommand{\pd}{\partial}
\newcommand{\donothing}[1]{}
\newcommand{\wbar}[1]{\overline{\rule{0pt}{2.4mm} {#1}}}

\newcommand{\svect}[1]
{\bkt{\begin{smallmatrix} #1\end{smallmatrix}}}
\newcommand{\myfrac}[2]{\stackrel {\scriptstyle #1}{\scriptstyle #2}}
\newcommand{\les}{\lesssim}

\newcommand{\ev}{\omega}

\newcommand{\lec}{\lesssim}
\newcommand{\gec}{\gtrsim}

\newcommand{\td}{\tilde}
\newcommand{\om}{\omega}
\newcommand{\si}{\sigma}
\newcommand{\De}{\Delta}
\newcommand{\mat}[1]{\begin{bmatrix}#1\end{bmatrix}}
\newcommand{\smat}[1]{\bkt{\begin{smallmatrix}#1\end{smallmatrix}}}
\newcommand{\bL}{\mathbf{L}}
\newcommand{\bj}{\pmb{\jmath}}
\newcommand{\ph}{\varphi}
\newcommand{\lcsp}{L^2_{r}}
\newcommand{\olcsp}{O_{L^\infty_{3r}}}
\newcommand{\tsum}{{\textstyle \sum}}

\newcommand{\idx}{I}%{\mathcal{A}}

\newcommand{\nn}{\nonumber}
\newcommand{\I}{\infty}
\newcommand{\Om}{\Omega}
\newcommand{\Psp}{P^\sharp}
\newcommand{\Pcsp}{P_c^\sharp}
\newcommand{\dist}{\mathop{\mathrm{dist}}}
\newcommand{\mmin}{\wedge}

\newcommand{\tcp}{{t_c^+}}
\begin{document}
\title{Small solutions of nonlinear Schr\"odinger equations
near first excited states}
%\date{\jobname .tex}
\date{}
\author{Kenji Nakanishi,\quad  Tuoc Van Phan, \quad Tai-Peng Tsai}
\maketitle
\begin{abstract}
Consider a nonlinear Schr\"odinger equation in $\R^3$ whose linear
part has three or more eigenvalues satisfying some resonance
conditions. Solutions which are initially small  in $H^1 \cap
L^1(\R^3)$ and inside a neighborhood of the first excited state
family are shown to converge to either a first excited state or a
ground state at time infinity. An essential part of our analysis
is on the linear and nonlinear estimates near nonlinear excited
states, around which the linearized operators have eigenvalues
with nonzero real parts and their corresponding eigenfunctions are
not uniformly localized in space.
\end{abstract}
\tableofcontents

\section{Introduction}

Consider the nonlinear Schr\"odinger equation in $\R^3$,
\begin{equation} \label{Sch}
i \partial _t \psi = H_0 \psi + \kappa  |\psi|^2 \psi, \qquad
\psi|_{t=0}= \psi_0,
\end{equation}
where $H_0= -\Delta + V$ is the linear Hamiltonian with a
localized real potential $V$, $\kappa =\pm 1$, and $\psi(t,x):\R
\times \R^3 \to \C$ is the wave function.  We often drop the $x$
dependence and write $\psi(t)$. We assume $\psi_0 \in H^1$ is
localized, say $\psi_0 \in L^1$, so that its dispersive component
decays rapidly under the evolution. For any solution $\psi(t)\in
H^1(\R^3)$ its $L^2$-norm and energy
\begin{equation} \label{1-2}
{\cal E}[\psi] = \int \frac 12 |\nabla \psi|^2 + \frac 12 V
|\psi|^2 + \frac 14 \kappa  |\psi|^4 \, d x
\end{equation}
are constant in $t$. The global well-posedness for small solutions
in $H^1(\R^3)$ can be proven using these conserved quantities no
matter what the sign of $\kappa $ is.

We assume that $H_0$ has $K+1$ simple eigenvalues $e_0<e_1<\dots
<e_K(<0)$ with normalized real eigenfunctions $\phi_k$,
$k=0,1,\ldots,K$, where $K\ge 2$. They are assumed to satisfy
\begin{equation} \label{resonant-cond}
e_0 < 2 \, e_1 < 4 \, e_2,
\end{equation}
and some generic conditions to be specified later. Through
bifurcation around zero along these eigenfunctions, one obtains
$K+1$ families of {\it nonlinear bound states} $Q_{k,n}= n \phi_k
+h$, $h= O(n^3)$ and%
\footnote{The $L^2$ inner product $( \; , \; )$ is $ (f,g) = \int
_{\R^3} \, \bar f \, g  \, d x$. For a function $\phi \in L^2$, we
denote by $\phi^\perp$ the $L^2$-subspace $\bket{g \in L^2:
(\phi,g)=0}$. }
 $(h,\phi_k)=0$ for $k=0,\ldots,K$, and $n> 0$ sufficiently small, which solve the
equation
\begin{equation}   \label{Q.eq}
    (-\Delta + V) Q + \kappa  |Q|^2 Q = EQ,
\end{equation}
for some $E=E_{k,n}=e_k + O(n^2)$, see Lemma \ref{th:2-1}. They are
real and decay exponentially at spatial infinity. Each of them gives
an exact solution $\psi(t,x) = Q(x) e^{-i E t }$ of \eqref{Sch}. The
family $Q_{0,n}$ are called the {\it nonlinear ground states} while
$Q_{k,n}$, $k>0$, are called the $k$-th {\it nonlinear excited
states}.

Our goal is to understand the long-time dynamics of the solutions at
the presence of nonlinear bound states. The first question is the
stability problem of nonlinear  ground states. It is well-known that
nonlinear ground states are {\it orbitally stable} in the sense that
the difference
\begin{equation}
   \inf_{n,\th} \norm{ \psi(t) -Q_{0,n} \, e^{i \theta} }_{H^1(\R^3)}
\end{equation}
remains uniformly small for all time $t$ if it is initially small.
On the other hand, the difference is expected to approach zero
locally since the majority of which is a dispersive wave that
scatters to infinity. Hence one expects that it is {\it
asymptotically stable} in the sense that
\begin{equation}
\label{asymp.stab} \norm{ \psi(t) -Q_{0,n(t)} \, e^{i \theta(t)}
}_{L^2 \loc} \to 0
\end{equation}
as $t \to \infty$, for a suitable choice of $n(t)$ and $\theta(t)$.
Here $\norm{\cdot }_{L^2 \loc}$ denotes a local $L^2$ norm, to be
made precise in \eqref{Lploc.def}. One is also interested in how
fast \eqref{asymp.stab} converges and whether $n(t)$ has a limit.

The second question is the asymptotic problem of the solution when
$\psi(0)$ is small but not close to ground states. It is delicate
since nonlinear excited states stay there forever but are expected
to be unstable from physical intuition. Thus, a solution may stay
near an excited state for an {\it extremely long time} but then
moves on and approaches another excited state.

We now review the literature, assuming $\psi_0$ is small in $H^1
\cap L^1$.

If $- \Delta + V$ has only one bound state, i.e., with no excited
states, the asymptotic stability of ground states is proved in
\cite{SW1,SW2}, with convergence rate $t^{-3/2}$. It is then shown
in \cite{PW} that all solutions with small initial data, not
necessarily near ground states, will locally converge to a ground
state.

Suppose $- \Delta + V$ has two bound states. the asymptotic
stability of ground states is proved in \cite{TY1}, with a slower
convergence rate $t^{-1/2}$ due to the persistence of the excited
state. The problem becomes more delicate when the initial data are
away from ground states. It is proved in \cite{TY3} that, near
excited states, there is a finite co-dimensional manifold of
initial data so that the corresponding solutions locally converge
to excited states. Outside of a small wedge enclosing this
manifold, all solutions exit the excited state neighborhood and
relax to ground states \cite{TY2}. It is further showed in
\cite{TY4} that for all small initial data in $H^1 \cap L^1$,
there are exactly three types of asymptotic profiles: vacuum,
excited states or ground states. The last problem is also
considered in \cite{SW4}.

Suppose $- \Delta + V$ has three or more bound states. The
asymptotic stability of ground states is proved in \cite{Tsai}. In
fact, it is shown that all solutions with \begin{equation}
\norm{\psi_0}_{H^1\cap L^1}^{3-\e}\le |(\phi_0,\psi_0)|\ll 1, \quad
0<\e\ll 1, \end{equation} relax to ground states. It ensures that
the solution is away from excited states but allows the ground state
component to be much smaller than other components.

We also mention a few related results on the asymptotic stability of
ground states of nonlinear Schr\"odinger equations with more general
nonlinearities. For small solutions, one extension is to replace the
resonance condition \eqref{resonant-cond} by weaker conditions, e.g.
those by \cite{Gang-Sigal} and by \cite{Cuc-Miz}. Another extension
is to assume $\psi_0 \in H^1$ without assuming $\psi_0 \in L^1$. It
is first proved in \cite{GNT} for $K=0$ and dimension $N=3$ and then
extended by \cite{Miz1,Miz2} for $K=0$ and $N=1,2$. It is also
extended by
%\cite{GNT2} for $K=1$ and $N=3$, and then by
\cite{Cuc-Miz} for $K \ge 1$
%, $N=1,2,3$,
with \eqref{resonant-cond} replaced by weaker conditions used  by
\cite{Gang-Sigal}. A third extension is to allow subcritical
nonlinearity $\pm |\psi|^{p-1}\psi$, $p< 1+4/N$, see e.g.
\cite{Kirr-Miz}. A fourth extension is to assume $K=1$ and $e_1$ has
multiplicity, see \cite{Gang-Wei, Gus-Phan}.

The stability of {\it large} solitary waves is considered for
$K=0,1$,  by \cite{BP1, BP2,BS} for $N=1$ and   by  \cite{C1,C2} for
$N=3$.

See \cite{KriSch, Schlag, Gus-Phan} and their references for
construction of stable manifolds similar to that in \cite{TY3}.

\medskip

In this paper, our goal is to continue the study of \cite{Tsai}
under the same assumptions, with initial data $\psi_0$ now inside
a neighborhood of the first excited state $Q_{1,n}$. This is the
easiest interesting case not covered in \cite{Tsai}. Guided by the
$K=1$ case, one expects that the solution should either converge
to a first excited state (with the ground state component always
negligible), or leave the excited state neighborhood after some
time (which may be extremely long, say greater than
$e^{e^{-1/n}}$), and then relax to a ground state.

The new difficulty of the $K>1$ case is the existence of higher
excited state components. If the solution is to converge to a
first excited state with the ground state component always
negligible, one can think that the ground state component is
absent and the first excited state as a new ground state. Thus, in
the $K>1$ case the convergence to a first excited state is
expected to be in the rate $t^{-1/2}$, much slower than $t^{-3/2}$
in the $K=1$ case.

When the difference is of order  $t^{-3/2}$, one can use {\it
centered orthogonal coordinates} as in \cite{PW,TY4},
\begin{equation}\label{centeredortho}
\psi(t) = Q_{1,n(t)}e^{i \th(t)} +  h(t), \quad h(t) =x_0(t)\phi_0
+ \xi(t), \quad \xi \in \bE_c(H_0).
\end{equation}
The equations of $\dot n(t)$ and $\dot \th(t)$ contain linear
terms in $h$. When $x_0(t)$ is negligible, these linear terms are
of order $t^{-3/2}$ and hence integrable in $t$, ensuring the
convergence of the parameters. However, when $K>1$, the difference
is order $t^{-1/2}$ and one cannot show the convergence of the
parameters if their equations contain linear terms. To remove
linear terms, one is forced to use {\it linearized coordinates}
around the first excited state, to be specified later in
\S\ref{S3.2}.

We now describe a few special properties of the linearized
operator around an excited state. When the function $\psi$ is
close to a nonlinear bound state $Q=Q_{m,n}$ with corresponding
frequency $E=E_{m,n}$, one writes $\psi = (Q(x)+h(t,x))e^{-iEt}$.
The perturbation $h(t,x)$ satisfies
\begin{equation}
\pd_t h = \L h + \text{nonlinear terms},
\end{equation}
where the linearized operator $\L$ around $Q$ is given by
\begin{equation}
\L h = -i \bket{ (H+\kappa  Q^2)\,h + \kappa  Q^2 \,\wbar h\, },
\quad H=-\Delta + V -E +  \kappa  Q^2.
\end{equation}
Note $HQ=0$. Since $\L$ does not commute with $i$, it is not useful
to consider its spectral properties. Instead one looks at its matrix
version acting on $\smat{ \Re h \\ \Im h}$:
\begin{equation}
\bL =\mat{ 0 & H \\ -H -2 \kappa  Q^2 & 0}.
\end{equation}%
The spectral property of $\bL$ for $m >0$ is studied in \cite{TY3}
and recalled in Proposition \ref{L-spectral}. It is a perturbation
of $J(H_0-e_m)$ with $J= \mat{ 0 & 1
\\ -1 & 0}$ which has eigenvalues  $\pm i (e_k-e_m)$, $k=0,\ldots,K$.
When $m>0$, $k<m$ and $e_k<2e_m$, the eigenvalues $\pm i
(e_k-e_m)$ are embedded in the continuous spectrum $\pm i\,
[|e_m|,\I)$. These embedded eigenvalues split into a quadruple of
eigenvalues of $\bL$, $\pm \LLLL_k$ and $\pm \bar \LLLL_k$, with
$\Im \LLLL_k = |e_k-e_m|+O(n^2)$ and $C^{-1} n^4 < \Re \LLLL_k < C
n^4$ (assuming the generic condition \eqref{gamma0.def}). The size
of their corresponding eigenvectors are roughly\footnote{Denote
$\bka{\xi}=(1+|\xi|^2)^{1/2}$ for $\xi \in \R^d$, $d \ge 1$. For
$r\in \R$, denote by $L^2_r$ the weighted $L^2$ spaces with
$\norm{f}_{L^2_r} = \norm{\bka{x}^r f(x)}_{L^2}$.}
\begin{equation}
  O_{L^2_{100}}(1) + \frac{O(n^2)}{\bka{x}} 1_{|x|<n^{-4}}.
\end{equation}
The second part is not localized; It is small in $L^\I \cap L^3$,
of order 1 in $L^2$, and of order $n^{6-12/p}$ in $L^p$ for $p<2$.
In particular, the projection $P_c^\bL$ onto the continuous
spectral subspace $\bE_c^\bL$ of $\bL$ is of order $n^{6-12/p}\gg
1$ in $L^p$ for $p<2$, giving an extra difficulty to the usual
analysis.

To overcome this difficulty, we prove decay estimates of the form
(see Lemma \ref{th:decay})
\begin{equation}
\label{S1:decay} \norm{e^{t\bL}\Pcsp  \ph}_{L^p} \le C_p t^{-\frac
32+\frac 3p} \bka{t}^{\frac 3{2p}} \norm{\ph}_{L^{p'}}, \quad (t
\ge 0)
\end{equation}
for $3\le p <6$, with constant $C_p$ independent of $n$. Here
$\Pcsp$ is an extended projection: It is  the sum of $P_c^\bL$ and
all projections onto eigenspaces whose corresponding eigenvalues
have negative real parts. As shown in Remark (iii) after Lemma
\ref{th:decay}, these estimates with $n$-independent constant are
false if $\Pcsp$ is replaced by $P_c^\bL$. Also note that
\eqref{S1:decay} is time-direction sensitive: it is true only for
$t \ge 0$. The decay exponent above is not as good as the usual
free Schr\"odinger evolution, but it is sufficient for us if we
take $p<6$ close to $6$. A side benefit of extending $\Pc$ to
$\Pcsp$ is that we no longer need to track the component
$(\Pcsp-\Pc)h$.

\bigskip
Our assumptions on the operator $H_0= -\Delta + V$ are as follows:

\noindent {\bf Assumption A0}. $H_0= - \Delta + V$ acting on
$L^2(\R^3)$ has $K+1$ simple eigenvalues $e_0<e_1< \ldots< e_K<0$,
$K\ge 2$, with normalized real eigenvectors $\phi_0, \ldots,
\phi_K$.

\noindent {\bf Assumption A1}. $V(x)$ is a real-valued function
satisfying $|\nabla^\al V(x)| \lec \bka{x}^{-5- s_1}$ for $|\al|\le
3 $, for some $s_1>0$. $0$ is not an eigenvalue nor a resonance for
$H_0$.

\noindent {\bf Assumption A2}. Resonance condition. We assume that
\begin{equation}\label{resonant-cond2}
%  e_{k-1} < 2 \, e_k, \quad \forall k=1,\ldots, K .
e_0 < 2 \, e_1 < 4 \, e_2.
\end{equation}
%
%Hence $2e_{k0}>|e_0|$ for all $k>0$.
We further assume that, for some small $s_0>0$,
\begin{equation} \label{gamma0.def}
\gamma_0 \equiv
%\inf_{\myfrac{0 \le m < k, l \le K}{|s|<s_0}}
\inf_{\myfrac{0 \le m \le 1, |s|<s_0}{ m < k, l \le K}}
\lim_{r \to 0+} \Im \bke{ \phi_m\phi_k^2 \, , \, \frac
{1}{-\Delta+V+e_m - e_k - e_l -s -r i} \Pc ^{H_0} \phi_m\phi_k^2 }
> 0.
\end{equation}

\noindent {\bf Assumption A3}. No-resonance condition (between
eigenvalues). Let $j_{\text{max}}=3$. For all
$j=2,\ldots,j_{\text{max}}$ and for all $k_1,\ldots, k_j
,l_1,\ldots, l_j \in \bket{0, \ldots, K}$,
if $e_{k_1} +\cdots+ e_{k_j} = e_{l_1} +\cdots+ e_{l_j}$, then there
is a permutation $s$ of $\{1,\ldots,j\}$ such that $
(l_1,\ldots,l_j) = (k_{s1},\ldots,k_{sj})$.
\bigskip

%with $\bket{k_1,\ldots, k_j} \not = \bket{l_1,\ldots, l_j}$ as sets
%with multiplicities, (e.g., $\bket{0,0,1} \not  = \bket{0,1,1}$),
%
%\begin{equation}
%e_{k_1}+ \cdots + e_{k_j} \not = e_{l_1} + \cdots + e_{l_j} .
%\end{equation}

Assumption A1 ensure several estimates for linear Schr\"odinger
evolution such as decay estimates and the $W^{k,p}$ estimates for
the wave operator $W_{H_0}=\lim_{t\to \infty} e^{i t
H_0}e^{it\Delta}$. They are certainly not optimal. The main
assumption in A2 is the condition $e_{k-1} < 2 e_k$. It ensures that
$H_0 + e_m - e_k-e_l$ is not invertible in $L^2$ for $m <k,l$, and
provides (for our cubic nonlinearity) the required resonance between
eigenvalues through the continuous spectrum. Since the expression
for $\gamma_0$ is quadratic, it is non-negative and $\gamma_0 > 0$
holds generically. Assumption A3 is a condition to avoid direct
resonance between the eigenvalues. It is trivial if $K=0,1$. It
holds true generically and is often seen in dynamical systems of
ODE's. If we relax the assumption \eqref{resonant-cond2}, we may
need to increase $j_{\text{max}}$.
Fix $r_1>10$ large enough.  We denote by $L^p \loc$ the local $L^p$
spaces given by the norm
\begin{equation} \label{Lploc.def}
\norm{\phi}_{L^p \loc(\R^3)} \equiv \bket{\int_{\R^3} \bka{x}^{-p
r_1} |\phi(x)|^p d x}^{1/p}.
\end{equation}

Now we are ready to state our main theorem.

\begin{theorem}\label{mainthm}
Assume Assumptions A0--A3 and fix $0<\de \le \frac 1{10}$. There
are constants $C_0,C_1>0$, and small $n_0>0$ such that the
following hold. If $n=(\phi_1, \psi_0) \in (0,n_0)$ and
$\norm{\psi_0 - n \phi_1}_{H^1 \cap L^1} \le n^{1+\de}$, then the
solution $\psi(t)$ of \eqref{Sch} with $\psi(0)=\psi_0$ satisfies
\begin{equation}\label{eq:mainthm}
\limsup_{ t \to \infty} \norm{\psi(t) - Q_{m,n_+} e^{i
\th(t)}}_{L^2_{loc}} t^{1/2} \le C_0/n
\end{equation}
for $m=0$ or $m=1$, for some $n_+\in (C_1^{-1} n,C_1 n)$ and
some %continuous real-valued function $\th(t)$.
$\th(t)\in C([0,\I),\R)$.
\end{theorem}

%{\it Comments on Theorem \ref{mainthm}}:

In fact we have more detailed estimates of the solution for all
time, see Propositions \ref{th:S-4}, \ref{tc-pro},
\ref{outside-in}, \ref{ti.lemma}, and \ref{gs.con}. In particular,
if the initial data $\psi_0$ is placed in the neighborhood of an
excited state $Q_{m,n}$ with $m \ge 2$, even if $K >2$,
Propositions \ref{th:S-4}, \ref{tc-pro}, \ref{outside-in},
\ref{ti.lemma} show that the solution will either converge to
$Q_{m,n_+}$ for some $n_+$, or eventually exits the neighborhood,
stays away from bound states for a time interval of order between
$n^{-4}\log \frac 1n$ and $n^{-4-2\de}$, until it reaches the
neighborhood of another bound state $Q_{m',n'}$, $m'<m$. If
$m'=0$, then Proposition \ref{gs.con} shows that $\psi(t)$ will
converge to some $Q_{0,n_+}$. However, if $m'>0$, our current
analysis is not sufficient to control its evolution after this
time.

We now sketch the structure of our proof and this paper.

In \S2 we give the linear analysis, including the decay estimates
\eqref{S1:decay}.

In \S3 we consider the decomposition of the solutions in different
coordinates and the normal forms of their equations.

In \S4 we start with the solution in a $n^{1+\de}$-neighborhood of
$Q_{1,n}$ and use linearized coordinates. We follow the evolution as
long as the ground state component $z_0$ is negligible,
characterized by $|z_0(t)| < n^{-3}(n^{-4-2\de} + t)^{-1}$. If it is
always negligible, we prove that the solution converges to an
excited state with convergence rate $t^{-1/2}$.

In \S5 we consider the case that $|z_0(t_c)| \ge
n^{-3}(n^{-4-2\de} + t_c)^{-1}$ in a first time $t_c \in [0,\I)$,
which may be $0$ or extremely large, say $>e^{e^{-1/n}}$. After an
initial layer, we show that $|z_0(t)|$ starts to grow
exponentially with exponent $Cn^4$ until it reaches the size
$2n^{1+\de}$ at time $t_o$. The time it takes, $t_o-t_c$, is of
order $n^{-4} \log \frac {2n^{1+\de}}{|z_0(t_c)|}$. Along the way
higher excited states may have size larger than $|z_0(t)|$
%, \rho_c)$ where $\rho_c =n^{-1}(n^{-4-2\de} + t_c)^{-1/2}$,
but can be controlled. This section is the most difficult part in the
nonlinear analysis because it involves estimates not previously
studied.

In \S6  we study the dynamics after $t_o$ when there are at least
two components of size greater than $2n^{1+\de}$, and change to
{\it orthogonal coordinates}
\begin{equation} \label{1-6}
\psi = x_0 \phi_0 + \cdots + x_K \phi_K + \xi, \quad \xi \in
\bE_c(H_0).
\end{equation}
Although $\xi(t_o)$ is already non-localized, we can prove
``outgoing estimates'' for $\xi(t_o)$, introduced in
\cite{TY2,TY4}, to capture the time-direction sensitive
information of the dispersive waves. We show that, after a time of
order between $n^{-4}\log \frac 1n$ and $n^{-4-2\de}$, the ground
state component $x_0$ grows to order $n$ while all other
components become smaller than $n^{1+\de}$. (This is called the
{\it transition regime}.)

In \S7 the ground state component becomes dominant and we change
to linearized coordinates around it. Again we need to keep track
of out-going estimates during the coordinate change. We show that
the solutions will converge to ground states with convergence rate
$t^{-1/2}$. The analysis is similar to \S4 but easier because it
has no unstable direction. (This is called the {\it stabilization
regime}.)

Analysis similar to \S6 and \S7 is done in \cite{Tsai}, (and in
the two-eigenvalue case near ground states in
\cite{BP2,TY1,TY2,C2,BS}). However, with weaker decay estimates
like \eqref{S1:decay}, we need more refined analysis. For example,
since the nonlinearity is of constant order $n^3$  in the
transition regime, we need to make this time interval as short as
possible by taking $\de>0$ small. We also take $p<6$ close to $6$
to minimize our loss in estimating the $L^p$-norm of the
dispersive component during this interval.

\subsection*{New proof of linear decay estimates for ground states}

We end this introduction by noting that, our linear analysis, Lemmas
\ref{th:decay} and \ref{th:sdecay}, in the case $m=0$, provide a new
proof of linear estimates for the linearized operators around ground
states, which is used to prove the stability of ground states in 3D,
see \cite{C1, TY1, Tsai}.  Proofs in these references either use the
wave operator between $\L$ and $-i (H_0 - E)$, or use a similarity
transform $\L = U (-iA)U^{-1}$ for some self-adjoint perturbation
$A$ of $H_0 -E$ and non-self-adjoint operator $U$. Our proof here
use simple perturbation argument and requires less assumptions on
the potential $V$. Moreover, this perturbation argument allows the
operator $V$ to be more general than a potential, as long as the
decay and singular decay estimates for $-\De + V$ hold.

%%%%%%%%%%%%%%%%%%%%%%%%%%%%%%%%%%%%%%%%%%%%%%%%%%%%%%%%%%%%%%%%%%%%%%%%%%

\section{Linear analysis}

In this section we will study various properties of the linearized
operator around a fixed bound state, in particular an excited
state. The starting point is the following lemma on the existence of
nonlinear bound states and their basic properties.

\begin{lemma}[Nonlinear bound states]
\label{th:2-1} Assume Assumptions A0--A1. There exists a small
$n_1>0$ such that for each $k =0,\ldots,K$ and $n \in [0,n_1]$,
there is a solution $Q_{k,n} \in H^2 \cap W^{1,1}$ of \eqref{Q.eq}
with $E=E_{k,n}\in \R$ such that
\begin{equation}
Q_{k,n} = n\phi_k + q(n), \quad (q,\phi_k)=0.
\end{equation}
The pair $(q, E)$ is unique in the class $\norm{q}_{H^2}+|E-e_k|
\le n^2$.   Moreover, $ \norm{q}_{H^2 \cap W^{1,1}}\lec n^3$ and
$\norm{\frac {\pd}{\pd n}q}_{H^2 \cap W^{1,1}} + |E-e_k| \lec
n^2$.  $|E-e_k - Cn^2 | \lec n^4$ where $C=\kappa  \int \phi_k^4$.
We also denote $R_{k,n} = \frac{\pd}{ \pd E_k} Q_{k,n} =
\frac{\pd}{ \pd n} Q_{k,n}/ \frac{\pd}{ \pd n} E_{k,n} = \frac
1{2C n} \phi_k + O_{H^2 \cap W^{1,1}}(n)$.
\end{lemma}

In the following we fix $m \in \{0,\ldots,K\}$ and $n \in
[0,n_1]$. Let $Q = Q_{m,n}$, $R = R_{m,n}$ and $E=E_{m,n}$. The
function $Q$ satisfies $HQ=0$ where
\begin{equation}
\label{H.def} H= H_0 - E + \kappa  Q^2.
\end{equation}
The following lemma collects useful properties of $H$.

\begin{lemma} \label{th:H}
Assume Assumptions A0-A1 and let $H$ be defined as in
\eqref{H.def}. The operator $H$ has $K+1$ real eigenvalues $\td
e_k = e_k - e_m + O(n^2)$ with normalized eigenfunctions $\td
\phi_k = \phi_k + O(n^2)$. In particular, $\td e_m=0$ and $\td
\phi_m = C Q_m$. The projection to its continuous spectral
subspace is $P_c^H f = f -\sum_k (\td\phi_k,f)\td \phi_k$.
Furthermore, we have the following decay estimates
\begin{equation}
\label{H-decay} \norm{e^{-itH}\Pc^H \ph}_{L^q} \le C
|t|^{-3/2+3/q} \norm{\ph}_{L^{q'}}, \quad (2 \le q \le \infty),
\end{equation}
and singular decay estimates: for sufficiently large $r_1>9/2$,
for $0\le N\le 3$, for $\al_j \in \C$ with $\Im \al_j >0$, $|\Re
\al_j+e_m| \in [a_1,a_2] \subset (0,\infty)$, $j \le N$,
\begin{equation}
\label{H-sdecay} \norm{\bka{x}^{-r_1} e^{-itH}\Pi_{j=1}^{N} (H
-\al_j)^{-1} \Pc^H \ph}_{L^2} \le C \bka{t}^{-3/2}
\norm{\bka{x}^{r_1} \ph}_{L^2}, \quad (t \ge 0).
\end{equation}
Here the constant $C$ is independent of $n$, $\ph$ and $ \al_j$.

\end{lemma}

Note that this lemma contains $H=H_0$ as a special case with $n=0$.
The proof of the first part is well-known by perturbation. Estimate
\eqref{H-decay} is by Journe-Soffer-Sogge \cite{JSS}.  Estimate
\eqref{H-sdecay} for $N=0$ is by Jensen-Kato \cite{JK} and Rauch
\cite{Rauch}. Estimate \eqref{H-sdecay} for $\al_1 = \cdots =
\al_N$, $N \ge 1$, was first proven by Soffer-Weinstein \cite{SW3}
for Klein-Gordon equations, then by Tsai-Yau \cite{TY1} and Cuccagna
\cite{C2} for (linearized) Schr\"odinger equations. The general case
is similar and a proof based on Mourre estimate is sketched below
for completeness.  (See \cite{C2} for a different approach).

Denote the dilation operator $D=x\cdot p+p\cdot x$ with $p=-i\nabla$,
and the commutators \newcommand{\ad}{{\mathrm{ad}}}
\begin{equation}
\ad_D^0(H)=H, \quad \ad_D^{k+1}(H)=[\ad_D^k(H),D], \quad k \ge 0.
\end{equation}
Fix $g_* \in C^\infty_c(\R)$ with $g_*=1$ on $[-1,1]$ and $\supp g_*
\subset (-2,2)$. For each $j$, let $g_j(t) = g_*((t - \Re z_j)/\e)$.
If $\e>0$ is sufficiently small, $g_j(H) \ad^k_D(H) g_j(H)$ are
bounded operators in $L^2$ for $k \le 3$ and all $j$, and the Mourre
estimate holds: For some $\th>0$,
\begin{equation}
g_j(H) [iH,D]  g_j(H) \ge \th  g_j(H)^2, \quad \forall j.
\end{equation}
See \cite{CFKS}. Thus the pair $H,D$ satisfies the assumptions of the
minimal velocity estimates in \cite{HSS} and Theorem 2.4 of
\cite{Ski}, and one has
\begin{equation}
\norm{\chi(D \le \th t/2) e^{-itH} g_j(H) \bka{D}^{-r_1}}_{L^2 \to L^2}
\le C \bka{t}^{-r_1+\e_1},
\end{equation}
where $0<\e_1\ll 1$ and $\chi(D \le a)$ is the spectral projection
of $D$ associated to the interval $(-\infty,a]$.  The same argument
of \cite{SW3} then gives \eqref{H-sdecay}.

\subsection{Linearized operator}
A perturbation solution $\psi(x,t)$ of
\eqref{Sch} of the exact solution $Q(x)e^{-iEt}$ can be written in the
form
\begin{equation}
\psi(x,t)=[Q(x) + h(x,t)]e^{-iEt }
\end{equation}
for some function $h$ which is small in a suitable sense. Then, $h$
satisfies
\begin{equation}
\pd_t h = \L h + \ \text{nonlinear terms},
 \end{equation}
where the operator $\L$ is defined as
\begin{equation} \label{L.def}
\L h = -i\{(H_0-E + 2\kappa  Q^2) h + \kappa  Q^2 \bar{h}\}.
\end{equation}
The operator $\L$ is linear over $\R$ but not over $\C$. As a result
it is not useful to consider its spectral properties.

Consider the injection from scalar functions to vector functions
\begin{equation}
\bj: L^2(\R^3,\C) \to L^2(\R^3,\C^2), \quad \bj(\ph) = [\ph]: =
\smat{\Re \ph \\ \Im \ph}.
\end{equation}
With respect to this injection, the operator $\L$ is naturally
extended to a matrix operator acting on $L^2(\R^3, \C^2)$ with the
following form
\begin{equation}
\bL = \mat{ 0 & L_- \\ -L_+ & 0}, \quad \text{where}
\left \{ \begin{array}{ll}
L_- & = H = H_0 - E + \kappa  Q^2, \\
L_+ & = H+2\kappa  Q^2 = H_0 - E +3 \kappa  Q^2.
\end{array} \right.
\end{equation}
We will use $\L = \bj^{-1} \bL \bj$ for computations involving $\L$.

The space $L^2(\R^3, \C^2)$ is endowed with the natural inner product
\begin{equation}
(f,g) = \int_{\R^3} (\bar{f}_1 g_1 + \bar{f}_2g_2)\, dx
\end{equation}
for $f = \smat{ f_1 \\ f_2}$ and $g =\smat{ g_1 \\ g_2}$.  We will use
the Pauli matrices
\begin{equation}
\sigma_1 = \mat{ 0 &1 \\ 1 & 0} , \quad \sigma_2 = \mat{ 0 & -i \\ i &
0 } , \quad \sigma_3 = \mat{ 1 & 0 \\ 0 & -1 }.
\end{equation}

\subsection{Invariant subspaces}
In this subsection we study the spectral subspaces of $\bL$. Since
$\bL$ is a perturbation of $JH$, we first give the following lemma for
comparison.

\begin{lemma}[Invariant subspaces of $JH$]  \label{JH}
Assume Assumptions A0--A2.  The space \\ $L^2(\R^3,\C^2)$ can be
decomposed as the direct sum of $JH$-invariant subspaces
\begin{equation}
\label{L2decomp0} L^2(\R^3,\C^2)= \bE_0^{JH} \oplus \cdots \oplus
\bE_K^{JH} \oplus \bE_c^{JH}.
\end{equation}
For each $k\in \{0,\ldots,K\}$, the space $\bE_k^{JH}$ is spanned
by 2 eigenvectors $\smat{1\\-i}\td \phi_k$ and $\smat{1\\i}\td
\phi_k$ with eigenvalues $-i\td e_k$ and $i\td e_k$, respectively.
Its corresponding orthogonal
projection is $P_k^{JH}\smat{f_1\\ f_2} = \smat{ (\td \phi_k, f_1) \\
(\td \phi_k, f_2) }\td \phi_k$.  The subspace $\bE_{c}^{JH}$
%:=\{\smat{f_1\\ f_2} : (\td \phi_k, f_j)=0, \forall k,j\}$
has projection
$P_c^{JH} f
%= f - \sum_{k=0}^K P_k^{JH} f
=\smat{P_c^H f_1 \\ P_c^H f_2}$.
\end{lemma}

The proof is straightforward and skipped. We next give the
corresponding statements for $\bL$.

%==================================================================
\begin{proposition}[Invariant subspaces of $\bL$] \label{L-spectral}
Assume Assumptions A0--A2. Fix $m\in \{0,\ldots,K\}$ and $n \in
(0, n_1]$. Let $Q = Q_{m,n}$, $R = R_{m,n}$ and $E=E_{m,n}$. The
space $L^2(\R^3,\C^2)$ can be decomposed as the direct sum of
$\bL$-invariant subspaces
\begin{equation}
\label{L2decomp} L^2(\R^3,\C^2)= \bE_0^{\bL} \oplus \cdots \oplus
\bE_K^{\bL} \oplus \bE_c^{\bL}.
\end{equation}
If f and g belong to different subspaces, then
\begin{equation}
(\si_1 f, g)=0.
\end{equation}
These subspaces and their corresponding projections satisfy the
following.
\begin{enumerate}

\item[(i)] $\bE_m^\bL$ is the $0$-eigenspace spanned by $\smat{ 0\\ Q }$
and $\smat{ R\\ 0}$, with $\bL \smat{ 0\\ Q } = \smat{0 \\ 0}$ and
$\bL \smat{R\\ 0} = -\smat{0\\ Q }$. Its projection is $P_m f= c_m
(\si_1 \smat{ R\\ 0 },f) \smat{ 0\\ Q } + c_m (\si_1 \smat{ 0\\ Q
},f) \smat{ R\\ 0 }$, $c_m = (Q,R)^{-1}$.

\item[(ii)] $\bE_{k}^\bL$ for $0 \le k < m$, if such $k$ exists, is
spanned by 4 eigenvectors $\Phi_k= \smat{ u_k\\ -iv_k }$, $\bar
\Phi_k$, $\si_3 \Phi_k$ and $\si_3 \bar \Phi_k$, with eigenvalues
$\LLLL_k$, $\bar \LLLL_k$, $-\LLLL_k$, and $-\bar \LLLL_k$,
respectively. Here $\LLLL_k = -i(e_k - e_m) + O(n^2)$, $ n^4 \lec
\Re \LLLL_k \lec n^4$, $u_k$ and $v_k$ are complex-valued
functions, $u_k = \bar u_k^+ + \bar u_k^-$ and $v_k = \bar u_k^+ -
\bar u_k^-$ , with
\begin{equation}
\label{u-pm} u_k^+=  \phi_k + \olcsp(n^2), \quad u_k^-= (H - i
\bar \LLLL_k)^{-1} \phi_k^* + \olcsp(n^2)
\end{equation}
where $ \phi_k^* = \Pc^H \phi_k^* = \olcsp (n^2)$. Furthermore,
$(u_k,v_k)=0$ and $(u_k,v_\ell)=(\bar u_k,v_\ell)=0$ for $k \not =
\ell$. All $(\bar u_k,v_k)$, $\norm{u_k^+}_{L^2}$ and
$\norm{u_k^-}_{L^2}$ are equal to $1+O(n^2)$ and
$\norm{u_k^-}_{L^2_{loc}} \lec n^2$.  The projection to
$\bE_k^\bL$ is $P_k + \Psp_k $ where
\begin{equation}
\label{Pk-def}
\begin{split}
P_k f &= c_k (\si_1 \bar \Phi_k, f)\Phi_k + \bar c_k (\si_1 \Phi_k, f)
\bar \Phi_k, \\
\Psp_k  f &= -c_k (\si_1 \si_3 \bar \Phi_k, f)\si_3\Phi_k - \bar c_k
(\si_1 \si_3 \Phi_k, f) \si_3 \bar \Phi_k,
\end{split}
\end{equation}
and $c_k = (\si_1 \bar \Phi_k, \Phi_k)^{-1}=i/(\int 2 u_kv_k)= i/2
+ O(n^2)$.

\item[(iii)] $\bE_{k}^\bL$ for $m<k \le K$, if such $k$ exists, is
spanned by 2 eigenvectors $\Phi_k= \smat{ u_k\\ -iv_k }$ and $\bar
\Phi_k$ with eigenvalues $\LLLL_k$ and $\bar \LLLL_k$,
respectively. Here $\R \ni i \LLLL_k = e_k - e_m + O(n^2)$, $u_k$
and $v_k$ are real-valued, both equal to $\phi_k + \olcsp(n^2)$,
and normalized by $(u_k,v_k)=1$.  Its projection is $P_k$, also
given by \eqref{Pk-def}, with $c_k = i/2$.

\item [(iv)] $\bE_{c}^\bL=\{g: (\si_1 f,g)=0, \forall f \in \bE_k, \forall
k =0,\ldots,K\}$. Its projection is $P_c^\bL f= f - \sum_{k=0}^K
P_k f - \sum_{k<m}\Psp_k  f$.

\end{enumerate}
\end{proposition}
%==================================================================

Note that $\LLLL_k$ is in the first quadrant and near the
imaginary axis for $k <m$, and in the lower imaginary axis for
$k>m$. They are all perturbations of $-i\td e_k$ of Lemma
\ref{JH}. When $k<m$, $-i\td e_k$ are inside the continuous
spectrum $\pm i[|E_m|,\infty)$ and their resonance make the
eigenvalues split.

\medskip

\setlength{\unitlength}{1mm}\noindent
\begin{center}
\begin{picture}(120,80)

\put (0,40){\vector(1,0){120}}

\put (60,0){\vector(0,1){80}}

\put (60.5,60){\line(0,1){18}}\put (59.5,60){\line(0,1){18}}

\put (60.5,20){\line(0,-1){20}}\put (59.5,20){\line(0,-1){20}}

\put(60,40){\circle*{2}}
\put(64,42){\makebox(0,0)[c]{\scriptsize$\bE_m^\bL$}}

\put(60,53){\circle*{2}}
\put(69,53){\makebox(0,0)[c]{\scriptsize$\bar\Phi_{j>m}, \bar
\LLLL_j$}}

\put(60,27){\circle*{2}} \put(67,27){\makebox(0,0)[c]{\scriptsize$
\Phi_{j},  \LLLL_j$}}

\put(63,70){\circle*{2}}
\put(73,70){\makebox(0,0)[c]{\scriptsize$\Phi_{k<m}, \LLLL_k$}}

\put(57,70){\circle*{2}}
\put(47,70){\makebox(0,0)[c]{\scriptsize$\si_3 \bar \Phi_{k},
-\bar \LLLL_k$}}

\put(63,10){\circle*{2}}
\put(71,10){\makebox(0,0)[c]{\scriptsize$\bar \Phi_{k},
\bar\LLLL_k$}}

\put(57,10){\circle*{2}}
\put(47,10){\makebox(0,0)[c]{\scriptsize$\si_3 \Phi_{k}, -
\LLLL_k$}}

\end{picture}

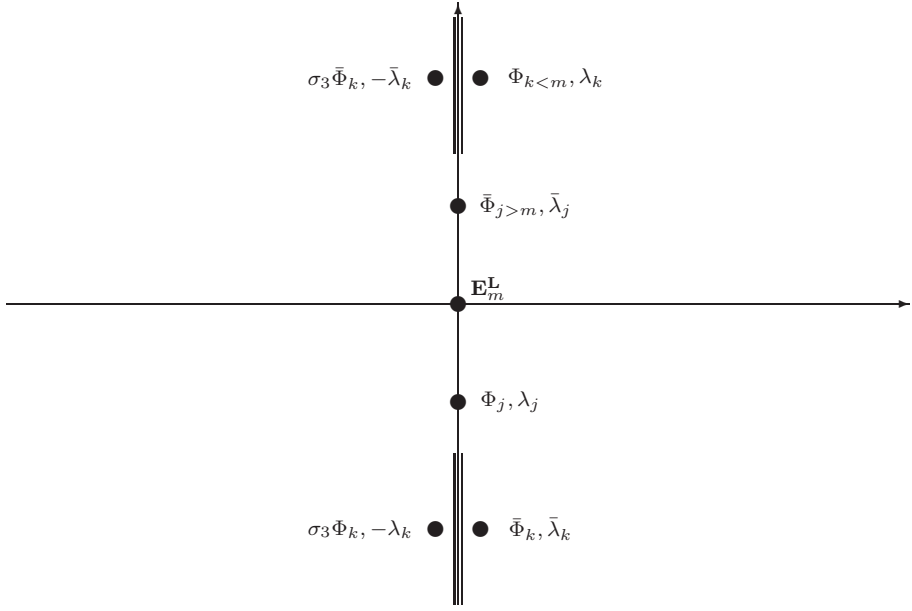
Figure 1: Spectrum of $\bL$ around $Q_m$, $0<m<K$.
\end{center}
\medskip

\myproof The same proof of \cite[Theorem 2.2]{TY3} works in our
many eigenvalue case. The only thing we need to check is the
properties of $u_k^+$ and $u_k^-$ when $k<m$.  Fix $k<m$. Denote
by $\Pi$ the orthogonal projection from $L^2$ onto $\{ \td \phi_k,
Q_m\}^\perp $, and $B= 2 \kappa  Q_m^2$. We omit the subscript $k$
below. By the defining equations $\L_m \Phi = \LLLL  \Phi$ and
$\Phi = \smat{u \\ -iv}$, $\bar u$ satisfies
\begin{equation}
(H^2 + H B)\bar u = -\bar \LLLL ^2 \bar u.
\end{equation}
By the same proof for the two-eigenvalue case in \cite[section
  2.1]{TY3} (in which $\Pi=P_c^H$), $\bar u$ can be solved in the form
\begin{equation}
\bar u = \tilde \phi + h, \quad h=\Pi h=- (H^2 + \Pi H B \Pi +
\bar \LLLL ^2)^{-1}  \Pi H B \tilde \phi.
\end{equation}
One can rewrite
\begin{equation}
h = (H^2+\bar \LLLL ^2)^{-1} \Psi, \quad \Psi = \Pi \Psi=[1+\Pi
HB\Pi (H^2+\bar \LLLL ^2)^{-1}]^{-1} \Pi H B \tilde \phi.
\end{equation}
By resolvent estimates and a power series expansion as in
\cite{TY3}, the function $\Psi$ is localized and
$\norm{\Psi}_{\lcsp } \le Cn^2$. Since $v =(i\LLLL )^{-1}(H+B)u$,
we have $u^\pm = \mp \frac 1{2z}(H\mp z +B) \bar u$ with $z= i\bar
\LLLL  = |e_k-e_m| + O(n^2)$. For $u^+$,
\begin{equation}
u^+ = -\frac 1{2z}(H- z) \tilde \phi -\frac 1{2z}(H+ z)^{-1} \Psi
-\frac 1{2z}B \bar u.
\end{equation}
The first term is equal to $(1+O(n^2)) \tilde \phi$. Since $(H+
z)^{-1}\Pi$ is order one, the remaining two terms are $\olcsp(n^2)$,
and so is $\phi - \tilde \phi$. This shows $u^+ = \phi +
\olcsp(n^2)$. For $u^-$,
\begin{equation}
u^- = \frac 1{2z}(H+ z) \tilde \phi +\frac 1{2z}(H- z)^{-1}  \Psi +
\frac 1{2z}B \bar u.
\end{equation}
The first term is $O(n^2) \tilde \phi$. Since $(H-z)^{-1}(\Pi -
\Pc^H)\Psi$ are sum of eigenfunctions with $O(n^2)$ coefficients, we
get \eqref{u-pm} with $ \phi_k^* = \frac 1{2z} \Pc^H \Psi =
\olcsp(n^2)$.

The orthogonality $(u,v)=0$ is equivalent to $(\si_1 \Phi,\Phi)=(\si_1
\si_3 \Phi,\Phi)=0$, which follow from the general fact shown in
\cite[\S 2.6]{TY3} that
\begin{equation}\label{eq2-11}
(\si_1 f,g)=0 \quad \text{if} \quad \bL f=\LLLL f, \ \bL g = \mu
g, \text{ and }\bar \LLLL  \not = \mu.
\end{equation}
It also follows from \eqref{eq2-11} that $(u_k,v_\ell)=(\bar
u_k,v_\ell)=0$ for $k \not = \ell$.  That $\norm{u^+}_{L^2}=1+O(n^2)$
and $\norm{u^-}_{L^2_{loc}} \lec n^2$ follow from \eqref{u-pm}.  Note
\begin{equation}
0 = (\bar u,\bar v) = (u^+ + u^-, u^+-u^-) = (u^+,u^+) - (u^-,u^-) +
(u^-,u^+)-(u^+,u^-).
\end{equation}
Since the last two terms are $O(n^2)$, we get $\norm{u^-}_{L^2}-
\norm{u^+}_{L^2}=O(n^2) $. Finally
\begin{equation}
(\bar u,v)=(u^+ + u^-, \bar u^+-\bar u^-) = (u^+,\bar u^+) - (u^-,\bar
u^-) + (u^-,\bar u^+)-(u^+,\bar u^-) .
\end{equation}
We have $(u^-,u^+)-(u^+,u^-)=O(n^2)$. By \eqref{u-pm} we also have
(denoting $ok=O(n^4)$)
\begin{equation}
\label{eq-228} (\bar u^-, u^-) = ((H-\bar z)^{-1} \bar
\phi^*_k,(H-z)^{-1} \phi^*_k)+ok = (\bar \phi^*_k,(H-z)^{-2}
\phi^*_k) +ok= ok
\end{equation}
by the singular decay estimate of Lemma \ref{th:H} with $t=0$. Thus
$(\bar u_k,v_k)=1 + O(n^2)$.  Similarly, $(\bar u_k,v_\ell)= O(n^2)$
for $k \not = \ell$.
\myendproof

\medskip

In the following lemma we provide more properties of $u_k^-$.
\begin{lemma}
\label{uv-est} Assume the same as in Proposition \ref{L-spectral}
and fix $k<m$. Then

(i) $\norm{u_k^-}_{L^p} \le C_p (n^2+n^{6- \frac {12}p})$ for $1 \le p
\le \infty$, in particular $\norm{u_k^-}_{L^2_{-r}} \le C n^2$.

(ii) $\norm{e^{-isH}\Pc^H u_k^-}_{L^2_{-r}} +
\norm{e^{-isH_0}\Pc^{H_0} u_k^-}_{L^2_{-r}} \le C n^2\bka{s}^{-3/2}$
for $s \ge 0$.

(iii) $\norm{u_k^-}_{H^1} \leq C$.

\end{lemma}

\myproof Denote $z = i \bar \LLLL_k$ and $\ph = \phi_k^*$. For
(i), it suffices to check $(H-z)^{-1} \ph$, the main part of
$u_k^-$ in \eqref{u-pm}. Write $H-z = -\De + \nu^2 + V_1$ where
$V_1 = V + \kappa  Q_m^2$, $\nu^2 = E_m + z$ with $\Im \nu>0$.
Thus $\Im \nu \sim +n^4$.  By resolvent expansion,
\begin{equation}
(H-z)^{-1} \ph = (-\De - \nu^2)^{-1}
\ph +  (-\De - \nu^2)^{-1} V_1 (H-z)^{-1} \ph.
\end{equation}
Since the resolvent $(-\De - \nu^2 )^{-1}$ has the convolution kernel
$G(x)=(4\pi|x|)^{-1}\exp(i \nu |x|)$,
\begin{equation}
\norm{(-\De - \nu^2 )^{-1} \ph}_{L^p} \lec
\norm{G*\ph}_{L^p} \lec (\norm{G}_{L^p(B_1^c)} +
\norm{G}_{L^2(B_1)})\cdot \norm{\ph}_{L^1\cap L^2}
\end{equation}
which is bounded by $(n^{4-12/p}+1)\cdot n^2$.  Since $\norm{V_1
(H-z)^{-1} \ph }_{L^1\cap L^2} \lec \norm{(H-z)^{-1} \ph
}_{L^2_{-r}}\lec n^2$, we have the same bound for the second term. The
above show (i).

For (ii), we only need to consider $e^{-isH_0}\Pc^{H_0}u_k^-$
since the other term follows from Lemma \ref{th:H}. By resolvent
expansion $R=(H- z)^{-1}= R_0 (1+\kappa Q_m^2 R)$ where $R_0=(H_0-
E_m-z)^{-1}$,
\begin{equation}
\Pc^{H_0} u^- =R_0 \ph' + \olcsp(n^2), \quad \ph' = \Pc^{H_0}
(1+\kappa  Q_m^2 R) \ph = \olcsp(n^2).
\end{equation}
Thus
\begin{equation}
e^{-isH_0}\Pc^{H_0} u^- = e^{-isH_0}R_0 \ph' + O_{L^2_{-r}}(n^2
\bka{s}^{-3/2}).
\end{equation}
By the singular decay estimate for $H_0$, the first term is also of
order $O_{L^2_{-r}}(n^2 \bka{s}^{-3/2})$.

To prove (iii), it suffices to prove that $\norm{\nabla v}_{L^2} = O(1)$
where $v = (H -z)^{-1}\varphi$. It can be shown by
multiplying the equation $(H - z) v = \varphi$ by $\bar{v}$
and then integrating it on $\mathbb{R}^3$.
\myendproof

\medskip

We will need the following lemmas for scalar functions.

\begin{lemma}
\label{uv-orth} Fix $0 \le k \le K$, $k \not = m$. Let $\ph\in
L^2(\R^3,\C)$ be a scalar function.

(i) $P_k [\ph] = \Re \al \Phi_k$, $\bj^{-1}P_k[\ph] = \al \bar u^+
+ \bar \al u^-$, where
\begin{equation}
\al= 2 c_k (\si_1 \bar \Phi_k, [\ph]) =-2c_k i[( u_k^+ , \ph) -(
u_k^-, \bar \ph)].
\end{equation}

(ii) $P_k \ph =0$ iff $(\si_1 \Phi_k, [\ph])=0$ iff $(u_k^+, \ph) =
(u_k^-, \bar \ph)$.

(iii) For $k<m$, $\Psp_k  \ph =0$ iff $(\si_1 \si_3\Phi_k, [\ph])=0$
iff $(u_k^+, \bar \ph) = (u_k^-, \ph)$.
\end{lemma}

\myproof Write $[\ph]=\smat{\ph_1\\ \ph_2}$. Since $[\ph]$ is real, we
have by \eqref{Pk-def} that $P_k [\ph] = \Re \al \Phi_k$ with $\al = 2
c_k (\si_1 \bar \Phi_k, [\ph])$. Omitting the subscript $k$, we have
\begin{align*}
(\si_1 \bar \Phi_k, [\ph]) & = (i\bar v,\ph_1) + (\bar u,\ph_2) = (
u^+ - u^-, -i \ph_1) + ( u^+ + u^-, \ph_2) = -i( u^+ , \ph) + i( u^-,
\bar \ph),
\end{align*}
which gives the formula for $\al$. Thus
\begin{equation}
\bj^{-1} P_k [\ph]= \bj^{-1}\Re \al \mat{  u
\\ -i  v}
=\frac12 \bket{(\al u + \bar \al \bar u) + i(-i \al v + i \bar \al
\bar v)} = \al \bar u^+ + \bar \al u^-.
\end{equation}
The claim (ii) follows from (i). For (iii), since $\si_3 \si_1
\si_3 = - \si_1$, $(\si_1 \si_3\Phi_k, [\ph])=0$ is equivalent to
$0=(\si_1 \Phi_k, \si_3 [\ph])=(\si_1 \Phi_k, [\bar \ph])$ and
hence to $(u_k^+, \bar \ph) = (u_k^-, \ph)$. \myendproof

%When $k > m$,  $u_k$ and $v_k$ are real and claim (iii) is redundant.

\medskip

The following lemma will be used to treat the linear term in the
$\eta$ equation.
\begin{lemma}
\label{JPk} (i) For $k <m$,
\begin{equation}
\label{JPk-1a} J\Phi_k  =  i \Phi_k - 2i \smat{1 \\ -i} \bar
u_k^+.
\end{equation}

(ii) If $f\in L^2(\R^3,\C^2)$ and $P_k f=0$, then $\norm{P_k
Jf}_{L^2}\lec \norm{f}_{L^2_{-r}}$.
\end{lemma}

\myproof  For (i), rewrite
\begin{equation}
\label{eq2-34} \Phi_k = \smat{u_k \\ -iv_k} = \smat{1 \\ -i} \bar
u_k^+ + \smat{1 \\ i} \bar  u_k^-.
\end{equation}
Applying $J$
\begin{equation}
\label{JPk-1aa} J\Phi_k =-i \smat{1 \\ -i} \bar u_k^+ +  i\smat{1 \\
i} \bar  u_k^- .
\end{equation}
Canceling $u_k^-$ we get \eqref{JPk-1a}.

For (ii), we have $(\si_1 \bar \Phi_k,f)=(\si_1 \Phi_k,f)=0$.
Using $J^* = -J$, $J\si_1=-\si_1 J$, and \eqref{JPk-1a},
\begin{equation}
\label{JPk-2} (\si_1 \bar \Phi_k, J f) = -(J\si_1 \bar \Phi_k,  f)
= (\si_1 J \bar \Phi_k,f) = (\si_1 (-i \bar \Phi_k + 2i \smat{1 \\
i} u_k^+),f)= (2i\smat{i \\ 1}u_k^+,f).
\end{equation}
Similarly $(\si_1 \Phi_k, J f) = (2i\smat{i \\
-1}\bar u_k^+,f)$. This shows (ii). \myendproof

Note, in deriving \eqref{JPk-1a} if we cancel $u_k^+$ instead of
$u_k^-$, we get
\begin{equation}
\label{JPk-1b} J\Phi_k  = - i \Phi_k + 2i \smat{1 \\ i} \bar
u_k^-.
\end{equation}

%==================================================================

\subsection{Decay estimate}

In the following two subsections we prove decay estimates for
$e^{t \bL}$ with the constant independent of $n$. This
independence is essential for our analysis of the nonlinear
dynamics both inside a neighborhood of an excited and away from
bound states. For example, it ensures that the time spent
traveling between bound states is no longer than $O(n^{-4-2\de})$.

This independence cannot be achieved if we restrict ourselves of
$\bE_c$, the continuous spectral subspace, because the projection
$\Pc^\bL$ as an operator acting on $L^1$ is of order $O(n^{-6})$
due to the presence of $u_k^-$. Suppose $F$ is the total
nonlinearity in the equation of the perturbation $h$. Our choice
of parameters $a(t)$ and $\th(t)$ makes $P_m F=0$, but does not
make $F \in \bE_c$.  To avoid the large constant problem, we
absorb the range of $\Psp_k $, $k<m$, which have exponential
decay, into $\bE_c$. The range of $P_k$ for $k<m$, which have
exponential growth, is left out and will be taken care of using
the evolution with correct time direction.

\medskip

Define $\bE_c^\sharp $ as the direct sum of
$\bE_c^{\bL}$ and eigenspaces whose eigenvalues have {\it negative}
real parts
\begin{equation}
\bE_c^\sharp  = \bE_c^{\bL} \oplus \mathrm{span}_{\C} \{ \si_3 \Phi_k,\si_3
\bar \Phi_k: 0\le k <m \}.
\end{equation}
Its corresponding projection is denoted as
\begin{equation} \label{Pc+.def}
\Pcsp  f = \Pc^{\bL} f + \sum_{k<m} \Psp_k  (f) = f -P_d f, \quad
P_d f= \sum_{k=0}^K P_k (f).
\end{equation}

We extend the definition of $\Pcsp $ to scalar functions by $\Pcsp
\ph = \bj^{-1}\Pcsp  [\ph]$, and similarly for $P_d$.  If a scalar
function $\ph$ satisfies $[\ph] \in \bE_c^\sharp $, then $(\si_1
\Phi_k,[\ph])=0$ for all $k$.

%% Lemma \ref{eta-wteta} was placed here, but no longer used and moved to the end.

The next lemma is on the uniform bound of $H^1$-norm of $e^{t \bf L}
\Pcsp  \ph$ for $t \ge 0$.
\begin{lemma}
\label{H1conserve} For any scalar function $\ph\in H^1$ we have
\begin{equation}
\norm{ e^{t \bf L} \Pcsp \ph}_{H^1} \le C \norm{ \ph}_{H^1}, \quad
(t \ge 0),
\end{equation}
where the constant $C$ is independent of $n$ and $t \ge 0$.
\end{lemma}
\myproof  From \eqref{Pc+.def} and \eqref{Pk-def}, we have
\begin{equation}
e^{t \bf L} \Pcsp \ph  = e^{t \bf L} \Pc^{\bL} \ph
- \sum_{k<m} \Big [\bar{c}_k (\sigma_1\sigma_3\Phi_k, \ph)
e^{-\bar \LLLL_k t} \sigma_3\bar{\Phi}_k  + c_k
(\sigma_1\sigma_3\bar{\Phi}_k, \ph) e^{- \LLLL_k t} \sigma_3
\Phi_k \Big ].
\end{equation}
By Lemma \ref{uv-est}, we have $\norm{\Phi_k}_{H^1} = O(1)$ for
all $ k <m$. From this and $\Re \LLLL_k >0$ for all $k <m$, we can
find a constant $C>0$ independent of $n$ such that
\begin{equation}  \label{conser1}
 \norm{ e^{t \bf L} \Pcsp \ph}_{H^1}  \leq \norm{e^{t \bf L} \Pc^{\bL} \ph}_{H^1} + C\norm{\ph}_{H^1}.
\end{equation}
Moreover, by following the proof of \cite[(2.6)]{TY3}, we see that
there exists a constant $C$ independent of $n$ such that
\begin{equation} \label{conser2}
 \norm{e^{t \bf L} \Pc^{\bL} \ph}_{H^1} \leq C \norm{\Pc^{\bL} \ph}_{H^1}.
\end{equation}
Again, since $\norm{\Phi_k}_{H^1} = O(1)$ for all $k$, we
also have $\norm{\Pc^{\bL} \ph}_{H^1} \leq C \norm{\ph}_{H^1}$ for
some constant $C$ which is independent of $n$. From this,
\eqref{conser1}, and \eqref{conser2}, Lemma \ref{H1conserve} follows.
\myendproof

\begin{lemma}
\label{Ecp} If a scalar function $\eta$ satisfies $[\eta] \in
\bE_c^\sharp $, then
\begin{equation}
\norm{ \eta - \Pc^H \eta}_{L^\infty_{3r}} \lec n^2 \norm{\eta}_{L^2_{loc}}+
\sum_{k<m} |(\bar u_k^-,\Pc^H \eta)|.
\end{equation}
\end{lemma}

\myproof
Write $\eta'=\Pc^H \eta$ and
\begin{equation}
\eta - \eta' = (1-\Pc^H)\eta = \tsum_k (\td \phi_k,\eta)\td \phi_k.
\end{equation}
For $k \ge m$, $|(\td \phi_k,\eta)|\le ok$ where $ok$ denotes $ O(n^2
\norm{\eta}_{L^2_{loc}})$.  For $k < m$, by Lemma \ref{uv-orth} (ii),
\begin{equation}
(\td \phi_k,\eta) + ok = ( u_k^+,\eta)= ( u_k^-,\bar \eta)  = (
u_k^-,\bar \eta') + ( u_k^-,\bar \eta-\bar \eta').
\end{equation}
Since $\norm{u_k^-}_{L^2_{loc}} \lec n^2$,
\begin{equation}
( u_k^-,\bar \eta-\bar \eta') = \tsum_{j=0}^K ( u_k^-,(\td
\phi_j,\eta)\td \phi_j) = ok.
\end{equation}
The above show the lemma.
\myendproof

The following lemma provides decay estimates for $e^{-itH}u_j^-$.
\begin{lemma}
\label{th2-12} Let $H_*$ be the self-adjoint realization of $
-\De$ on $L^2(\R^3)$. Let $V$ be a localized real potential so
that $H_* + V$ satisfies the decay and singular decay estimates
\eqref{H-decay} and \eqref{H-sdecay}. Let $0<n<n_0 \ll 1$, $a >0$,
and $z=a + n^4 i$. Let $\ph(t) =  n^2 (H_* + V-z)^{-1}e^{-it(H_* +
V)} \Pc g$ with $\norm{g}_{L^1}\lec 1$ and $\Pc=\Pc^{H_*+V}$. Then
for all $p \in (3,\infty]$, $m=\frac 12 - \frac 3{2p}\in [0,1/2]$,
\begin{equation}
\label{eq:2.44} \|\ph(t)\|_{L^p} \lec
t^{-m}(1+t)^{-m-\min(m,1/4)}, \quad \forall t>0.
\end{equation}
Above the $p$-dependent constant is uniform in $a \in [a_1,a_2]
\subset (0,\infty)$ and independent of $t$ and $n$.
\end{lemma}

\myproof The case $V=0$ is postponed to Subsection
\ref{freedecay}. For general $\td H = H_* + V$, by resolvent
expansion and Duhamel's formula,
\begin{align*}
\ph(t)&= n^2(H_*-z)^{-1} e^{-itH_*} \Pc g +n^2(H_*-z)^{-1}\int_0^t
e^{-i(t-s)H_*}Ve^{-is\td H} \Pc g
\\
&\quad + n^2(H_*-  z)^{-1}V (\td H-z)^{-1} e^{-it\td H} \Pc  g.
\end{align*}
Denote $\td\al_p(t):=t^{-m}(1+t)^{-m-\min(m,1/4)}$.
 By the estimate for $H_*$, \eqref{H-decay}, \eqref{H-sdecay},
and the proof for Lemma \ref{uv-est} (i), the above is bounded by
\begin{equation}
\norm{\ph(t)}_{L^p}\lec \td \al_p(t) + \int_0^t \td\al_p(t-s)
\bka{s}^{-3/2} ds + n^2 \cdot \bka{t}^{-3/2} \lec \td \al_p(t).
\end{equation}
\myendproof

The following is the main result of this subsection.

\begin{lemma}[Decay estimate]
\label{th:decay} For any scalar function $\ph \in L^{9/8} \cap
L^{3/2}$,
\begin{equation}
\label{th:decay-eq1} \norm{e^{t\bL}\Pcsp  [\ph]}_{L^\infty+ L^2}
\le C\al_\I(t) \norm{\ph}_{L^{9/8} \cap L^{3/2}}, \quad (t \ge 0).
\end{equation}
For $3< p <6$ and any scalar function $\ph\in L^{p'}$,
\begin{equation}
\label{th:decay-eq2} \norm{e^{t\bL}\Pcsp  [\ph]}_{L^p} \le C_p
\al_p(t) \norm{\ph}_{L^{p'}}, \quad (t \ge 0).
\end{equation}
Above the constants are independent of $n$ and $\ph$, and
\begin{equation}
\label{th:decay-eq3} \quad \al_\I(t) := t^{-1/2}\bka{t}^{-2/3},
\quad \al_p(t) :=t^{-\frac 32 + \frac 3p} \bka{t}^{\frac 3{2p}} .
\end{equation}
\end{lemma}
\myremark (i) For \eqref{th:decay-eq1} we could have chosen $\ph
\in L^{q} \cap L^{3/2}$, $\frac{12}{11} \le q < \frac65$. Then
$\al_\I(t) = t^{-1/2}\bka{t}^{-s}$, with $s=3/q -2 \in (1/2, 3/4]$
by the same proof. The exponent $q=\frac{12}{11}$ gives the
optimal decay rate that Lemma \ref{th2-12} provides for
$e^{-itH}P_cu_j^-$. However, when we estimate $\norm{\eta^3}_{L^q}
\lec \norm{\eta}_{L^2}^{3-3\th} \norm{\eta}_{L^p}^{3\th}$, we
prefer a larger $q$. For convenience we choose $q=9/8$.

(ii) Suppose we keep $q=\frac{12}{11}$ with $\al_\I(t)=t^{-1/2}
\bka{t}^{-3/4}$, and estimate $\norm{\eta^3}_{L^{12/11}} \lec
\norm{\eta}_{L^2}^{3-3\th} \norm{\eta}_{L^p}^{3\th}\lec \al_\I(t)$,
we need $\frac {11}2 < p <6$.

(iii) These estimates are false if $\Pcsp$ is replaced by $\Pc$.
Suppose the contrary, then they would be also true if $\Pcsp$ is
replaced by $P_d^\sharp = \Pcsp - \Pc$. Consider the case $m=1$ and
$\ph = \phi_0$ the $e_0$-eigenfunction of $-\De+V$. Then
\begin{equation}
\norm{e^{t\bL}P_d^\sharp   [\ph]}_{L^p} \sim e^{-c n^4 t}, \quad
\norm{\ph}_{L^{p'}} \sim 1.
\end{equation}
However the former is not bounded by $C t^{-k}$ for all $t>0$, for
any $k>0$ and $C$ independent of $n$.

\medskip

\myproof Denote $\eta(t) = e^{t\bL}\Pcsp [\ph]$ and $\eta' =
\Pc^{JH} \eta$. Lemma \ref{Ecp} implies
\begin{equation}
\norm{\eta}_X \lec \norm{\eta'}_X + {\textstyle \sum_{k<m}} |(\bar
u_k^-,\eta')|,\quad X=L^\infty+L^2.
\end{equation}
Denote $\bL = JH + W_1$ with $W_1=\smat{0 & 0 \\ - 2 \kappa  Q_m^2
& 0}$. By Duhamel's formula,
\begin{equation}
\eta'(t) = e^{tJH} \Pc^{JH}  \Pcsp  [\ph] + \int_0^t \Pc^{JH}
 e^{(t-s)JH} W_1 \eta(s)\,ds.
\end{equation}

By Lemma \ref{uv-orth} (i),
\begin{equation}
 \bj^{-1} \Pcsp  [\ph] =  \ph - \bj^{-1} \Re \tsum_{j=0}^K
z_j\Phi_j = \ph - \tsum_{j=0}^K (z_j \bar u_j^+ +  \bar z_j u_j^-)
\end{equation}
where $z_j\in \C$ are bounded by $\norm{\ph}_{L^q}$ for any $q\le
2$. Using \eqref{u-pm} for $j<m$ in particular $u_j^- = (H- i \bar
\LLLL_j)^{-1} \phi_j^* + \olcsp(n^2)$, $\Im i \bar \LLLL_j\sim
n^4$, and by Lemma \ref{th:H} (with $p=9,3$) and Lemma
\ref{th2-12} (with $p=\I$),
\begin{equation}\label{etap-est}
\norm{\eta'(t)}_X \lec  \al(t) \norm{\ph}_{Y} + \int_0^t
\bka{t-s}^{-3/2} n^2 \norm{\eta(s)}_X \, ds.
\end{equation}
where $\al(t) =t^{-1/2}\bka{t}^{-2/3}$ and $Y= L^{9/8} \cap
L^{3/2}$. By the same reasons,
\begin{equation}
|(\bar u_k^-,\eta')| = (\bar \phi_k^*, (H-i \bar \LLLL_k)^{-1}
\eta')+ O(n^2
 \norm{\eta'}_X),
\end{equation}
and
\begin{equation}
| (\bar \phi_k^*, (H-i \bar \LLLL_k)^{-1} \eta')|\lec n^2
\norm{(H-i \bar \LLLL_k)^{-1} \eta'}_{L^2_{loc}}\lec n^2 \cdot
\text{ RHS
  of } \eqref{etap-est}.
\end{equation}
Summing the estimates, we get $\norm{\eta(t)}_X \lec  \text{ RHS
  of } \eqref{etap-est}$, which implies \eqref{th:decay-eq1}.

The estimate \eqref{th:decay-eq2} is proved similarly with $X=L^p$,
$Y=L^{p'}$  and $\al(t) = \al_p(t)\sim \max(\td \al_p(t),
t^{-3(\frac 12 - \frac 1p)})$. \myendproof

\subsection{Decay estimate for free evolution with resonant data}
\label{freedecay}

In this subsection we prove Lemma \ref{th2-12} for $H_* = -\De$,
i.e. decay estimate for $\ph(t) = n^2 (H_* - z)^{-1}e^{-tH_*}g$
where $z=a+n^4 i$, $a\sim 1$, and $g\in L^1$. The operator
$(H_*-z)^{-1} e^{-itH_*}$ has symbol $(\xi^2 -
z)^{-1}e^{-it\xi^2}$ and thus its Green's function $G$ is radial
and, for $r=|x|$,
\begin{align*}
G(r,t)& = (2\pi)^{-3} \int_0^\infty (p^2 - z)^{-1}e^{-itp^2}
\int_{|\om|=1}
  e^{i pr \om_1} dS(\om) \,p^2dp
\\
& = (2\pi)^{-3} \int_0^\infty (p^2 - z)^{-1}e^{-itp^2} 4\pi
\frac{\sin
  (rp)}{rp}\,p^2dp
\\
&= \frac 1{4\pi^2 i r} \int _\R  (p^2 -
z)^{-1}e^{-itp^2}e^{irp}\,pdp.
\end{align*}

It is well known that $G(r,0) = \frac {1}{4\pi r}e^{i \sqrt z r}$.
We are not aware of an explicit formula for $G(r,t)$. Because for
$3< p\le \infty$ we have
\begin{equation}
\norm{\ph(t)}_{L^p} = \norm{n^2G(t)*g}_{L^p}\lec
n^2\norm{G(t)}_{L^{p}}\norm{g}_{L^1} \lec
n^2\norm{G(t)}_{L^{3,\I}}^{3/p}\norm{G(t)}_{L^{\I}}^{1-3/p}\norm{g}_{L^1}
,
\end{equation}
estimate \eqref{eq:2.44} follows from \eqref{Green-fcn-decay} of
the following lemma.

\begin{lemma}
\label{Green-fcn} Let $H_*$ be the self-adjoint realization of $
-\De$ on $L^2(\R^3)$. Let  $G(x,t)$ be the Green's function of the
operator $(H_*-z)^{-1} e^{-itH_*}$ where $z$ is the same as in
Lemma \ref{th2-12}. Then $G(x,t)=G(|x|,t)$ and
\begin{equation}
|G(r,t)| \lec \left \{
\begin{aligned}
\label{eq20-5} &  \frac{ r^{-1/2}}{n^4r + r^{1/2} + (t-r)_+},
\quad & r>1,\frac t{100},
\\
& t^{-3/2}, \quad & 1<r<\frac t{100},
\\
& \min(t^{-1/2}(1+t)^{-1},\ r^{-1}) , \quad & r<1.
\end{aligned}\right .
\end{equation}
In particular,
\begin{equation}
\label{Green-fcn-decay}
 \|G(\cdot ,t)\|_{L^\infty_x}\lec t^{-1/2}(1+t)^{-1/2}(1+n^4
 t^{1/2})^{-1},
 \quad \|G(\cdot ,t)\|_{L^{3,\infty}_x}\lec 1.
\end{equation}
\end{lemma}

\myproof We may assume $a=1/4$. The general case follows from change
of variables and is uniform for $a \in [a_1,a_2]$. Introduce a
regularizing factor $e^{-\de p^2}$ and write $(p^2-z)^{-1}$ as a
time integral (using $\Re z>0$)
\begin{align}
\nn G(r,t) & = \lim_{\de \to 0_+} \frac 1{4\pi^2  r} \int _\R
\int_0^\infty e^{-itp^2-\de p^2 - is(p^2-z)+irp}ds\,pdp
\\
& = \lim_{\de \to 0_+} \frac 1{4\pi^2  r} \int_0^\infty e^{isz+
\frac {ir^2}{4\al}} \int _\R e^{-i\al(p- \frac r{2\al})^2 }
\,pdp\,ds  , \quad \al = s+t -i\de.
\end{align}
Using $\int_\R e^{-p^2}dp = \sqrt \pi$ and
\begin{equation}
\int_\R e^{-i\al (p-\be)^2}p dp = \int_\R e^{-i\al (p-\be)^2}\be dp
= \be \int_\R e^{-i\al p^2} dp = \be  (i\al)^{-1/2}\sqrt \pi,
\end{equation}
we get
\begin{align}
G(r,t) &= \lim_{\de \to 0_+} \frac 1{4\pi^2  r} \int_0^\infty
e^{isz+ \frac {ir^2}{4\al}} \frac r{2\al} (i\al)^{-1/2}\sqrt \pi
\,ds \nonumber
\\
&=\frac 1{8\pi^{3/2}\sqrt i} \int_0^\infty e^{isz+ \frac
{ir^2}{4(s+t)}} (s+t)^{-3/2}\,ds \nonumber
\\
&=\frac 1{8\pi^{3/2}\sqrt i} \int_t^\infty e^{i\Phi} s^{-3/2}\,ds,
\label{Grt.formula}
\end{align}
where the phase $\Phi$ is
\begin{equation}
\Phi(r,s) = sz -tz + \frac {r^2}{4s}, \quad \Phi_s = z - \frac
{r^2}{4s^2},\quad \Phi_{ss} =  \frac {r^2}{2s^3}.
\end{equation}
Note $z=\frac 14 + n^4 i$, $\Phi_s$ vanishes at $s=r/(2\sqrt z)\sim
r$, and $\Re i \Phi <0$ for $s>t$.

First note
\begin{equation}\label{eq.2.81}
|G(r,t)| \lec   \int_t^\infty  s^{-3/2}\,ds = Ct^{-1/2},
\end{equation}
which is valid for all $r>0$ and $t>0$. We will use a stationary
phase argument to get a better estimate. The main contribution
should come from $I\equiv r(1-\mu,1+\mu)$ where $0<\mu \le \frac
1{200}$ will be chosen. Comparing \eqref{eq.2.81} and \eqref{eq20-1}
below, it is clear we do not get a better estimate unless $\mu$ is
small.

We first consider the case $r>1$.

Suppose $t \in I$. The contribution from $s\in (t,r+\mu r)$ is
bounded by
\begin{equation}
\label{eq20-1} |\int _{t}^{r+\mu r} e^{i\Phi}s^{-3/2}ds| \lec
\int_I r^{-3/2}ds \lec \mu r^{-1/2}.
\end{equation}

The contribution from $(r+\mu r,\infty)$ is, with $t_1= r+\mu r$,
\begin{equation}
\label{eq2.83} \int _{t_1}^\infty e^{i\Phi} s^{-3/2}ds = \int
_{t_1}^\infty \pd_s (e^{i\Phi}) \frac 1{i\Phi_s} s^{-3/2}ds =
\frac 1{i\Phi_s} e^{i \Phi} s^{-3/2}|_{s=t_1} + \int_{t_1}^\infty
e^{i \Phi} J ds,
\end{equation}
where
\begin{equation}
J= -\frac \pd{\pd s}(\frac 1{i\Phi_s} s^{-3/2}) = \frac
{\Phi_{ss}}{i(\Phi_s)^2} s^{-3/2} + \frac {3}{2i\Phi_s s^{5/2}}.
\end{equation}
For $s \ge t_1$, we have $|\Phi_s|\sim n^4 + (s-r)/r$ and
$|\Phi_{ss}| \lec s^{-1}$. Thus $|J|\lec
(|\Phi_s|^{-1}+|\Phi_s|^{-2})s^{-5/2}$, and the boundary term is
bounded by
\begin{equation}
\label{eq20-2} |\frac 1{i\Phi_s} e^{i \Phi} s^{-3/2}|_{s=t_1}|\lec
\frac 1{|\Phi_s(t_1)|} t_1^{-3/2}\lec \frac {r^{-3/2}}{n^4+\mu}.
\end{equation}
Decompose $(t_1,\infty) = (t_1, 100 r) \cup (100r,\infty)$. On
$(t_1,100r)$, we have
\begin{equation}
\label{eq20-3} |\int_{t_1}^{100r} e^{i \Phi} J ds| \lec
\int_{t_1}^{100r} \frac{ r^{2-5/2}}{(n^4r + s-r)^2} ds \lec \frac{
r^{-1/2}}{ n^4r + t_1 -r} = \frac {r^{-3/2}}{n^4+\mu} .
\end{equation}
For $s> 100 r$, we have $|\Phi_s| \gec 1$ and
\begin{equation}
\label{eq20-4} |\int_{100r}^\infty e^{i \Phi} J ds| \lec
\int_{100r}^\infty s^{-5/2}
  ds \lec r^{-3/2}.
\end{equation}
We now choose $\mu\le \frac 1{200}$ so that $\mu r^{-1/2} \sim \frac
{r^{-3/2}}{n^4+\mu}$. If $r\ge 1$, we can choose $\mu = \frac 1{200}
r^{-1/2}(1+n^8 r)^{-1/2}$ and get for $t/r \in (1-\mu,1+\mu)$
\begin{equation}
\label{eq20-4-1} |G(r,t)| \le \frac {r^{-1/2}}{n^4 r+ r^{1/2}}.
\end{equation}

If $t\in (r+\mu r,100r)$, we can take $t_1 = t$ in the above
estimates and ignore the contribution from \eqref{eq20-1} to get the
bound for $r>1$
\begin{equation}
\label{eq20-4a}|G(r,t)| \lec \frac {r^{-1/2}}{n^4r +|t-r|} .
\end{equation}

If $t>100r$, we can replace $100r$ by $t$ in \eqref{eq20-4}  and
ignore the contribution from \eqref{eq20-1} and \eqref{eq20-3} to
get (also true for $r<1$),
\begin{equation}
\label{eq20-4b} |G(r,t)| \lec t^{-3/2}.
\end{equation}

If $t\in ( \frac r{100}, r-\mu r)$ and $r>1$, the additional
contribution from $s\in (t,r-\mu r)$ is estimated as in
\eqref{eq2.83}--\eqref{eq20-3} with $t_1=r-\mu r$ and $100r$
replaced by $r/100$, and bounded by \eqref{eq20-4a}, which is
smaller than \eqref{eq20-4-1} for $r>1$.

If $t\in (0, \frac r{100})$, we have $|\Phi_s|\sim r^2 s^{-2}$ and
$|\Phi_{ss}|\sim r^2 s^{-3}$ for $s\in (t,\frac r{100})$. The
additional contribution from $s\in (t,\frac r{100})$ is estimated as
in \eqref{eq2.83}--\eqref{eq20-3} and bounded by
\begin{equation}
\label{eq20-4c}\bkt{  r^{-2} s^{1/2}}_{s=t}^{r/100} + \int
_{s=t}^{r/100} r^{-2} s^{-1/2}ds \le r^{-3/2}
\end{equation}
which is smaller than \eqref{eq20-4-1} for $r>1$.

We now consider the case $r<1$. Let $\al>0$ be a small number to be
chosen. The contribution from $s \ge \max(t,\al r)$ is bounded by
\begin{equation}%\label{}
|\int_{\al r}^\infty e^{i \Phi}s^{-3/2}1_{s>t}\,ds| \le |\int_{\al
r}^\infty s^{-3/2} ds|  =C (\al r)^{-1/2}.
\end{equation}
If $t<\al r$, we have $|\Phi_s|^{-1} \sim r^{-2} s^2$,
$|\Phi_{ss}|/|\Phi_s| \lec s^{-1}$,  and the contribution from $s <
\al r$ is
\begin{equation}%\label{}
\int_t^{\al r} e^{i \Phi}s^{-3/2}\,ds =\bkt{\frac 1{i\Phi_s} e^{i
\Phi} s^{-3/2}}_{s=t}^{s=\al r} + \int_t^{\al r} e^{i \Phi} J ds,
\end{equation}
which is bounded by
\begin{equation}%\label{}
%\al^2 (t)^{-3/2}.\quad
r^{-2} (\al r)^{1/2}.
\end{equation}
We want $(\al r)^{-1/2}\sim r^{-2} (\al r)^{1/2}$ and we can choose
$\al = \frac r{100}$, which gives $r^{-1}$ bound for $r<1$.

In conclusion, we have proved \eqref{eq20-5} for all $r>0$ and
$t>0$. \myendproof

\medskip

\myremark (i) Lemma \ref{th2-12} for the free case can be considered
an estimate of $(f, n^2 G(t)g)$. If \eqref{eq20-5} cannot be
improved, then Lemma \ref{th2-12} cannot be improved, even if
assuming further that one of $f,g$ is in $L^2_r$ (but not both). To
see it, let $g$ be the characteristic function of the unit ball.
Note $|I|\sim \mu r \gg 1$ for $r\gg 1$, thus $(n^2 G g)(r,t)$ has
the optimal size at $r\sim t$.  Since translation does not change
the $L^1\cap L^2$-norm of $f$, we can put the support of $f$ at $r
\sim t$, showing the optimality of Lemma \ref{th2-12}.

 (ii) Although the real
part of the phase, $e^{-n^4(s-t)}$, is decaying, it does not seem to
improve our estimate. In the case $t \sim r \sim n^{-8}$, we have
$|I|\sim \mu r \sim n^{-4}$ and the estimate \eqref{eq20-1} does not
improve because of the factor $e^{-n^4(s-t)}$, in view of the
identity $ \int_0^{n^{-4}} e^{-n^4 s} ds = C \int_0^{n^{-4}} ds $.

(iii) Since $|\Im \Phi_s| \sim |s-r|/r \lec \mu$ for $s \in I$,
$e^{i \Phi}$ almost has no oscillation on $I$ if $\mu^2 r \sim \mu
\cdot |I| \ll 1$. Thus, if $\mu = \e r^{-1/2}$ with $0<\e \ll 1$,
then the upper bound in \eqref{eq20-1} is also a lower bound.  In
the case $t\sim r \gg \e^{-2} n^{-8}$, we have $\mu \ll n^4$ and
$\mu r^{-1/2} \gg \frac {r^{-3/2}}{n^4 } \sim \frac {r^{-3/2}}{n^4 +
\mu}$. Thus \eqref{eq20-5} is optimal in this case.

\subsection{Singular decay estimate}

We will need to identify the main part of
\begin{equation}
\eta(t) = \int_0^t e^{(t-s)\L} \Pcsp  e^{-i\al s}f(s)ds
\end{equation}
where $\al \in \C$ with $\Im \al >0$ and $f(s)$ is an $L^2$-valued
function of $s$ with $\dot f$ smaller than $f$ in a suitable sense. We
will rewrite it in matrix form in order to integrate by parts. Using
\begin{equation} \label{scalar-vector}
[\ph] = \smat{\Re \ph \\ \Im \ph} = \Re\ph \smat{1\\ -i},
\end{equation}
and denoting $R=(\bL+i\al)^{-1}$, we have
\begin{align}
\eta(t) &= \bj^{-1}\Pc^{\sharp} \int_0^t  e^{(t-s)\bL} \Re
e^{-i\al
    s} f(s) \smat{1\\ -i} ds \nonumber
\\
&= \bj^{-1}\Pc^{\sharp} \Re \bigg ( -R  e^{-i\al
    t} f(t) \smat{1 \\ -i}  \nonumber
\\
&\hspace{28mm} + e^{t\bL}R f(0) \smat{1 \\ -i}
+ \int_0^t  e^{(t-s)\bL} R  e^{-i\al  s} \dot f(s) \smat{1\\ -i} ds \bigg).
\end{align}
To estimate the last two terms, we need the following lemma.

\begin{lemma}[Singular decay estimate]
\label{th:sdecay} There is a constant $C>0$ independent of $\al\in
\C$ with $\Im \al>0$, $n \in [0,n_0]$, and vector function $\Psi
\in \lcsp $ so that
\begin{equation}
\norm{\bj^{-1} \Re e^{t \bL}(\bL + i \al)^{-1}\Pcsp
\Psi}_{L^2_{loc}} \le C \bka{t}^{-3/2} \norm{\Psi}_{\lcsp }, \quad
(t \ge 0).
\end{equation}
\end{lemma}

\myproof Denote by $\eta$ the scalar function to be estimated,
$\eta(t)= \bj^{-1} \Re e^{t \bL}R \Pcsp  \Psi$, and $\eta'=\Pc^H
\eta$. Lemma \ref{Ecp} implies
\begin{equation}
\label{sdecay-eq1} \norm{\eta}_{L^2_{loc}} \lec
\norm{\eta'}_{L^2_{loc}} + {\textstyle \sum_{k<m}} |(\bar
u_k^-,\eta')|.
\end{equation}
Denote $\bL = JH + W_1$ with $W_1=\smat{0 & 0 \\ - 2 \kappa  Q_m^2
& 0}$, $R=(\bL +i \al)^{-1}$ and $R_0= (JH+ i\al)^{-1}$. By
Duhamel's formula and resolvent expansion,
\begin{equation}
%\label{etap-exp}
\eta'(t)= \Pc^{H}\bj^{-1} \Re \bke{ {e^{t JH}} R_0(1+ W_1 R)\Pcsp
\Psi + \int_0^t e^{(t-s)JH} W_1 \eta(s) \,ds}.
\end{equation}

Denote the first term on the right side by $\eta'_1(t)$.  Using $
\Pcsp  \Psi = \Psi - \tsum _{k} P_k(\Psi)$,
\begin{equation}
\eta'_1(t) = \bj^{-1} \Re {e^{t JH}} R_0\Pc^{JH} \bke{\Psi - \tsum
_{k<m} P_k(\Psi) + \Psi_1},
\end{equation}
where $\Psi_1=\Pc^{JH}[- \sum_{k \ge m} P_k\Psi + W_1R \Pcsp  \Psi ]$
is localized with
\begin{equation}
\label{Psi1-est} \norm{\Psi_1}_{L^2_r} \lec n^2 \norm{\Psi}_{L^2}
+ n^2 \norm{ R
  \Pcsp  \Psi}_{L^2_{loc}}\lec n^2 \norm{\Psi}_{L^2_r}.
\end{equation}

Note that
\begin{equation}
e^{tJH} = \mat{ \cos(tH) & \sin(tH) \\ - \sin(tH) & \cos(tH) }
=\sum_{\e=\pm 1} e^{i\e t H} \tfrac 12(I - i\e J),
\end{equation}
\begin{equation}
(JH+i\al)^{-1} = (H^2-\al^2)^{-1} (-JH + i \al),
\end{equation}
and
\begin{equation}
(I - i\e J) (-JH + i \al) = -\e i(H-\e \al) ( I -\e i J).
\end{equation}
We conclude, for $R_0=(JH+i\al)^{-1}$,
\begin{equation}
\label{etJHR0} e^{tJH}R_0 = \sum_{\e=\pm 1} e^{i\e t H}
(H+\e\al)^{-1} \tfrac {-\e i}2(I - \e i J).
\end{equation}

By \eqref{etJHR0}, \eqref{Psi1-est}, Lemma \ref{th:H}, and $\Im \al
>0$,
\begin{equation}
\norm{ \bj^{-1} \Re {e^{t JH}} R_0\Pc^{JH} \bke{\Psi +
\Psi_1}}_{L^2_{loc}} \lec \bka{t}^{-3/2} \norm{\Psi}_{L^2_r} .
\end{equation}

For $k<m$, note
\begin{equation}
(I+ iJ) \Phi_k =  2 \bar u_k^+ \smat{1 \\ -i}, \quad
(I+ iJ) \bar \Phi_k =  2 u_k^- \smat{1 \\ -i}.
\end{equation}
Using \eqref{etJHR0} and writing $P_k \Psi = a \Phi_k + b \bar
\Phi_k$, we have
\begin{align*}
 \Re  {e^{t JH}} R_0P_k \Psi
&=   \Re \tsum_{\e = \pm 1}  e^{i\e t H} (H+\e\al)^{-1} \tfrac {-\e
i}2(I - \e i J)P_k \Psi
\\
& =   \Re   e^{-i t H} \bket{
(H-\al)^{-1} i ( a \bar u^+  + b u^-)
+ (H+\bar \al)^{-1} i ( \bar b \bar u^+  + \bar  u^-) }\smat{1 \\  -i} .
\end{align*}
By \eqref{scalar-vector},
\begin{equation}
\bj^{-1} \Re {e^{t JH}} R_0P_k \Psi = e^{-i t H} \bket{ (H-\al)^{-1}
i ( a \bar u^+ + b u^-) + (H+\bar \al)^{-1} i ( \bar b \bar u^+ +
\bar a u^-) }.
\end{equation}
Note $\Im (- \bar \al) = \Im \al >0$. By Lemma \ref{th:H} and \eqref{u-pm},
\begin{equation}
\norm{ \Pc^H \bj^{-1} \Re {e^{t JH}} R_0P_k \Psi }_{L^2_{loc}} \lec
n^2 \bka{t}^{-3/2}\norm{\Psi}_{L^2}.
\end{equation}
Thus
\begin{equation}\label{etap-est2}
\norm{\eta'(t)}_{L^2_{loc}} \lec \bka{t}^{-3/2} \norm{\Psi}_{L^2_r} +
\int_0^t \bka{t-s}^{-3/2} n^2 \norm{\eta(s)}_{L^2_{loc}} \,ds.
\end{equation}
On the other hand, for $j<m$, by \eqref{u-pm} again,
\begin{equation}
|(\bar u_j^-,\eta')| = (\bar \phi_j^*, (H-i \bar \LLLL_j)^{-1}
\eta')+ O(n^2
 \norm{\eta'}_{L^2_{loc}}).
\end{equation}
Note $\Im i \bar \LLLL_j>0$. By Lemma \ref{th:H} and the previous
decomposition of $\eta'$,
\begin{equation}
| (\bar \phi_j^*, (H-i \bar \LLLL_j)^{-1} \eta'(t))|\lec n^2
\norm{(H-i \bar \LLLL_j)^{-1} \eta'}_{L^2_{loc}}\lec n^2 \cdot
\text{ RHS
  of } \eqref{etap-est2}.
\end{equation}
By \eqref{sdecay-eq1} and summing the estimates, we get
$\norm{\eta(t)}_{L^2_{loc}} \lec \text{ RHS of } \eqref{etap-est2}$,
which implies the lemma.
\myendproof
%===================================================================================================

\subsection{Upper and lower spectral projections}
\label{upperlower}

In this subsection we prove various estimates for the spectral
projections $\Pi_\pm$ which are defined in \eqref{Pipm.def} and
corresponds to $\pm \Im z \ge |E|$ in the spectrum of $\bL$. In
particular, Lemma \ref{Pipmest} allows us to replace $\Pcsp$ by
$P_{\pm} = \Pcsp\Pi_\pm $ in Lemmas \ref{th:decay} and
\ref{th:sdecay}.

Decompose $\bL = JA + W_2 = JH + W_1$ where $A=-\De + |E|$, $W_2 =
J(V+\kappa  Q^2)+W_1$, and $W_1 = \mat{ 0 &  0 \\ -  2\kappa  Q^2
& 0}$. Let $R(z)=(\bL-z)^{-1}$, $R_0(z)=(JA-z)^{-1}$ and
$R_1(z)=(JH-z)^{-1}$ be their resolvents. Note $R_0(z)$ can be
decomposed as
\begin{equation}
\label{R0.dec}
\begin{split} R_0(z) &= (JA -z)^{-1}= \mat{-z & A \\
-A & -z}^{-1} = (A^2+z^2)^{-1} \mat{-z & -A \\ A & -z}
\\
& = (A-iz)^{-1} M + (A+iz)^{-1} \bar M, \quad M = \frac 12\mat{i &
-1 \\ 1 & i}.
\end{split}
\end{equation}
$R_1(z)$ has a similar formula with $A$ replaced by $H$.

Let $\Ga_{c\pm}$ be contours about the upper and lower continuous
spectra $\Si_\pm= \pm [|E|i,+\I i)$, respectively. For an
eigenvalue $\LLLL $ of $\bL$, let $\Ga_\LLLL $ be a small circle
centered at $\LLLL $ with radius $\sim n^{4}$. All contours are
oriented clockwise and do not intersect. Let $P_*=\frac 1{2\pi
i}\int _{\Ga_*} R(z)dz$, $*=c\pm, \LLLL $, be their corresponding
spectral projections. Decompose $\Pcsp$ as the sum of its upper
and lower half plane components:
\begin{equation}
\label{Pcsp.dec} \Pcsp = P_+ + P_-, \quad P_\pm =
P_{c\pm}+P_{L\pm}, \quad P_{L+} =\sum_{k<m}P_{-\bar \LLLL_k},
\quad P_{L-}=\sum_{k<m}P_{- \LLLL_k}.
\end{equation}
 Also denote
\begin{equation}
\label{Pipm.def} \Pi_\pm = P_\pm + P_{R\pm}, \quad
P_{R+}=\sum_{k<m}P_{\LLLL_k}, \quad P_{R-}=\sum_{k<m}P_{\bar
\LLLL_k}.
\end{equation}
Note $P_\pm = \Pcsp \Pi_\pm$.

Let \begin{equation} \label{tau0.def} \de_0 = \frac 14 \min
\{|e_K|, |e_k - e_{k-1}|  : 1 \le k \le K\}, \quad \tau_0 = \frac
12 e_K -e_m.
\end{equation}
Note $\Im \LLLL_N < \tau_0 - \de_0 < \tau_0+\de_0 < |E|$.

We collect a few estimates for $R_0(z)$ and $R(z)$.
\begin{lemma}
\label{resolvent-est} Let  $\si_{d}^0 = \{ \pm i(e_k-e_m): 0\le k
\le K\}$, $s>\tfrac 12$ and $1\le p <\I$. We have
\begin{equation}
\label{R0.est}
\begin{split}
&\norm{R_0(z)}_{L^2_s \to L^2_{-s}} \le C \bka{z}^{-1/2}, \quad
 z \not \in i\R, \\
&\norm{R_1(z)}_{L^2_s \to L^2_{-s}} + \norm{R(z)}_{L^2_s \to
L^2_{-s}} \le C \bka{z}^{-1/2}, \quad z \not \in i\R,\,
\dist(z,\si_d^0)\ge \de_0,
\\
&\norm{R(z)}_{L^2_s \to L^2_{-s}} \le C n^{-4}, \quad 0<|\Re z| <
\frac 14\ga_0 n^4,\, \dist(z,\si_d^0)<\de_0,
\\
&\norm{R_0(z)}_{L^p \to L^p} +  \norm{R_1(z)}_{L^p \to L^p} +
\norm{R(z)}_{L^p \to L^p} \le C_p \bka{z}^{-1+\e_p}, \quad  |\Im z|
=  \tau_0.
\end{split}
\end{equation}
Above $\e_p=0$ for $p>1$ and $0<\e_1 \ll 1$, and the constants are
uniform in $n\in [0,n_0]$.
\end{lemma}

\myproof The first estimate is by the scalar case proved in
\cite[Remark 2 in Appendix A]{Agmon} and by \eqref{R0.dec}. The
second estimate is valid if $R(z)$ is replaced by $R_1(z) =
(JH-z)^{-1}$, which is by the scalar case proved in \cite[Theorem
9.2]{JK} and  by \eqref{R0.dec} with $A$ replaced by $H$. It is true
for $R(z)$ using the resolvent series $R(z)=R_1(z) \sum_{k=0}^\I
[W_1 R_1(z)]^k$ and the fact $W_1$ is a small localized matrix
potential. The third estimate is proved in \cite[Lemma 2.5]{TY3}.

The last estimate for $R_0(z)$ is by the scalar case proved in
\cite[Lemma 7.4]{C2} and by \eqref{R0.dec}. It is true for $R_1(z)$
because $\norm{(H-z)^{-1}}_{L^p\to L^p} \lec \bka{z}^{-1+\e_p} $ for
$|\Im z| = \tau_0$, which follows from
\begin{equation}
%\label{R0.est}
\begin{split}
(H-z)^{-1} f &= (H-z)^{-1} P_c f +  (H-z)^{-1} \tsum_{k=0}^K (\td
\phi_k,f)\td \phi_k
\\
&= W^{-1}(A-z)^{-1} W P_c f + \tsum_{k=0}^K (\td e_k - z)^{-1} (\td
\phi_k,f)\td \phi_k,
\end{split}
\end{equation}
where $W$ is the wave operator between $H$ and $A$ and $\td \phi_k$
are normalized eigenfunctions of $H$ with eigenvalues $\td e_k$.
Finally, the estimate for $R(z)$ follows from the resolvent series
$R(z)=R_1(z) \sum_{k=0}^\I [W_1 R_1(z)]^k$ again. \myendproof

\begin{lemma}
\label{Kpmest} Let $K_\pm = \Pi_\pm (J \mp i) $, initially defined
from $L^2_s $ to $L^2_{-s}$, $s>1$. For any $1 \le p \le q < \I$,
there is a constant $c$ so that $\norm{K_\pm u}_p \le c
\norm{u}_q$ for any $u \in L^2_s \cap L^q$.
\end{lemma}

This is clear for the reference self-adjoint operator $JA$, for
which $K_{\pm}=0$.

\myproof Recall $R_0$ is decomposed in \eqref{R0.dec}, and  $MJ = -i
M$ and $\bar M J = i \bar M$. As $z$ approaches $\Sigma_+ = [ |E|i,
+\I i)$, the upper continuous spectrum of $A$, the resolvent
$(A+iz)^{-1}$ is unbounded, and we write
\begin{equation}\label{RJiR}
R_0(z)J - i R_0(z) = -2i M(A-iz)^{-1} , \quad (z \sim \Sigma_+).
\end{equation}
Note right side is bounded. Similarly, as $z$ approaches $\Sigma_- =
-\Sigma_+$, we write
\begin{equation} R_0(z)J + i R_0(z) = 2i\bar M(A+iz)^{-1} , \quad (z
\sim \Sigma_-).
\end{equation}

We now prove the bound for $K_+$. The case of $K_-$ is similar.
Let $\Ga = \Ga_{c+} \cup \Ga_{p}$ and $\Ga_{p}=\cup_{k<m}
(\Ga_{\LLLL_k} \cup \Ga_{-\bar\LLLL_k})$. By spectral projection
formula and resolvent expansion,
\begin{equation}
\Pi_+ =\frac 1{2 \pi i}\int_\Ga R(z) dz = \frac 1{2 \pi i}\int_\Ga
[1+ R_0(z) W_0 + R_0(z) W_0 R(z) W_0]R_0(z) dz .
\end{equation}
By \eqref{RJiR},
\begin{equation}
\Pi_+( J - i) = \frac{-1}{\pi}\int_\Ga [ 1+ R_0(z) W_0 + R_0(z) W_0
R(z) W_0]M (A-iz)^{-1}dz = K_0 + K_1 + K_2.
\end{equation}
The above sum is well-defined as operators from $L^2_s$ to
$L^2_{-s}$.

Note that $K_0$ is zero since $(A-iz)^{-1}$ is regular inside $\Ga$
and the rest of the integrand of $K_0$ does not depend on $z$.

For $K_1$, the integral over $\Ga_{c+}$ is bounded from $L^q$ to
$L^p$ by Lemma 7.6 of Cuccagna [C2] using Coifman-Meyer multi-linear
estimates. The integral over $\Ga_{p}$ is also bounded from $L^q$ to
$L^p$ since
\begin{equation}\label{K1est}
\begin{split}
&\int _{\Ga_{p}}\norm{R_0(z) W_0 M(A-iz)^{-1}}_{L^q \to L^p} |dz|\\
& \le \int _{\Ga_{p}}\norm{R_0(z)}_{L^p \to L^p}
\norm{(A-iz)^{-1}}_{L^q \to L^q}   \lec \int _{\Ga_{p}} n^{-4} \cdot
1 \lec 1.
\end{split}
\end{equation}

For $K_2$, the integrand is analytic in $z$ and has enough decay in
$B(L^2_s \to L^2_{-s})$ in $|z|$ by Lemma \ref{resolvent-est}. Thus
we can change the contour to $\Ga_1 = \R + \tau_0 i$,
%We have
%\begin{equation}
%K_2 = \frac{-1}{\pi}\int_{\Ga_1} R_0(z) W R(z) W M (A-iz)^{-1}dz .
%\end{equation}
By Lemma \ref{resolvent-est}, $\norm{K_2}_{L^q \to L^p}$ is bounded
by
\begin{equation}
\int_{\Ga_1} \norm{R_0(z )}_{L^p \to L^p} \cdot \norm{R(z )}_{L^q
\to L^q}\cdot \norm{R_0(z )}_{L^q \to L^q} \, |dz| \le C.
\end{equation}

Summing the estimates we get the lemma. \myendproof

\begin{lemma}
\label{Pipmest} The projection operators $\Pi_\pm$ are bounded
from $L^2_s$ to $L^2_{-s}$, $s>1$, and from $L^p$ to $L^p$ for any
$1 \le p \le \I$.
\end{lemma}

\myproof From the definition of $K_\pm$ in Lemma \ref{Kpmest}, we
have
\begin{equation}
K_+ = \Pi_+ (J - i), \quad K_- = (1-\Pi_+-\Pi_0) (J+i),
\end{equation}
where $\Pi_0 = \sum_{j \ge m} P_{j}$ is bounded in $L^p$. Thus
\begin{equation}
\Pi_+ = \frac i2[K_+ + K_- - (1-\Pi_0) (J+i)],
\end{equation}
where shows $\Pi_+$ is bounded in $L^p$ for $p<\I$ by Lemma
\ref{Kpmest}. Similarly $\Pi_-$ and $\Pi_\pm^*$ are  bounded in
$L^p$ for $p<\I$. The boundedness of $\Pi_\pm$ in $L^\I$ follows
from that of $\Pi_\pm ^*$ in $L^1$ and duality. \myendproof

As a corollary, Lemmas \ref{th:decay} and \ref{th:sdecay} hold
with $\Pcsp$ replaced by $P_{\pm}$ since $P_{\pm} = \Pcsp\Pi_\pm $.
%===============================================================
%%============================================================================================================
\subsection{Fermi Golden Rule}

In this subsection we prove Corollary \ref{FRG-cor2}, which gives the
key resonance coefficients in the normal form equations in Lemmas
\ref{L-NF} and \ref{b.NF}.

For any $k \not=m$, recall \eqref{eq2-34} that
\begin{equation} \Phi_k = \smat{1 \\ -i} \bar{u}_k^+ + \smat{1\\i} \bar{u}_k^-.\end{equation}
From \eqref{u-pm}, we introduce $\Phi_k^+$ and $\Phi_k^-$ which
satisfy the equation $\Phi_k = \Phi_k^+ + \Phi_k^-$ where $\Phi_k^+$
is localized and
\begin{equation}\label{Phik-}
 \Phi_k^- = \smat{1\\i} (H-\bar{\alpha}_k)^{-1}\bar{\phi}_k^*, \quad \Phi_k^+ = \smat{1 \\ -i}\phi_k + O_{L^2_r}(n^2),
\end{equation}
Note that $\phi_k^* =O_{L^2_r}(n^2)$ is defined in \eqref{u-pm}
and $\alpha_k = i\bar{\LLLL }_k = |e_k -e_m| + O(n^2)$ with $\Im
\al_k
>0$. Moreover, since $\Phi_k = \Phi_k^+ + \Phi_k^-$, from
\eqref{Pk-def}, we see that for all function $f \in
L^2(\mathbb{R}^2,\mathbb{C}^2)$
\begin{equation} \label{Pk-ad}
\begin{split}
P_kf & = c_k(\sigma_1\bar{\Phi}_k,f)\Phi_k^+ + \bar{c}_k(\sigma_1\Phi_k,f)\bar{\Phi}_k^+  + c_k(\sigma_1\bar{\Phi}_k,f)\Phi_k^- + \bar{c}_k(\sigma_1\Phi_k,f)\bar{\Phi}_k^-, \\
(P_k)^*f & = c_k(\bar{\Phi}_k,f)\sigma_1 \Phi_k^+ +
\bar{c}_k(\Phi_k,f)\sigma_1 \bar{\Phi}_k^+  +
c_k(\bar{\Phi}_k,f)\sigma_1\Phi_k^- +
\bar{c}_k(\Phi_k,f)\sigma_1\bar{\Phi}_k^-.
\end{split}
\end{equation}
Since $\Phi_k^+$ is localized and $\Phi_k^- = O_{L^2\loc}(n^2)$, it
follows from Lemma \ref{JH} and \eqref{Pk-ad} that for all
function $f$ such with $\norm{f}_{L^2_r} = O(\delta)$
\begin{equation} \label{Pk-PkJH}
\begin{split}
(P_k - P_k^{JH}) f & = O(n^2\delta)\Phi_k + O(n^2\delta)\bar{\Phi}_k + O(\delta) \Phi_k^{-} + O(\delta) \bar{\Phi}_k^{-} + O_{L^2_r}(n^2\delta)\\
(P_k - P_k^{JH})^* f & = O(n^2\delta)\sigma_1 \Phi_k +
O(n^2\delta)\sigma_1 \bar{\Phi}_k + O(\delta)  \sigma_1 \Phi_k^{-}
+ O(\delta) \sigma_1 \bar{\Phi}_k^{-} + O_{L^2_r}(n^2\delta).
\end{split}
\end{equation}
Throughout this subsection, let $\omega$ and $\ep$ be two fixed
numbers such that
\begin{equation} \label{al-ome}
 \omega \pm \Im \LLLL_k =O(1) \not=0, \quad  0 < \ep \ll 1.
\end{equation}
Let $\al = -i \omega + \ep$ and
\begin{equation} \label{R-R0}
R = (\bL -\al)^{-1}, \quad R_0 = (JH -\al)^{-1}.
\end{equation}
Note that we have
\begin{equation} \label{R-exp}
R = R_0 + R_0 WR_0 + R_0 WRWR_0,
\end{equation}
where $W$ is a localized potential which is of order $\norm{Q}^2$.
\begin{lemma} \label{R0-R0*} For any $k \not= m$, there exist $C>0$ independent of $\ep$ and $n$ such that
\begin{equation}
\norm{R_0 \Phi_k^-}_{L^2\loc}, \quad \norm{R_0
\bar{\Phi}_k^-}_{L^2\loc}, \quad \norm{(R_0)^*
\sigma_1\Phi_k^-}_{L^2\loc}, \quad  \norm{(R_0)^*
\sigma_1\bar{\Phi}_k^-}_{L^2\loc} \leq  C n^2.
\end{equation}
\end{lemma}
\myproof We write
\begin{equation}
R_0 = (JH -\al)^{-1} = (H^2 + \al^2)^{-1} \begin{bmatrix} -\al & - H
\\ H & -\al \end{bmatrix}.
\end{equation}
Then, it follows that
\begin{equation} \label{R0-Phik}
\begin{split}
& R_0 \Phi_k^-   = \smat{-i \\1} (H + i\al)^{-1}(H -\bar{\al}_k)^{-1} \bar{\phi}_k^*\\
& R_0 \bar{\Phi}_k^-   = \smat{i \\1} (H - i\al)^{-1}(H -\al_k)^{-1} \phi_k^* \\
& (R_0)^* \sigma_1 \Phi_k^-   = \smat{1 \\ -i} (H + i\bar{\al})^{-1}(H -\bar{\al}_k)^{-1} \bar{\phi}_k^*, \\
& (R_0)^*\sigma_1 \bar{\Phi}_k^-   = \smat{1 \\ i} (H -
i\bar{\al})^{-1}(H -\al_k)^{-1} \phi_k^*.
\end{split}
\end{equation}
Since $\Re \al >0$ and $\Im (\al_k) >0$ and $\phi_k^* \in
O_{L^2_r}(n^2)$, our claim follows. \myendproof
\begin{lemma}\label{FG-err} There exists $C>0$ such that for any function $f, g \in L^2(\mathbb{R}^2, \mathbb{C}^2)$ with $f, g = O_{L^2_r}(n)$:
\begin{equation} \label{Pc-PJH-est}
\begin{split}
& |(f, (\bL -\al)^{-1} \Pcsp (\Pcsp  - P_c^{JH}) g)| \leq C n^4,\\
& |(f, (\Pcsp  -P_c^{JH})  (\bL -\al)^{-1}\Pcsp g)| \leq C n^4.
\end{split}
\end{equation}
\end{lemma}
\myproof Since the proof of both estimates in \eqref{Pc-PJH-est} are
similar, we shall only prove the first estimate of
\eqref{Pc-PJH-est}. From \eqref{Pk-PkJH}, we have
\begin{equation}(\Pcsp  - P_c^{JH}) g = \sum_{k=0}^{K}\{ O(n^3)\Phi_k + O(n^3)\bar{\Phi}_k + O(n) \Phi_k^{-} + O(n) \bar{\Phi}_k^{-}\} + O_{L^2_r}(n^3). \end{equation}
Since $\bL \Phi_k = \LLLL_k \Phi_k$ and $\LLLL_k -\al$,
$\bar{\LLLL }_k -\al$ are all non-zero order one, we get
\begin{equation}
 (f, \Pcsp R (\Pcsp  - P_c^{JH})g) =  O(n^4) + (f, \Pcsp  R [O(n) \Phi_k^- +O(n) \bar{\Phi}_k^{-}]).
\end{equation}
By similarity, we only need to show that $|(f, \Pcsp
R\Phi_k^-)| \leq Cn^3$. Let $\wt{g} = [WR_0 + WRWR_0]\Phi_k^-$. By
Lemma \ref{R0-R0*}, $\norm{\wt{g}}_{L^2_r} \leq C n^4$. Then, using
\eqref{R-exp}, \eqref{Pk-ad} and Lemma \ref{R0-R0*}, we have
\begin{equation} \label{err.Pk}
\begin{split}
& |(f, \Pcsp R\Phi_k^-)|  = |((\Pcsp )^*f, R_0 \Phi_k^- +  R_0\wt{g})|\\
& \leq |((P_d)^*f,  R_0\Phi_k^-) + ((R_0)^*(P_d)^*f, \wt{g})| + C n^3 \\
& \leq C\Big \{n\sum_{j\not=m}|(\sigma_1\Phi_j^-, R_0\Phi_k^-)| + n\sum_{j\not=m}^K|(\sigma_1\bar{\Phi}_j^-, R_0\Phi_k^-)| + n^3 \Big \}\\
& \leq C\Big \{n\sum_{j\not=m}^K|(\sigma_1\Phi_j^-, R_0\Phi_k^-)| +
n\sum_{j\not=m}^K|(\sigma_1\bar{\Phi}_j^-, R_0\Phi_k^-)| + n^3 \Big
\}
\end{split}
\end{equation}
Note that from \eqref{Phik-} and \eqref{R0-Phik}, we get
\begin{equation}
 (\sigma_1\Phi_j^-, R_0\Phi_k^-)| \leq Cn^4 , \quad (\sigma_1\bar{\Phi}_j^-, R_0\Phi_k^-) =0.
\end{equation}
So, from \eqref{err.Pk}, we obtain
\begin{equation}
||(f, \Pcsp  R\Phi_k^-)|| \leq Cn^3.
\end{equation}
This completes the proof of Lemma \ref{FG-err}. \myendproof
%========================================
\begin{corollary} \label{FGR.cor} For any function $f, g \in L^2(\mathbb{R}^2, \mathbb{C}^2)$ with $f,g = O_{L^2_r}(n)$, we have
\begin{equation}
(f, \Pcsp (\bL -\al)^{-1}\Pcsp  g) = (f, P_c^{JH_0} (J(H_0 -
E)-\al)^{-1}P_c^{JH_0} g) + O(n^4).
\end{equation}
\end{corollary}
\myproof Using \eqref{R-exp} and Lemma \ref{FG-err}, we have
\begin{equation}
(f, \Pcsp (\bL -\al)^{-1}\Pcsp  g) = (f, P_c^{JH} (JH
-\al)^{-1}P_c^{JH} g) + O(n^4).
\end{equation}
Now, since that $H - (H_0 -E) = \kappa  Q^2 =O(n^2)$ and $P_c^{JH}
- P_c^{JH_0} = O_{L^2_r}(n^2)$, we can use the same method as in
Lemma \ref{FG-err} to obtain
\begin{equation}
(f, \Pcsp (\bL -\al)^{-1}\Pcsp  g) = (f, P_c^{JH_0} (J(H_0 - E)
-\al)^{-1}P_c^{JH_0} g) + O(n^4).
\end{equation}
This completes the proof of Corollary \ref{FGR.cor}. \myendproof
\begin{corollary} \label{FRG-cor2} Let $f, g \in L^{2}(\mathbb{R}^3, \mathbb{C})$ be a localized real functions of order $n$ and let $f_1 = \smat{1 \\ i} f$ and $g_1 = \smat{i\\ 1} f$. We then have
\begin{equation}
\begin{split}
& \wei{\smat{1 \\ -i} f, (\bL -\al)^{-1}\Pcsp  \smat{i \\ 1} g} =-2(f, P_c^{H_0}(H_0-E -i\al)^{-1}P_c^{H_0}g) +O(n^4),\\
& \wei{\smat{1 \\ i} f, (\bL -\al)^{-1}\Pcsp  \smat{-i \\ 1} g} =-2(f, P_c^{H_0}(H_0-E +i\al)^{-1}P_c^{H_0}g) +O(n^4),\\
& \wei{\smat{1 \\ i} f, (\bL -\al)^{-1}\Pcsp  \smat{i \\ 1} g} = O(n^4), \\
& \wei{\smat{1 \\ -i} f,(\bL -\al)^{-1}\Pcsp \smat{-i \\ 1} g} = O(n^4).
\end{split}
\end{equation}
\end{corollary}
\myproof By Corollary \ref{FGR.cor}, we have
\begin{equation}  \wei{\smat{1 \\ -i} f, (\bL -\al)^{-1}\Pcsp  \smat{i \\ 1} g} =  \wei{\smat{1 \\ -i} f, P_c^{JH_0}(J(H_0-E) -\al)^{-1}P_c^{JH_0}  \smat{i \\ 1} g} + O(n^4).\end{equation}
On the other hand,
\begin{flalign*}
(J(H_0-E) -\al)^{-1}P_c^{JH_0}  \smat{i \\ 1} g & = (H_0-E +\al^2)^{-1}P_c^{H_0} \begin{bmatrix} -\al & -(H_0-E) \\ H_0 -E & -\al \end{bmatrix} \smat{i \\ 1} g \\
& = \smat{-1 \\ i} \begin{bmatrix} (H_0 -E -i\al)^{-1} \\ (H_0 -E -i\al)^{-1} \end{bmatrix} P_c^{H_0}g. \end{flalign*}
So, the first identity of our corollary follows. Similarly, we can prove all of the last three identities of the corollary.
\myendproof
%=========================================================================
% START OF SEC.3
%=========================================================================
\section{Equations and main terms}
In our analysis we use different coordinate systems. When the
solution is away from bound states, we use the {\it orthogonal
coordinates} \eqref{1-6}, i.e., we decompose the solution as a sum
of different spectral components with respect to $-\De+V$. When
the solution is near a nonlinear bound state, we use the {\it
linearized coordinates} \eqref{psi.dec}, i.e., decomposition with
respect to the corresponding linearized operator instead.  In
subsection \ref{orthcoor} we recall the equations and normal forms
in orthogonal coordinates from \cite{Tsai}. The rest of this
section is devoted to analysis in linearized coordinates. We will
not use  {\it centered orthogonal coordinates}
\eqref{centeredortho}, which is also mentioned in \S1.
\subsection{Orthogonal coordinates}
\label{orthcoor}

Let $t_0$ be a fixed initial time. For $t \geq t_0$ we may decompose
the solution with respect to $H_0$ as
\begin{equation}
\psi(t) = \sum_{j = 0}^{K} x_j(t) \phi_j + \xi, \quad \xi \in
\mathbf{H}_c(H_0), \ \forall \ t \geq t_0.
\end{equation}
Then for $t \geq t_0$, as in \cite[Section 4]{Tsai} we have
\begin{equation}
\begin{split}
i \dot{x}_j & = e_j x_j + (\phi_j, G),  \quad (j =0, \ldots , K), \\
i \partial_t \xi & = H_0 \xi + \Pc^{H_0} G,  \quad G := \kappa
\psi^2 \bar{\psi}.
\end{split}
\end{equation}
Let
\begin{equation} \label{G.def}
 G_3 :=  \kappa  \left |\sum_{j=0}^{K} x_j\phi_j \right |^2 \left (\sum_{j=0}^{K} x_j \phi_j \right )
= \kappa  \sum_{l,m,j=0}^K x_lx_m\bar{x}_j \phi_l\phi_m\phi_j.
\end{equation}
We then decompose $\xi$ as (for details, see \cite[Section 4]{Tsai})
\begin{equation} \label{xi.dec}
\xi(t) = \xi^{(2)}(t) + \xi^{(3)}_1(t) +  \xi^{(3)}_2(t) + \cdots + \xi^{(3)}_5(t), \ \forall \ t \geq t_0,
\end{equation}
where
\begin{equation}  \label{xi2}
\begin{split}
 & \xi^{(2)}(t)  : = \sum_{l,m,j =0}^{K} x_lx_m\bar{x}_j(t) \xi_{lm}^j, \ \text{with}\\
 & \xi_{lm}^j : = -\kappa  \lim_{r \rightarrow 0^+} [H_0 - e_l - e_m + e_j -ri]^{-1} \Pc^{H_0} \phi_m\phi_l\phi_j,
\end{split}
\end{equation}
and, with $u_j(t) = e^{ie_jt} x_j(t)$ which have less oscillation
than $x_j(t)$,
\begin{equation} \label{xi.all}
\begin{split}
& \xi_1^{(3)}(t) := e^{-iH_0(t-t_0)} \xi(t_0), \qquad \xi_2^{(3)}(t) := -e^{-iH_0(t -t_0)} \xi^{(2)}(t_0), \\
& \xi^{(3)}_3(t) := - \int_{t_0}^t e^{-iH_0(t-s)} \Pc^{H_0} \sum_{l,m,j =0}^K e^{i(-e_l-e_m + e_j)s} \frac{d}{ds} (u_lu_m\bar{u}_j) \xi_{lm}^j ds, \\
& \xi^{(3)}_4(t) := \int_{t_0}^ t e^{-iH_0 (t-s)} \Pc^{H_0} i^{-1} (G - G_3 - \kappa  \xi^2 \bar{\xi}) ds, \\
& \xi^{(3)}_5(t) := \int_{t_0}^ t e^{-iH_0 (t-s)} \Pc^{H_0} i^{-1}
(\kappa  \xi^2 \bar{\xi}) ds.
\end{split}
\end{equation}
We  recall the following two lemmas from \cite{Tsai}:
\begin{lemma}[Lemma 4.1 \cite{Tsai}] \label{nonlin.est} Let $p, p'$ such that $4\leq p <6,\ (p)^{-1} + (p')^{-1} = 1$. Suppose that for a fixed time $t \geq t_0$ and for $0 < n \leq n_0 \ll 1$, we have
\begin{equation} \max_{j} |x_j(t)| \leq 2n, \quad \norm{\xi(t)}_{L^2_{loc} \cap L^p} \leq 2n, \quad \norm{\xi(t)}_{L^2} \ll 1. \end{equation}
Then for $u_j(t) = e^{ie_jt} x_j(t)$,
\begin{equation}
\begin{split}
\norm{G}_{L^1_{loc}} + \max_{j}|\dot{u}_j|  \les n^3 \quad
\text{and} \quad \norm{G - G_3 - \kappa  \xi^2\bar{\xi}}_{L^1 \cap
L^{p'}}  \les n^2 \norm{\xi}_{L^2_{loc}}.
\end{split}
\end{equation}
\end{lemma}
%=====================================================================
%\myproof The proof of this Lemma is similar to that of \cite[Lemma 4.1]{Tsai}. So, we skip it.
%\myendproof \\ %\ \\
%======================================================================
%Also, we need the following normal form for the equation of $\dot{u}_j$.
\begin{lemma}[Lemma 4.2 \cite{Tsai}] \label{normal-form-u} Let $p, u_j$ be as in Lemma \ref{nonlin.est}. Suppose that for some $0 < n \leq n_0$ and for some $t \geq t_0$,
\begin{equation} \max_{j}|x_j(t)| \leq 2n, \quad \norm{\xi(t)}_{L^2_{loc} \cap L^p} \leq 2n \quad \text{and} \quad \norm{\xi(t)}_{L^2} \leq \alpha \ll 1.  \end{equation}
Then, there are perturbations $\mu_j(t)$ of $u_j(t)$, $ j \in I$,
such that
\begin{equation} \label{muj.eqn}
\dot{\mu}_j(t) = \sum_{l=0}^K c_l^{j} |\mu_l|^2 \mu_j + \sum_{a,b=0}^K d_{ab}^j |\mu_a|^2 |\mu_b|^2 \mu_j + g_j,
\end{equation}
and
\begin{equation} \label{error.est.mu}
\begin{split}
& |u_j(t) - \mu_j(t)| \les n^3, \\
& |g_j(t)| \les n^7 + n^2 \norm{\xi^{(3)}}_{L^2_{loc}} + n \norm{\xi}_{L^2_{loc}}^2 + \norm{\xi}_{L^2_{loc}}^{\frac{2(p-3)}{p-2}}\norm{\xi}_{L^p}^{\frac{p}{p-2}}.
\end{split}
\end{equation}
Moreover, all of the coefficients $c_l^j$ and $d_{ab}^j$ are of order one. The coefficients $c_l^j$ are all purely imaginary and
\begin{equation}
\label{Redabj}
\Re d_{ab}^j = (2 - \delta_a^b)\gamma_{ab}^j -
2(2-\delta_j^b)\gamma_{jb}^a,
\end{equation}
with $\delta_a^b=1$ if
$a=b$ and $\delta_a^b=0$ if $a \not = b$, and
\begin{equation}
\gamma_{ab}^l = \kappa ^2 \Im \bke{\phi_a\phi_b\phi_l, (H_0 - e_a - e_b + e_l - i0^+)^{-1}
\Pc^{H_0} \phi_a\phi_b\phi_l}, \ \forall \ a,b,l \in I.
\end{equation}
\end{lemma}

By the resonance condition Assumption (A2), the number
$\gamma_{ab}^l \geq 0$ and it is positive if and only if $ l < a,b$.

%=================================================
\subsection{Linearized coordinates}
\label{S3.2} When the solution $\psi$ lies in a neighborhood of an
excited state $Q=Q_{m,n}$, $m \in J$, it is natural to decompose
$\psi-Q$ into invariant subspaces of  the linearized operator
around $Q$, see Lemma \ref{L-spectral}. The collection of these
components is called the {\it linearized coordinates}.

\begin{lemma}\label{Linearized.dec}
There are small positive constants $n_0$ and $\e_3$ such that the
following hold. Suppose $\norm{\psi}_{H^1} \leq n_0$ satisfies
$\norm{\psi - (\psi,\phi_m) \phi_m}_{L^2}\le \e_3
|(\psi,\phi_m)|$.
\begin{itemize}
\item[\textup{(i)}]
For any $0 < n < n_0$, there exist unique $a, \theta \in \R$ such
that
\begin{equation}
\psi = [Q_{m,n} + a R_{m,n} + h ]e^{i\theta},
\end{equation}
where $Q_{m,n}$ and $R_{m,n}$ are given by Lemma \ref{th:2-1},
$P_m h = 0$, and $|n^{-1}a| + \norm{h}_{H^1} \leq \e_3 n$.

\item[\textup{(ii)}]
There exist unique $n(\psi) \in (0,n_0)$ and $\theta \in \R$ such
that $a=0$. Moreover, if $\psi$ is decomposed as in (i) with respect
to another $n$, then
\begin{equation}
n(\psi) = n + \frac a{2 C n} + O(n^3), \quad C = \kappa  \int
\phi_m^4.
\end{equation}

\item[\textup{(iii)}] If $\psi$ is decomposed as in (i) with respect
to $n_1$ and $n_2$ with $\norm{h_j} \le \rho \le \e_3 n$,
$|a_j|\le C \rho^2$, and $|n_1-n_2|\lec n^{-1}\rho^2$, then
\begin{equation}
C (n_1^2 - n_2^2)+ a_1 - a_2 = O(\rho |n_1-n_2|).
%E_{1,n_1} + a_1 - E_{1,n_2} - a_2 = O(\rho \al).
\end{equation}
\end{itemize}
\end{lemma}

The proof of Lemma \ref{Linearized.dec} is similar to those for
\cite[Lemmas 2.1--2.4]{TY1}.

%\subsubsection{Equations}

By Lemma \ref{Linearized.dec}, when $\psi(t)$ is a sufficiently small neighborhood of an excited state $Q=Q_{m,n}$,
there is a unique choice of real $a(t)$ and $\th(t)$ so that
\begin{equation} \label{psi.dec}
\psi(t)  =[Q + a(t)R + h(t)]e^{-iEt + i \th(t)}, \quad P_m h(t) = 0.
\end{equation}
Here $R= R_{m,n}$ and $E=E_{m,n}$. We can further decompose
\begin{equation}
h = \zeta + \eta, \quad \zeta = \sum_{ k \not= m} \zeta_k, \quad
\eta \in \bE_c^\sharp ,
\end{equation}
where, for each $k \not = m$,
\begin{equation}
\zeta_k : = \bj^{-1}\Re (z_k \Phi_k) = z_k \bar{u}_k^+ + \bar{z}_k
u_k^-, \quad u_k^{\pm} : = \frac{1}{2}(\bar{u}_k \pm \bar{v}_k).
 \end{equation}
%Thus we have
%\begin{equation}
%\psi(t)  =[Q + a(t)R + \zeta(t) + \eta(t) ]e^{-iEt + i \th(t)}.
%\end{equation}

Substituting \eqref{psi.dec} into \eqref{Sch} and using $\L iQ = 0$
and $\L R = -iQ$, we get
\begin{equation} \label{h.eqn}
\partial_t h  - \L h = F_h \equiv
i^{-1}(F + \dot{\theta}(Q+aR +h)) - aiQ - \dot{a}R,
\end{equation}
where
\begin{equation} \label{F.def}
F = \kappa  Q(2|h_\sigma|^2 + h_\sigma^2) + \kappa  |h_\sigma|^2
h_\sigma ,\quad  h_\sigma = aR +h.
\end{equation}
We choose $\dot \th$ and $\dot a$ so that $P_m F_h=0$. Thus $F_h =
(1-P_m)i^{-1}(F + \dot{\theta}(aR +h))$ and
\begin{equation} \label{eq:all}
\left\{
\begin{aligned}
&\dot a = (c_m  Q, \Im (F +\dot \theta h) ) ,
\\
& \dot \theta = F_\theta \equiv -\bkt{a+ \bke{ c_m R,\, \Re F } }
\cdot \bkt{1+ (c_m R,R)a +(c_m R,\Re h)}^{-1}.
\end{aligned}
\right .
\end{equation}
Taking $\Pcsp $ of \eqref{h.eqn}, we get
\begin{equation} \label{eta.eqn}
\partial_t \eta  - \L \eta = \Pcsp  i^{-1}(F + \dot{\theta}(aR +h)).
\end{equation}
Note $z_k = 2c_k (\si_1 \bar \Phi_k ,[h])$. Taking $2c_k (\si_1
\bar \Phi_k ,[\cdot])$ of \eqref{h.eqn}, $k \not = m$, we get
\begin{equation}
\dot z_k - \LLLL_k z_k = Z_k: =2 c_k (\si_1 \bar \Phi_k , [F_h]).
\end{equation}
A direct computation using \eqref{eq2-34} shows\footnote{Note
$-2c_k \sim -i$ which is the coefficient of \cite[page 242, line
5]{Tsai}.}
\begin{equation}
\label{Z.def} Z_k  = -2c_k \left \{(u_k^+, F)+(u_k^-, \wbar F)+
\bkt{(u_k^+, h)+(u_k^-, \wbar h)+ (\bar{u}_k,R) a} \dot \theta
\right \}.
\end{equation}
Let $\om_k: =- \Im \LLLL_k$ and let $p_k(t) = z_k(t) e^{i\ev_k
t}$. We have
\begin{equation} \label{z.eqn}
\dot p_k =  (\Re \LLLL_k) p_k + e^{i\ev _k t}Z_k .
\end{equation}
Also, for any $k \not=m$, let $r_k := e^{-\LLLL_k t}z_k$, we have,
\begin{equation} \label{r.eqn}
\dot{r}_k = e^{-\LLLL_kt}Z_k.
\end{equation}
Note that $r_k = p_k$ for all $k > m$ and $r_k =
e^{-\Re(\LLLL_k)t}p_k$ for $k <m$. We shall use $r_k$ in computing
the normal form for the equation of $a$.
%=======================================================
%========================================================================
\begin{definition} \label{convention} Denote $I = \{0,1,\cdots, K\}$,
$I^*= \{0^*, 1^*,\cdots, K^*\}$. For all $m \in I$, let $I_{>m} =
\{m+1, \cdots , K\}$, $I_{<m} = \{0,\cdots, m-1\}$,  $I_m = I
\setminus \{m\}$,  $I_m^* = I^* \setminus \{m^*\}$ and $\Omega_m: =
I_m \cup I^*_{m}$. For $j \in I_{m}$, let
\begin{equation}  %\label{S-14}
\LLLL_{j^*} = \bar{\LLLL }_{j}, \quad \ev_{j^*} = - \ev _j , \quad
z_{j^*} = \bar z_j, \quad r_{j^*} =\bar{r}_j, \quad  p_{j^*} =
\bar p_j, \quad u_{j^*}^{\pm} = \bar{u}_{j}^{\pm}, \quad
\text{and} \quad  v_{j^*}^{\pm} = \bar{v}_{j}^{\pm}.
\end{equation}
\end{definition}
It then follows that for all $j \in \Omega_m$, we have $z_{j}(t) =
e ^{-i \ev_j t} p_j(t)$ and $r_{j} = e^{-\LLLL_j t}z_j$.
%=====================================================
%=====================================================
%=====================================================
\subsection{Decomposition of $a$}
Recall $\dot a = (c_m  Q, \Im (F +\dot \theta h) )$. Let $F_1 :=
\kappa Q ( 2|\zeta|^2 + \zeta^2 )$, $A^{(2)} :=  c_m (Q, \Im F_1
)$ and $A^{(3)} :=   c_m  (Q, \Im (F-F_1 + \dot{\theta} h))$.
Then, we have $\dot{a} = A^{(2)} + A^{(3)}$. We shall impose the
boundary condition of $a$ at $t= T$, which is in fact the
condition imposed on the choice of $E = E(T)$. Hence, we have
\begin{equation} a(t) = a(T) + \int_T^t [A^{(2)}(s) + A^{(3)}(s)] ds. \end{equation}
Recall that
\begin{equation} \zeta = \sum_{k \in \idx_{m}} \zeta_k, \quad \zeta_k = z_k \bar{u}_k^{+} + \bar{z}_k u_k^-.  \end{equation}
Therefore,
\begin{equation} \Im \zeta_k \zeta_l = \Im [(z_k z_l)(\bar{u}_k^+\bar{u}_l^+ -\bar{u}_k^-\bar{u}_l^-) + (z_k \bar z_l)(\bar{u}_k^+u_l^- -\bar{u}_k^-u_l^+)]. \end{equation}
Let
\begin{equation}a_{kl,1} : = \kappa  c_m  (Q^2, \bar{u}_k^+\bar{u}_l^+ -\bar{u}_k^-\bar{u}_l^-), \quad a_{kl,2} = \kappa  c_m  (Q^2,\bar{u}_k^+u_l^- -\bar{u}_k^-u_l^+ ). \end{equation}
Note that $a_{kl,1}, a_{kl,2} = O(n^2)$, $a_{kl,1}, a_{kl,2}$ are
real if both $k,l >m$, and $a_{kk,2}$ are purely imaginary. In
particular $a_{kk,2}=0$ if $k>m$. We have
\begin{equation} \label{A(2).exp}
\begin{split}
A^{(2)} & = \kappa  c_m  (Q^2, \Im \sum_{k,l \in \idx_{m}} \zeta_k
\zeta_l )
 =  \Im \sum_{k,l \in \idx_{m}} \left \{a_{kl,1} z_k z_l
 + a_{kl,2}z_k \bar{z}_l \right \} \\
 & = b_0(t) + \Im (A_1^{(2)}),
\end{split}
\end{equation}
where
\begin{equation} \label{b0.def}
b_0(t)= \sum_{k {<m}}  b_{0k}|z_k|^2, \quad b_{0k}: =\Im a_{kk,2},
\quad \wt{b}_0(t) : =\int_{T}^tb_0(s)ds,
\end{equation}
\begin{equation}
A_1^{(2)} : =\sum_{k,l \in \idx_{m}} a_{kl,1} z_k z_l + \sum_{k \not=l}  a_{kl,2}z_k \bar{z}_l.
\end{equation}
%Due to the $\Im$-operator, the middle term in the right hand side of
%\eqref{A(2).exp} equation is killed.
Note  $|b_{0k}| \lec  n^2\norm{u_k^-}_{L^2_{loc}} =O(n^4)$ for $k
<m$ by Lemmas \ref{L-spectral} and \ref{uv-est}.

We shall integrate  $A_1^{(2)}$  by parts. Note that for all
$\LLLL_k + \LLLL_l = -i(\om_k + \om_l) + O(n^4)$ and $\LLLL_k +
\bar{\LLLL }_l = -i(\om_k - \om_l) + O(n^4)$. Therefore, $\LLLL_k
+ \LLLL_l = O(1)$ for all $k, l \in \idx_m$ and $\LLLL_k +
\bar{\LLLL }_l = O(1)$ for all $k, l \in \idx_m$ and $k \not=l$.
We then write
\begin{equation} \label{A(2).int}
\begin{split}
A_1^{(2)} & =\sum_{k,l \in \idx_{m}} a_{kl,1} e^{(\LLLL_k +\LLLL_l)t}r_k r_l + \sum_{k \not=l}  e^{(\LLLL_k +\bar{\LLLL }_l)t}a_{kl,2}r_k \bar{r}_l \\
& = \sum_{k,l \in \idx_{m}} \frac{a_{kl,1}}{\LLLL_k +\LLLL_l} \left[\frac{d}{dt}(z_k z_l)- e^{(\LLLL_k +\LLLL_l)t}\frac{d}{dt}(r_k r_l) \right ]
\\
& \quad + \sum_{k \not=l} \frac{a_{kl,2}}{\LLLL_k +\bar{\LLLL }_l}
\left[\frac{d}{dt}(z_k \bar{z}_l)- e^{(\LLLL_k +\bar{\LLLL
}_l)t}\frac{d}{dt}(r_k \bar{r}_l) \right ].
\end{split}\end{equation}
Now, define
\begin{equation} \label{a(2).def}
\begin{split}
a^{(2)} (t) & : = \Im \sum_{k,l \in \idx_{m}}
\frac{a_{kl,1}}{\LLLL_k +\LLLL_l} z_k z_l +\Im \sum_{k \not=l}
\frac{a_{kl,2}}{\LLLL_k +\bar{\LLLL }_l}z_k\bar{z}_l
\\
A_{2,rm} & :=   \Im \sum_{k,l \in \idx_{m}}
\frac{a_{kl,1}e^{(\LLLL_k +\LLLL_l)t}}{\LLLL_k
+\LLLL_l}\frac{d}{dt}(r_k r_l) + \Im \sum_{k \not=l}
\frac{a_{kl,2}e^{(\LLLL_k +\bar{\LLLL }_l)t}}{\LLLL_k +\bar{\LLLL
}_l}\frac{d}{dt}(r_k\bar{r}_l).
\end{split}
\end{equation}
We shall get
\begin{equation}
\begin{split}
\Im (A_1^{(2)})& = \frac{d}{dt} a^{(2)} (t) - A_{2,rm}(t).
\end{split}
\end{equation}
So, we have $A^{(2)} = \frac{d}{dt} a^{(2)}(t) + b_0(t) - A_{2,rm}(t)$. Therefore,
\begin{equation} \label{a.dec}
a(t) = a^{(2)}(t) +  b(t),
\end{equation}
where $b(t)$ satisfies
\begin{equation} \label{eq:b1}
\dot {b} = {b}_0 +  c_m (Q, \Im (F -F_1 + \dot{\theta}h)) -
A_{2,rm}, \quad b(T) = a(T) - a^{(2)}(T).
\end{equation}
Moreover, let $a_{kl,3}:= 2a_{kl,1}(\LLLL_k +\LLLL_l)^{-1}$ and
$a_{kl,4}:= 2a_{kl,2}(\LLLL_k +\bar{\LLLL }_l)^{-1}$. Since
$a_{kl,1}$ and $a_{kl,2}$ are of order $n^2$, so are $a_{kl,3}$
and $a_{kl,4}$.  Moreover, $a_{kl,3}, a_{kl,4}$ are purely
imaginary for $k, l \in \idx_{>m}$. Using \eqref{r.eqn}, $a_{kl,1}
= a_{lk,1}$ and $a_{kl,2} = -\bar{a}_{lk,2}$, we obtain
\begin{equation} \label{A2rm.def}
A_{2,rm} =  \Im \sum_{k,l \in \idx_{m}} a_{kl,3}z_k Z_l
+ \Im \sum_{k \not=l} a_{kl,4} Z_k \bar{z}_l.
\end{equation}
It worths noting that the benefits from using $r_k$ instead of $p_k$
in \eqref{a(2).def} is that we do not have terms of order $zz_k$ for
$k \in \idx_{<m}$ in \eqref{A2rm.def}. This is very essential in the
normal forms.
%===============================================================================================
%========================================================================================

\subsection{Decomposition of $\eta$}
We shall single out the
main terms in $\eta$. Recall from \eqref{eta.eqn} that
\begin{equation} %\label{eta.eqn}
\partial_t \eta  - \L \eta = \Pcsp  i^{-1}(F + \dot{\theta}(aR +\zeta +\eta )).
\end{equation}
In the vector form, we have
\begin{equation} \label{eta.eq2}
\pd_t [\eta ]= \bL [\eta] + \Pcsp  J\dot \th [\eta] + \Pcsp J[ (F +
\dot{\theta}(aR +\zeta))].
\end{equation}
We first deal with the non-localized linear term $J \dot \th [\eta]$
using Lemma \ref{Kpmest}, following Buslaev-Perelman \cite{BP2},
also see \cite{BS, C2}.\footnote{The term $i \dot \th \eta$ is not a
problem in \cite{TY1} in which $\L$ is factorized in the form $\L =
U^{-1}JAU$ for some scalar self-adjoint operator $A$. Such
factorization does not exist for linearized operators near excited
states. In \cite{TY3}, the term $i \dot \th \eta$ is removed by
introducing $\td \eta = \Pcsp  e^{i \th} \eta$ and using Strichartz
estimates to control the (small) commutator term. This last method
is not suitable for $L^p$-decay approach since the commutator term,
although smaller, has the same decay rate as $\eta$ itself. The
approach of Buslaev-Perelman has the further benefit of being
applicable to large soliton case.} We need to revise their original
statement and proof to take care of eigenvalues near the continuous
spectrum.

Recall $P_\pm$ are defined in subsection \ref{upperlower}. Taking
projection $P_\pm$ of \eqref{eta.eq2}, and using
\begin{equation}
P_\pm J \mp i P_\pm =P_\pm (P_\pm J \mp i P_\pm)= P_\pm [K_\pm -
(P_{R\pm} J \mp i P_{R\pm})] = P_\pm K_\pm,
\end{equation}
 we get
\begin{equation}
\pd_t P_\pm [\eta] = \bL P_\pm [\eta] \pm i\dot \th P_\pm [\eta] +
P_\pm K_\pm \dot \th [\eta] + P_\pm J[ (F + \dot{\theta}(aR
+\zeta))].
\end{equation}
Denote
\begin{equation}
\eta_\pm := e^{\mp i \th} P_\pm [ \eta].
\end{equation}
We have
\begin{equation} \label{eta.eq3}
\pd_t \eta_\pm = \bL \eta_\pm  + e^{\mp i \th} P_\pm \bkt{K_\pm
\dot \th [\eta] + J[ (F + \dot{\theta}(aR +\zeta))]}.
\end{equation}
Recall that $[\zeta_k] = (z_k \Phi_k + \bar{z}_k \bar{\Phi}_k)/2$.
Note the term $e^{\mp i \th} P_\pm  J \dot{\theta} [\zeta ]$ is
not localized. However, by formula \eqref{JPk-1a}
\begin{equation} \label{eq3-26}
\Pcsp  J \Phi_k = \Pcsp   \Phi_k ', \quad \Pcsp  J\bar{\Phi}_k = \Pcsp   \bar{\Phi}_k ' \quad \Phi_k' = \smat{-2i\\ -2}
\bar u_k^+
\end{equation}
and note $\Phi_k '$ is localized. Thus we can rewrite the linear
terms in \eqref{eta.eq3} as
\begin{equation}\label{FL}
F_{L\pm} := e^{\mp i\theta} \dot{\theta}\bket{K_\pm [\eta] + J[aR]
+ \tsum_{j \in \idx_m} (z_j\Phi_j' + \bar{z}_j\bar{\Phi}_j')},
\end{equation}
where all functions are localized, and \eqref{eta.eq3} becomes
\begin{equation} \label{eta-pm.eqn}
\pd_t \eta_\pm = \bL \eta_\pm  + P_\pm \bkt{e^{\mp i\th}J[F]+
F_{L\pm}}.
\end{equation}
In other words, for some $t_0 \geq 0$ and for all $t \geq t_0$, we have
\begin{equation} \label{wt-eta.eqn}
\eta_\pm (t) = e^{\bL (t-t_0)}\eta_\pm (t_0)+ \int_{t_0}^t e^{\bL (t-s)} P_{\pm}
\{e^{\mp i\th}J[F] + F_{L\pm}\}(s) ds.
\end{equation}
We will decompose $\eta_\pm$ as follows.  Denote
\begin{equation}\label{wt-eta14}
\begin{split}
& \eta_{\pm,1}^{(3)}(t): = e^{\bL(t-t_0)}\eta_{\pm}(t_0), \\
& \eta_{\pm,4}^{(3)}(t):= \int_{t_0}^t e^{\bL (t-s)} P_{\pm}\{
F_{L\pm} + e^{\mp i \th} J[F - F_1] \}(s)ds .
\end{split}
\end{equation}
Then, we have
\begin{equation} \label{wt-eta.eqn2}
\eta_\pm(t) = \eta_{\pm,1}^{(3)}(t) + \eta_{\pm, 4}^{(3)}(t)+ \int_{t_0}^t
e^{\bL (t-s)} P_{\pm} \{ e^{\mp i\th} J[F_1]\}(s) ds.
\end{equation}
We shall integrate the last term in \eqref{wt-eta.eqn2}. Recall
that $F_1 = \kappa  Q(2|\zeta|^2 + \zeta^2)$ is the main term in
$F$ with
\begin{equation}
\zeta = \sum_{k \in \idx_{m}} \zeta_k = \sum_{k \in \idx_{m}}
(z_k\bar{u}_k^+ + \bar{z}_ku_k^-), \quad u_k^+ = \phi_k + O_{L^2_r}(n^2),\quad u_k^- = O_{L^2\loc}(n^2).
\end{equation}
So,
\begin{equation} \label{F1.dec}
F_1 = \sum_{k,l \in \idx_m} F_{kl}[z_kz_l + 2z_k\bar{z}_l] +
\sum_{k,l \in \Omega_m} \wt{F}_{kl}z_kz_l, \quad F_{kl} = \kappa
Q \phi_k \phi_l, \quad \wt{F}_{kl} =  O_{L^\infty_{3r}}(n^3).
\end{equation}
In other words, we can write
\begin{equation} \label{F1.com}
F_1 = \kappa  \sum_{k,l \in \Omega_m} z_kz_l \Phi_{kl},
\end{equation}
for some localized functions $\Phi_{kl}$ which can be computed explicitly.
 In particular, $\Re \Phi_{kl} = O(n)$ and $\Im\Phi_{kl} = O(n^3)$ for all $k, l \in \Om_m$.

To integrate $P_{\pm}  e^{\mp i\th}J[F_1]$ in $\eta_{\pm}$
equation, we want to integrate terms of the form
\begin{equation}
I_{\pm}(t) = \int_{t_0}^t e^{(t-s)\bL}  e^{-i\om s}P_{\pm} f(s)
\,ds,
\end{equation}
where $\om \in \R$, $f(s) \in L^2(\R^3, \C^2)$ and $\dot f(s)$
decays faster than $f$. We re-write $I_{\pm}$ as
\begin{equation}
I_{\pm}(t) = e^{t \bL}\int_{t_0}^t e^{-s(\bL + i\om)} P_{\pm}
f(s)ds.
\end{equation}
Denote $R= \lim_{\e \to 0+}(\bL + i \om -\e)^{-1}$. Integration by
parts gives
\begin{equation}
I_\pm(t) = - e^{-i \om t} RP_{\pm} f(t) +  e^{(t-t_0) \bL}e^{-i
\om t_0} RP_{\pm} f(t_0) + \int_{t_0}^t e^{(t-s) \bL} RP_{\pm}
e^{-i \om s}\dot f(s)ds.
\end{equation}
The choice of the sign of $\e$ ensures that $e^{t\bL} RP_\pm$ has
singular decay estimate according to Lemma \ref{th:sdecay}. We can
now identify the main term of $\eta_{\pm}$. Since $i^{-1}F_1 = -i
\kappa  \sum z_k z_l \Phi_{kl}$ with summation over $k, l \in
\Omega_m$,
\begin{equation}
J[F_1] = -\Re \tsum i \kappa  z_k z_l \Phi_{kl} \smat{ 1 \\ -i} =
-\Re  \tsum f_{kl}(s) e^{-i(\om_k + \om_l)s},
\end{equation}
where $f_{kl} = i \kappa  p_k p_l \Phi_{kl} \smat{ 1 \\
-i}$. We decompose $P_\pm = \Pi_\pm \Pc^\sharp$ since $\Pi_\pm$
does not commute with $\Re$. Denote $R_{ kl} = \lim_{\e \to
0+}(\bL + i(\ev_k + \ev_l)-\e)^{-1} \Pcsp $ and $\om_{kl}=\om_k +
\om_l$. We get
\begin{equation}
\int_{t_0}^t e^{(t -s)\bL }P_{\pm} e^{\mp i\theta(s)}J[F_1]ds =
\eta_{\pm}^{(2)} +  \eta^{(3)}_{\pm,2} + \eta^{(3)}_{\pm,3}
\end{equation}
where
\begin{equation}\label{wt-eta23.def}
\begin{split}
\eta_{\pm}^{(2)} & =  e^{\mp i\theta(t)}\Pi_\pm\Re  \tsum_{k, l
\in \Omega_m} R_{ kl} e^{-i\om_{kl}t} f_{kl} (t)
\\
\eta^{(3)}_{\pm,2} &= - e^{(t-t_0) \bL}e^{\mp i\theta(t_0)}
\Pi_\pm \Re  \tsum_{k,l \in \Omega_m} R_{ kl}  e^{-i\om_{kl}t_0}
f_{kl} (t_0)
\\
\eta^{(3)}_{\pm,3} & =- \int_{t_0}^t e^{(t-s) \bL} e^{\mp
i\theta(s)} \Pi_\pm  \tsum_{k,l \in \Omega_m} \big(\Re R_{
kl}e^{-i\om_{kl}s} \dot f_{kl}  \mp i\Re R_{ kl}e^{-i\om_{kl}s}
\dot{\theta}f_{kl} \big) (s)ds .
\end{split}\end{equation}
%Here $\eta^{(2)}_{\pm},  \eta^{(3)}_{\pm,2},  \eta^{(3)}_{\pm,3}$ respectively denotes
%the first, the second and the last sum in \eqref{wt-eta23.def}.
Observe that
\begin{equation} \label{fkl-th.est}
\||\dot f_{kl}| +|\dot{\theta}f_{kl}| \|_{L^2_r} \lec n|\dot \th|
\beta^2+ n\beta \max|\dot p_k|, \quad \beta=\max |p_k|.
\end{equation}
Now, let
\begin{equation} \label{wt-eta3.def}
\eta_{\pm}^{(3)}(t) := \sum_{j=1}^4\eta^{(3)}_{\pm,j}(t), \quad
\eta^{(j)} :=e^{i\theta} \eta^{(j)}_{+} + e^{-i\theta}
\eta^{(j)}_{-}, \quad j = 2,3.
\end{equation}
Then, from \eqref{wt-eta.eqn2} and \eqref{wt-eta23.def}, we obtain
the decomposition of $\eta_{\pm}$ and $\eta$ as
\begin{equation} \label{wt-eta.dec}
\eta_{\pm} = \eta_{\pm}^{(2)} + \eta_{\pm}^{(3)}, \quad [\eta] =
e^{i\theta}\eta_{+} + e^{-i\theta}\eta_{-} = \eta^{(2)} +
\eta^{(3)}.
\end{equation}

We now compute the explicit form of $\eta^{(2)}$ which will be
used in the computation of the key coefficients in the normal
forms of $z_k$. By \eqref{wt-eta.dec}, \eqref{wt-eta23.def},
$\Pi_+ + \Pi_-=\Pcsp$, and \eqref{F1.dec},
\begin{equation}\label{wt-eta2.dec}
\begin{split}
\eta^{(2)} & =e^{i\th}\eta^{(2)}_+ +e^{-i\th} \eta^{(2)}_- \\
&= \Re \tsum_{k, l \in \Omega_m} R_{ kl} e^{-i(\om_k + \om_l)t}
f_{kl} (t)
\\
&= \Re \sum_{k,l \in \idx_m} \Big \{ R_{kl} z_kz_l \smat{i \\ 1}
F_{kl} + 2R_{kl^*}z_k\bar{z}_l \smat{i \\ 1} F_{kl} \Big \} +
\sum_{k,l \in \Omega_m } z_kz_l R_{kl} O_{L^2_r}(n^3).
%\\
%& = \eta^{(2)}_1 + \eta^{(2)}_2.
\end{split}
\end{equation}
Recall $F_{kl} = \kappa  Q\phi_k\phi_l$. Thus the first sum
contains terms of order $O(nz^2)$.

%========================================================================================
\subsection{Decomposition of $F$}
We now decompose $F$ into appropriate terms of the same order. We write %As in \cite{Tsai}, we decompose $F$ into
\begin{equation} F = F_1 + F_2 + \cdots + F_5,\end{equation}
where
\begin{equation}\label{F.dec}
\begin{split}
& F_1 = \kappa  Q( 2|\zeta|^2+ \zeta^2 ), \\
& F_2 =  2 \kappa  QR b(2\zeta+\bar \zeta)+ 3 \kappa   QR^2 b^2 + \kappa  (\zeta+bR)^2 (\bar \zeta+bR), \quad \\
& F_3  =   2 \kappa  QR a^{(2)} (2\zeta+\bar \zeta), \quad F_4 = 2\kappa  Q[(\zeta +\bar \zeta) \eta +\zeta \bar\eta], \\
%& F_{5,1}  =2\kappa  Q[(\zeta +\bar \zeta) \eta^{(3)}+  \zeta \bar\eta^{(3)}]\\
& F_{5}  = \kappa  Q \bkt{ 2 |\eta_a|^2 + \eta_a^2} + 2 \kappa
QR b(
2\eta_a + \bar \eta_a) \\
&\quad \quad \quad + \kappa  (aR+ h)^2 (aR+\bar h)- \kappa
(\zeta+bR)^2 (\bar \zeta+bR),
\end{split}
\end{equation}
with $\eta_a = \eta + a^{(2)}R$. Note that $F_1$ consists of terms
of order $nz^2$; $F_2$, $F_3$ and $F_4$ consist of terms no smaller
than $n^2 z^3$; and $F_5$ higher order terms.
%===================================================================================================
\subsection{Basic estimates and normal forms}
In this subsection, we first give some basic estimates in Lemmas
\ref{Basic.lemma}, \ref{au.est.lemma}  and \ref{dotb.est.lemma}.
We then give the normal forms of the equations of $z_k$ and $b$ in
Lemmas \ref{L-NF} and \ref{b.NF}.
%=====================================================
\begin{lemma}[Basic Estimates] \label{Basic.lemma}
Suppose, for a fixed time, for some $\beta \ll n\le n_0$ and $p \geq 5$,
\begin{equation} \label{S-39}
\begin{split}
&\norm{Q}=n , \qquad %\norm{\eta}_{L^2 \loc \cap L^p} \le n,
\norm{\eta}_{L^2 \cap L^p} \ll 1, \qquad \norm{\eta}_{L^2 \loc} \le n, \\
&\max_{j \not = m}|z_j| \le \beta, \qquad
%\max_{j <m} |z_j| \leq z_L, \qquad
|a| \le C \beta^2.
\end{split}
\end{equation}
For all $1 \leq r \leq 2$, denote
\begin{equation}\label{XwtX.def}
\begin{split}
& X  :=  n \beta  \norm{\eta}_{L^2\loc} + n \norm{\eta}_{L^2 \loc} ^2 + \norm{\eta^3}_{L^1 \loc}, \\
& \wt X : = \beta^2  \norm{\eta}_{L^2\loc} + n \norm{\eta}_{L^2 \loc} ^2 + \norm{\eta^3}_{L^1 \loc}, \quad  Y(r,p) :=n\norm{\eta}^2_{L^p}  + \norm{\eta^3}_{L^{r}}.
\end{split}
\end{equation}
We have
\begin{equation}\label{basic.est}
\begin{split}
& \norm{F_5}_{L^1 \loc} \les n \beta^4 +  \wt X, \quad \norm{F_3+F_4+F_5}_{L^1 \loc }  \les n^2 \beta^3 + X, \\
& \norm{F-F_1}_{L^1 \loc} \les \beta^3 + X, \qquad \norm{F}_{L^1 \loc} \les n \beta^2  + X, \\
& |F_\theta| \les \beta^2 + n^{-1}X, \quad \norm{F -F_1}_{L^r }  \les \beta^3 + n\beta\norm{\eta}_{L^2\loc} + Y(r,p), \\
& \norm{F}_{L^r }  \les n \beta^2 + n\beta\norm{\eta}_{L^2\loc} + Y(r,p).
%& \norm{F-F_1}_{L^{5/4}}  \les \beta^3 + X + \wt{Y} , \quad  \norm{F}_{L^{5/4}}  \les n \beta^2 + X + \wt{Y}.
\end{split}
\end{equation}
\end{lemma}
%=======================================================
%\subsection{Proof of Lemma \ref{Basic.lemma}.}
\myproof First note that the proof of the first five estimates of
\eqref{basic.est} can be found in \cite[Lemma 3.2]{Tsai}. Although
\cite{Tsai} is for $m=0$ case, for $L^1_{loc}$ bounds the new
non-localized terms for $m>0$ are similarly estimated.

Now consider the last two $L^r$-estimates of \eqref{basic.est}.
The only non-localized terms of $F$ are of order $(u_k^{-}z_k)^3$,
$(u_k^-z_k)^2\eta$, $u_k^-z_k \eta^2$, and $\eta^3$ for $k< m$.
Since $|(u_k^-z_k)^2\eta|+|u_k^-z_k \eta^2|   \les |u_k^-z_k|^3 +
|\eta|^3$, we have
\begin{equation}\begin{split}
& \norm{(u_k^{-}z_k)^3}_{L^r}   + \norm{(u_k^-z_k)^2\eta}_{L^r}  + \norm{u_k^-z_k \eta^2}_{L^r}+ \norm{\eta^3}_{L^r} \\
\ &  \les |z_k|^3 \norm{u_k^-}_{L^{3r}}^{3} + \norm{\eta^3}_{L^r}
\les \beta^3 + \norm{\eta^3}_{L^r}.
\end{split}\end{equation}
Then for $1 \leq r \leq 2<p/2$, the estimates of $\norm{F}_{L^r}$
is the same as the estimates of $\norm{F}_{L^1_{loc}}$ except the
non-localized terms we just estimated and
\begin{equation}
\norm{Q z_k u_k^\pm \eta}_{L^r}  \lec n\beta
\norm{\eta}_{L^2_{loc}}, \quad \norm{Q\eta^2}_{L^r}  \lec n
\norm{\eta}_{L^p}^2.
\end{equation}
So, we obtain the last two estimates. \myendproof

\medskip
For some fixed $\frac{9}{2} < p <6$ which will be chosen, let us
define
\begin{equation} \label{zLzH.def}
  z_L = \bke{\tsum_{k =0}^{m-1} |z_k|^2}^{1/2}, \quad
  z_H = \bke{\tsum_{k=m+1}^K |z_k|^2}^{1/2}.
\end{equation}
If $m =0$, we set $z_L =0$. We also denote
\begin{equation} \label{hat.X.def}
\hat{X} = \hat{X}_p := n^4z_L\norm{\eta}_{L^p}^2 +  n^{6}z_L^2 \norm{\eta}_{L^p} + m \cdot n^{\frac{6(6-p)}{p}}\norm{\eta}_{L^p}^3.
\end{equation}
Note that if $m=0$, then $\hat{X} =0$. Let
\begin{equation} \label{D.def}
D = 6K c_{\max} \ga_0^+/\ga_0 = O(1)
\end{equation}
where $c_{\max} = \max_k 2 \int \phi_k^4$ and
\begin{equation} %\label{hat.X.def}
\ga_0^+ = \max_{k,l,m\in I, |s|<s_0} \lim_{r \to 0_+} \Im \bke{
\phi_k \phi_l \phi_m, \, \frac 1{H_0 +e_k - e_l - e_m - s -r i}
\Pc^{H_0}\phi_k \phi_l \phi_m }.
\end{equation}
Note that $(Q_{k,n}, R_{k,n})^{-1} =2\kappa  \int \phi_k^4 +
o(1)$. We have the following lemma on the normal forms of $z_k$.
%===============================================================================================================
\begin{lemma} \label{au.est.lemma}
Assume as in the Lemma \ref{Basic.lemma}, then for all $ k \not =m$,
we have
\begin{equation} \label{Z.est}
\begin{split}
|Z_k| & \les n\beta^2 +  \hat{X}_p + X, \quad \text{if} \quad k <m, \quad |Z_k|  \les n\beta^2 +  X, \quad \text{if} \quad k > m,\\
|R_k| & \les \beta^3  + \hat{X}_p +   X, \quad \text{if} \quad k <m, \quad |R_k| \les \beta^3  + X, \quad \text{if} \quad k > m.
\end{split}
\end{equation}
Here $Z_k$ is defined in \eqref{Z.def}  and  $R_k = R_{k,1} +
R_{k,2}$ is part of $Z_k$, where
\begin{equation} \label{Rk.def}
\begin{split}
R_{k,1} & := -2 c_k \big[ (u_k^+, F-F_1)+(u_k^-, \wbar F-\wbar F_1) \big ], \\
R_{k,2} & : = -2 c_k  \bkt{(u_k^+, h)+(u_k^-, \wbar h)+
(\bar{u}_k,R) a} F_ \theta .
\end{split}
\end{equation}
\end{lemma}
%========================================================================================
\myproof Recall \eqref{Z.def} that %$Z_j = Z_{j,1} + Z_{j,2}$ where
\begin{equation}\begin{split}
Z_{k} & : = -2 c_k \left \{(u_k^+, F)+(u_k^-, \wbar F) +
\bkt{(u_k^+, h)+(u_k^-, \wbar h)+ (\bar{u}_k,R) a} \dot
\theta\right \}.
\end{split}\end{equation}
For $m < k \leq K$, since $u_k^+, u_k^-$ are both real and localized, $P_k\eta =0$, using Lemma \ref{uv-orth} we have
\begin{equation} \label{Hlin.term.est}
|(u_k^+,\eta) + (u_k^-,\bar{\eta})| = 2|(u_k^-,\bar{\eta})| \leq Cn^2\norm{\eta}_{L^2\loc}.
\end{equation}
Therefore,
\begin{equation}\begin{split}
|Z_k| & \leq  \norm{F}_{L^1_{loc}} + |\dot \theta|[|a| + |z| + n^2\norm{\eta}_{L^2\loc}] \\
\ & \les n \beta^2 + X + [\beta^2 + n^{-1}X] (\beta + \norm{\eta}_{L^2_{loc}}) \les n\beta^2 + X.
\end{split}\end{equation}
Now, we consider the case when $k < m$. We first consider the term
$2c_k [(u_k^+, F)+(u_k^-, \wbar F)]$. As we already see in the
proof of Lemma \ref{Basic.lemma}, the only non-localized terms in
$F$ are bounded by  $|\eta^3|+\sum_{j,l,h < m} |u_j^- u_l^{-}
u_h^-| z_L^3$. Thus for $k < m$, using H\"{o}lder's inequality and
Lemma \ref{uv-est},
\begin{equation}\begin{split}
& |[(u_k^+, F)+(u_k^-, \wbar F)]| \les \norm{F}_{L^1_{loc}} +
(|u_k^-|, |\eta^3|+\tsum_{j,l,h < m} |u_j^- u_l^{-} u_h^-| z_L^3)
\\
 & \les n\beta^2  +  \hat{X}_p + X.
\end{split}\end{equation}
On the other hand, using \eqref{JPk-2}, we have
\begin{equation}\label{Llin.term.est}
|(u_k^+,\eta) + (u_k^-,\bar{\eta})| = |(\sigma_1\bar{\Phi}_k,
J[\eta])| \lec \norm{\eta}_{L^2\loc}, \quad (k<m).
\end{equation}
Then, it follows from Lemmas \ref{Basic.lemma} and \ref{uv-orth}
that
\begin{equation}
| [ (u_k^+, h) + (u_k^-, \bar{h}) + (\bar{u}_k, R)a] F_\theta|
\les [|z| + n^{-1}|a| + \norm{\eta}_{L^2_{loc}}] |F_\theta| \les
\beta^3 + X.
\end{equation}
This completes the proof of the
estimates of $Z_k$. By a similar way, we can obtain the estimates
of $R_k$. \myendproof

%====================================================================================
\begin{lemma}
\label{dotb.est.lemma} Assume as in the Lemma \ref{Basic.lemma},
then we have
\begin{equation} \label{dot.b.est}
|\dot{b}| \leq C[n^4z_L^2 + n\beta^3 + nX + n^2\beta\hat{X}].
\end{equation}
Above $\hat{X} = \hat{X}_p$ is defined in \eqref{hat.X.def} and can
be omitted if $m=0$.
\end{lemma}
%========================================================================================
\myproof
Recall \eqref{eq:b1} that
\begin{equation} \dot{b} = b_0 + c_m(Q, \Im (F-F_1 +\dot{\theta}h)) - A_{2,rm}. \end{equation}
It follows from \eqref{b0.def}, \eqref{A2rm.def} and Lemma \ref{au.est.lemma} that
\begin{equation} |b_0| \leq Cn^4 z_L^2 , \quad |A_{2,rm}| \leq n^2\beta [n\beta^2 + X + \hat{X}].  \end{equation}
On the other hand, we have
\begin{equation} |c_m(Q, \Im (F-F_1 +\dot{\theta}h))| \les n\norm{F-F_1}_{L^1\loc} + n^3\beta^2 + |\dot{\theta}|[n^3\beta + n\norm{\eta}_{L^2\loc}] \les n\beta^3 + nX. \end{equation}
So, \eqref{dot.b.est} follows.
\myendproof
%========================================================================================
\begin{lemma} \label{L-NF}
Fix $0 \leq m \le K$ and $0<n_1 \sim n \le n_0$. Let $Q=Q_{m,n_1}$
and $\L = \L_{m,n_1}$. Suppose $\psi$ is decomposed as in
\eqref{psi.dec} with respect to $\L$, and for some $0 < \beta \ll
n$
\begin{equation}\label{ass.nf}
\norm{\eta}_{L^2 \loc} \le \beta, \quad \norm{\eta}_{L^2 \cap L^p}
\ll 1, \quad  \max_{k \not = m } |z_k| \le \beta, \quad  |a| \le C
\beta^2.
\end{equation}
Then there exist functions $q_k$, $g_k, Y_k$ and constants $D_{kl}$
for $l \not = m$ such that
\begin{equation} \label{q-p.est}
\begin{split}
& \dot q_k -\Re(\LLLL_k)q_k =\tsum_{l {>m}} D_{kl} |q_l|^2 q_k + Y_k q_k + g_k, \quad \text{with} \quad |q_k - p_k| \les  n \beta^2, \\
& |D_{kl}| \leq Dn^2,  \quad \Re (D_{kl}) \le - \gamma_0  n^2, \quad \forall  k, l >m, \quad \text{and} \\
& |\Re (Y_{k})| \les n^2z_L^2, \quad(k>m); \quad |\Re (Y_{k})|
\les n^2\beta^2, \quad (k<m).
\end{split}
\end{equation}
Recall $\Re \LLLL_k\gec n^4$ if $k<m$ and $\Re \LLLL_k=0$ if
$k>m$. Moreover, we have
%\begin{equation} %\label{Yk.est}
% |\Re (Y_{k>m})| \les n^2z_L^2,  \quad \quad |\Re (Y_{k<m})| \les n^2\beta^2,
%\end{equation}
%and
\begin{equation} \label{gk.est}
\begin{split}
|g_{k}| & \lec n\beta^4 + n^4\beta z_L^2
+n^3\beta\norm{\eta}_{L^2\loc} +n\beta\norm{\eta^{(3)}}_{L^2\loc}
+ n\beta\hat{X}_p + \wt{X}, \quad(k>m),
\\
|g_{k}| & \lec n^5\beta^2 + n^4\beta z_L^2  + n\beta^4
+n^3\beta\norm{\eta}_{L^2\loc} +n\beta\norm{\eta^{(3)}}_{L^2\loc}
+ \hat{X}_p + \wt{X}, \quad(k<m).
\end{split}
\end{equation}
Above $\hat{X}_p$ is defined in \eqref{hat.X.def} and can be omitted
if $m=0$.

\end{lemma}
%==============================================================
\myproof In case $m=0$, Lemma \ref{L-NF} is identical to \cite[Lemma
3.4]{Tsai}. So, it suffices to assume $m>0$. The main difference in
case $m>0$ is that $u_l^-$ are not localized and $u_l^\pm$ are
complex for $l <m$. For those new terms involving $z_l$ with $l <m$,
we either integrate them using integration by parts and equations of
$r_l$, as in \eqref{A(2).int}, or include them in the error terms.
We sketch the proof here. Recall that for $k \in I_{m}$,
\begin{equation} \dot{z}_k - \LLLL_k z_k = Z_k, \quad \dot{p}_k - \Re(\LLLL_k)p_k = e^{i\omega_k t}Z_k, \quad \dot{r}_k = e^{-\LLLL_k t}Z_k, \end{equation}
where $Z_k$ is define in \eqref{Z.def}. For $F_j, j = 1,2,\cdots,
5$, defined in \eqref{F.dec}, we let
\begin{equation}
\begin{split}
& T_{k,1} : = -2c_ke^{i\omega_k t}[(u_k^+, F_4) + (u_k^-,
\bar{F}_4)], \quad
T_{k,2}: = -2c_k e^{i\omega_k t}[(u_k^+, F_1) + (u_k^-, \bar{F}_1)], \\
& T_{k,3} := -2c_k e^{i\omega_k t}\{(u_k^+, F_2 + F_3) + (u_k^-, \bar{F}_2 + \bar{F}_3) + [(u_k^+, \zeta) + (u_k^-,\zeta) \dot{\theta}]\}, \\
& T_{k,4}:= -2c_ke^{i\omega_k}\{(u_k^+, F_5) + (u_k^-, \bar{F}_5)
+ [(u_k^+, \eta) + (u_k^-,\bar{\eta})
+(\bar{u}_k,R)a]\dot{\theta}\}.
\end{split}
\end{equation}
Then, we can write
\begin{equation}
\dot{p}_k - \Re(\LLLL_k)p_k = e^{i\omega_k t}Z_k = T_{k,1} + T_{k,2} + T_{k,3} + T_{k,4}.
\end{equation}

The term $T_{k,1}$ contains the key terms with resonant
coefficients. Recall
\begin{equation}
F_ 4 = 2\kappa  Q \tsum_{j \not=m} \{(z_j \phi_j + \bar{z}_j
\phi_j)\eta + z_j \phi_j \bar{\eta} \} + O (n^3
|z|\norm{\eta}_{L^2\loc}),
\end{equation}
$[\eta] = \eta ^{(2)} + \eta ^{(3)}$, and denote $\eta^{(2)}_1$
the first sum for $\eta^{(2)}$ in \eqref{wt-eta2.dec}, which is
the main term of $\eta$. Since $\eta^{(2)}_1$ involves matrix
operators $R_{kl}$, we rewrite $T_{k,1}$ in vector form and get
$T_{k,1} = T_{k,1,1} +\ \text{error term}$, where
\begin{equation} \label{Tk11.def}
  T_{k,1,1} := -4\kappa  c_ke^{-i\omega_k t} \tsum_{j \not=m} [\wei{\svect{1 \\ -i}Q\phi_k\phi_j, \eta^{(2)}}(z_j +\bar{z}_j) + \wei{\svect{1 \\ i}Q\phi_k\phi_j, \eta^{(2)}} z_j]. \end{equation}
The error term is controlled by $n^3\beta\norm{\eta}_{L^2\loc} +
n\beta\norm{\eta^{(3)}}_{L^2\loc}$. By the explicit formula of
$\eta^{(2)}$ in \eqref{wt-eta2.dec}, we see that
\begin{equation} T_{k,1,1} = \sum_{l > m}D_{kl}|p_l|^2 p_k +\ \text{(non-zero phase\ $n^2z^3$-terms)} \ + Y_{k,1}p_k, \quad
\end{equation}
where $Y_{k,1}: =\sum_{l < m}D_{kl}|p_l|^2$. The first term in
$T_{k,1,1}$ we will keep. The middle term we integrate using
integration by parts. The last term is part of $Y_kp_k$ term. So,
we get
\begin{equation} \label{T1.eqn}
T_{k,1} = \frac{d}{dt}\wt{T}_{k,1} + \tsum_{l >m}D_{kl}|p_l|^2 p_k
+ Y_{k,1}p_k + g_{k,1}, \quad |\wt{T}_{k,1}| \les n^2\beta^3.
\end{equation}
Moreover, $Y_{k,1}$ and $g_{k,1}$ satisfy the estimates as those
of $Y_k$ and $g_k$ in the lemma. To compute $\Re (D_{kl})$, we use
\eqref{wt-eta2.dec}, \eqref{Tk11.def}, and Corollary
\ref{FRG-cor2}. The leading terms of $D_{kl}$ are from
$\wei{\svect{1 \\ -i} Q\phi_k \phi_j ,\eta^{(2)}_1}\bar{z}_j$.
%where $\eta^{(2)}_1$ is the first sum for $\eta^{(2)}$ in \eqref{wt-eta2.dec}.
We get
\begin{equation} %\label{Dkl.com}
\begin{split}
 \Re D_{kl} &= -C\Im (Q\phi_k\phi_l,\, (H_0 -E -i[-i(\om_l +\om_k) + 0^+])^{-1}P_c^{H_0}Q\phi_k\phi_l)  + O(n^4)\\
& = -C n^2\Im (\phi_m\phi_k\phi_l,\, (H_0 +
s-i0^+)^{-1}P_c^{H_0}\phi_m\phi_k\phi_l) +O(n^4),
\end{split}
\end{equation}
where $C=2\kappa ^2(2-\delta_{l}^k)\ge 2$, $s\in \R$ and $|s-(e_m
-e_l -e_k)| \lec n^2$. In particular $\Re D_{kl} =O(n^4)$ if $k<m$
or $l<m$.

For terms in $T_{k,2}$, we integrate as in \eqref{A(2).int}. Since
$e^{i\omega_k t}Z_k = T_{k,2} + R_k$ (see Lemma
\ref{au.est.lemma}), we get
\begin{equation} T_{k,2} = \frac{d}{dt}(\wt{T}_{k,2,1}) + T_{k,2,1} + O(n\beta\max_{l}|R_l|), \quad
|\wt{T}_{k,2,1}| \les n\beta^2.  \end{equation} The last term is
moved to the error term and it can be estimated by using Lemma
\ref{au.est.lemma}. The term $T_{k,2,1}$ are terms of $n^2 z^3$
and can be written as
\begin{equation}
T_{k,2,1} = Y_{k,2} p_k + \ \text{(non-zero phase\
$n^2z^3$-terms)}, \quad Y_{k,2}: = \tsum_{l \not =m}d_{kl}|p_l|^2.
\end{equation} The constants $d_{kl}$ can be computed explicitly.
In particular, $\Re(d_{kl}) = 0$ for $k, l >m$. So, we get $|\Re
(Y_{k,2})| \leq n^4 z_L^2$ if $k > m$. Again, the non-zero phase
terms are integrated by using integration by parts. So, we get
\begin{equation} \label{T2.eqn}
T_{k,2} = \frac{d}{dt} \wt{T}_{k,2} + Y_{k,2}p_k + g_{k,2}, \quad
|\wt{T}_{k,2}| \les n\beta^2.
\end{equation}
Moreover, $Y_{k,2}$ and $g_{k,2}$ satisfy the estimates as those of $Y_k$ and $g_k$ in the statement of the lemma.

For terms in $T_{k,3}$, we integrate terms which are not smaller than $n^2\beta^3$. The main difficulty is from terms with $\dot{\theta}$ since this term is not a polynomial expression in $z$ and $b$. The analysis is the same as in \cite{TY1, TY2, Tsai} with the use of Lemma \ref{au.est.lemma} and \eqref{dot.b.est}. We can write
\begin{equation} \label{T3.eqn}
T_{k,3} = \frac{d}{dt} \wt{T}_{k,3} + Y_{k,3}p_k + g_{k,3}, \quad |\wt{T}_{k,3}| \les \beta^3.
\end{equation}
Again, $g_{k,3}$ and $Y_{k,3}$ satisfy the estimates as those of
$g_k$ and $Y_k$ in the statement of the lemma. The only difference
between our case here and \cite{TY1, TY2, Tsai} is the computation
of $\Re(Y_{k,3})$. This can be done with simple calculation and
the attention that $u_l^\pm$ is complex only for $l<m$ and $\Im
u_{l}^{\pm} = O(n^2)$.

The term $T_{k,4}$ is part of the error term and  can be estimated
as $Z_k$ in Lemma \ref{au.est.lemma}:
\begin{equation} \label{T4.eqn}
|T_{k,4}| \les n\beta^4 + \wt{X}, \quad (k>m); \quad |T_{k,4}|
\les n\beta^4 + \wt{X} + \hat{X}_p, \quad (k<m).
\end{equation}
Now, let
\begin{equation} q_k := p_k - \wt{T}_{k,1} - \wt{T}_{k,2} - \wt{T}_{k,3}, \quad Y_k := Y_{k,1} + Y_{k,2} + Y_{k,3}.\end{equation}
Our lemma follows from \eqref{T1.eqn}, \eqref{T2.eqn}, \eqref{T3.eqn} and \eqref{T4.eqn}.
\myendproof
%==============================================================
\begin{lemma} \label{b.NF} Assume as in Lemma \ref{L-NF}. Then,
there exist functions $\wt{b}, g_b$ and numbers $B_{kl}$ for $k,
l\in \idx_{m}$ such that
\begin{equation}
\begin{split}
& \dot{\wt{b}} = b_0+\sum_{k,l \in \idx_{>m}}B_{kl}|z_k|^2|z_l|^2 + g_b, \quad |b - \wt{b}| \leq Cn\beta [\beta^2 + n\norm{\eta}_{L^2\loc}],\\
& |g_b| \leq C[n^3\beta^4 + n^5\beta z_L^2 +n^2\beta^2z_L^2 + n\beta^5 + n^2z_L\norm{\eta}_{L^2\loc} \\
& \quad \quad \quad  + n^2\norm{\eta}_{L^2\loc}^2 + n\norm{\eta^3}_{L^1\loc} + n\beta^2\norm{\eta^{(3)}}_{L^2\loc} + n\beta\hat{X}_p].
\end{split}
\end{equation}
Above $b_0$ is define in \eqref{b0.def} and can be omitted if $m=0$.
Moreover, we also have $|B_{kl}|\leq Cn^2$ and $B_{kl} =
-\frac{c_m}{2} \Re D_{kl} + O(n^4)$ where $D_{kl}$ is defined in
Lemma \ref{L-NF} and $c_m = (Q_m, R_m)^{-1} =O(1) >0$. Moreover,
$\max_{kl} (|B_{kl}|)/(K^{-1} \gamma_{0} n^2) \leq \frac{D}{2}$.
\end{lemma}
%==================================================================
\myproof Recall \eqref{eq:b1} that for $b_0$ defined in
\eqref{b0.def} and $A_{2,rm}$ defined in \eqref{A2rm.def},
\begin{equation} \label{bNF.eqn}
\dot{b} = b_0 +  c_m(Q, \Im(F-F_1 +\dot{\theta}h))
-A_{2,rm}.\end{equation} In case $m =0$ we have $z_L =0$, $b_0
=0$, $\hat{X}_p =0$, and Lemma \ref{b.NF} follows from Lemma
\cite[Lemma 3.4]{Tsai}. Therefore, we may assume $m >0$. All of
the work here is similar to that of \cite{TY1, TY2, Tsai}, so we
only give a sketch. Define
\begin{equation}
\begin{split}
& b_1: = c_m(Q, \Im(F-F_1)) - b_3 - b_4, \quad b_2 : = c_m(Q, \Im h)\dot{\theta}, \\
& b_3: = 2\kappa  c_m\Im(Q, Q[(\zeta + \bar{\zeta})\eta + \zeta
\bar{\eta}]), \quad b_4 = c_m\kappa  \Im (Q, 2|\zeta|^2\eta +
\zeta^2\bar{\eta}).
\end{split}
\end{equation}
Then, we have
\begin{equation} \dot{b} = \tsum_{j=0}^4 b_j - A_{2,rm}. \end{equation}

For terms in $b_1 + b_2$, we treat them as in \cite{TY1, TY2, Tsai}. Note that we can write
\begin{equation} b_2 + b_2 \sim nz^3 + nbz + n^2z^4 + bz^2 + n^{-1}b^2z +\ \text{error}. \end{equation}
The non-zero phase terms are of the form $C_{kl}|z_k|^2|z_l|^2$
and $C_kb |z_k|^2$. Here $|C_{kl}| \les n^2$ and $|C_k| \leq C$
for all $k, l \in I_m$. Moreover, $\Im C_{kl} = \Im C_k =0$ for
all $k, l \in I_{>m}$. Therefore, due to the $\Im$-operator, the
terms of the forms $C_{kl}|z_k|^2|z_l|^2$ with $C_kb |z_k|^2$ with
$k, l \in I_{>m}$ are killed. If one of $k, l$ is in $I_{<m}$, we
get new non-zero phase terms of the forms $C_{kl}|z_k|^2z_L^2$
with $|\Im C_{kl}| \les n^4$ and $C_Lb z_L^2$ with $|\Im C_L| \les
n^2$. We move these new terms into the error term. For all of the
non-zero phase terms, we integrate them by using integration by
parts. We get
\begin{equation} \label{ba12.eqn}
b_1 + b_2 = \frac{d}{dt}[\wt{b}_1 + \wt{b}_2] + g_{b,1} + g_{b,2},
\quad |\wt{b}_1| + |\wt{b}_2| \les n\beta^3,
\end{equation}
with $g_{b,1} + g_{b,2}$ satisfying the estimates as that of $g_b$
in the statement of the lemma.

For terms in $b_3$, we have $b_3 = 2\kappa  c_m\Im(Q^2, \zeta
\eta)$. We move terms $2\kappa  c_m\Im(Q^2, \zeta_k \eta)$ to
error term if $k<m$. So, we have $b_{3} = b_{3,1} + b_{3,2}$ with
\begin{equation}
b_{3,1} := \Im \tsum_{k {>m}} [\wei{f_{1k}, [\eta]}\bar{z}_k +
\wei{f_{2k}, [\eta]}\bar{z}_k], \quad |b_{3,2}| \les
n^2z_L\norm{\eta}_{L^2\loc}.
\end{equation}
Here, $f_{1k}, f_{2k}$ are some explicit localized functions of
order $n^2$. We need to integrate $b_{3,1}$ using equations
\eqref{wt-eta.eqn} of $\eta_\pm$ as in \cite{TY1}. We get
\begin{equation} \label{ba3.eqn}
b_3 = \frac{d}{dt}\wt{b}_3 + g_{b,3}, \quad |\wt{b}_{3}| \leq C[n^2\beta\norm{\eta}_{L^2\loc} + n^3\beta^3].
\end{equation}
Also, $g_{b,3}$ satisfies the estimates as that of $g_b$ in the
statement of the lemma.

For terms in $b_4$, as in \eqref{Tk11.def}, we write $b_{4} = b_{4,1} + b_{4,er}$ with
\begin{equation}
\begin{split}
& b_{4,1} = c_m\kappa  \Im \sum_{k,l \in I_m}[\wei{\svect{1 \\ -i}Q\phi_k\phi_l, \eta^{(2)}}z_k\bar{z}_l + \wei{\svect{1 \\ i}Q\phi_k\phi_l, \eta^{(2)}}z_kz_l], \\
& |b_{4,er}| \les n^3\beta^3\norm{\eta}_{L^2\loc} + n\beta^2 \norm{\eta^{(3)}}_{L^2\loc}.
\end{split} \end{equation}
From \eqref{wt-eta2.dec} and as in \eqref{T1.eqn}, we have
\begin{equation} \label{ba4.eqn}
b_4 = \frac{d}{dt}\wt{b}_4 + \sum_{k,l\in I_{>m}}B_{kl}|z_k|^2|z_l|^2 + g_{b,4}, \quad |\wt{b}_4| \leq Cn^2\beta^4.
\end{equation}
The term $g_{b,4}$ satisfies the estimate as that of $g_b$ in the statement of the lemma. The computation of $B_{kl}$ is exactly the same as that of $\Re D_{kl}$. The main leading terms come from
\begin{equation} c_m\kappa  \wei{\svect{1 \\ -i}Q\phi_k\phi_l, \eta^{(2)}_1}z_k\bar{z}_l,  \quad c_m \kappa \wei{\svect{1 \\ i}Q\phi_k\phi_l, \eta^{(2)}_1}z_kz_l.\end{equation}
Using \eqref{wt-eta2.dec}, Corollary \ref{FRG-cor2} and by direct computation, we get
\begin{equation} B_{kl} = -\frac{c_m}{2}\Re (D_{kl}) + O(n^4). \end{equation}

For terms in $A_{2,rm}$, we can write
\begin{equation} A_{2,rm} = n^3z^3 + n^2z^4 + n^2bz^2 + \text{error}. \end{equation}
Again, zero-phase terms are of the form $c_{lk}|z_k|^2|z_l|^2$ and $c_{1k}b|z_k|^2$. By direct computation, it follows that $\Im c_{kl} = \Im c_{1k} = 0$ if $k, l \in I_{>m}$. So, those terms $c_{lk}|z_k|^2|z_l|^2$ and $c_{1k}b|z_k|^2$ are killed by the $\Im$-operator if $k, l \in I_{>m}$. The other terms of the form $c_{lk}|z_k|^2|z_l|^2$ and $c_{1k}b|z_k|^2$ with $k$ or $l$ in $I_{<m}$ are moved to the error term. For terms of non-zero phase, we integrate them. We get
\begin{equation} \label{baArm.eqn}
A_{2,rm} = \frac{d}{dt}\wt{A}_{2,rm} + g_{b,5}, \quad |\wt{A}_{2,rm}| \les n^3\beta^3.
\end{equation}
The error term $g_{b,5}$ satisfies the estimate as that of $g_b$ in the statement of the lemma.
Finally, let
\begin{equation} \wt{b} = b-[\wt{b}_1 + \wt{b}_2 + \wt{b}_3 + \wt{b}_4 - \wt{A}_{2,rm}], \quad g_b = \sum_{j=1}^5 g_{b,j}. \end{equation}
Lemma \ref{b.NF} follows from \eqref{ba12.eqn}, \eqref{ba3.eqn}, \eqref{ba4.eqn} and \eqref{baArm.eqn}.
\myendproof
%==============================================================
%==================================================================

%==============================================================
\section{Converging to an excited state}
%==================================================================

In this and the next sections, we study the dynamics when the
solution is in a neighborhood of some excited states $Q_{1}$ at $
t=0$. We want to show that the solution either converges to an
excited state, or exits the neighborhood eventually. In the first
case, the ground state component is always bounded by other states.
In the second case, the ground state component becomes significant
after some time, denoted $t_{c}$ below. In this section we study the
dynamics for $t < t_c$. In next section we study the dynamics for
$t>t_c$ if $t_c$ is finite.

Denote $x_j(t)=(\phi_j,\psi(t))$ and $\xi(t)=\Pc^{H_0}\psi(t)$.
The assumption of Theorem \ref{mainthm} states that, at time
$t=0$,
\begin{equation}
\label{eq6.1}  |x_{1}(0)|= n, \quad \|\tsum_{j \not = 1}
x_j(0)\phi_j + \xi(0) \|_{H^1\cap L^1} \le  \rho_0, \quad \rho_0 =
n^{1+ \delta}.
\end{equation}
 Denote
\begin{equation}
\label{eq6.3} T_e := \sup_{T >0}\bket{T:\ \frac 1{\e_3}
\norm{\psi(t) -x_1(t)\phi_1}_{L^2} \leq |x_1(t)|  \in ((0.9)n,
(1.1)n), \quad 0 \leq \forall\ t \leq T }.
\end{equation}
Above $\e_3>0$ is the small constant in Lemma \ref{Linearized.dec}
and $T_e > 0$ by \eqref{eq6.1}. $T_e$ is the time the solution
exits the neighborhood of first excited state family. Note that
\eqref{eq6.1}--\eqref{eq6.3} are in terms of the orthogonal
coordinates. For the majority of this section we will use
linearized coordinates which depend on the choice of $Q$, but
\eqref{eq6.1}--\eqref{eq6.3} are independent of such a choice.

From Lemma \ref{Linearized.dec} and the definition of $T_e$, for
each $0 \leq T < T_e$, we can find a unique $ n(T)=n(\psi(T)) \in
(0,n_0)$ such that the solution $\psi(t)$ can be decomposed as
\begin{equation} \label{psi.best.dec}
 \psi(t) = [Q + a(t)R + \zeta(t) + \eta(t)]e^{-iEt +i\theta}, \
 \quad \forall \ 0 \leq t < T_e,
\end{equation}
with $a(T) =0$, where $Q = Q_{1,n(T)}$, $R = R_{1,n(T)}$ and $E =
E_{1,n(T)}$. The components $\zeta $ and $\eta $ are in the
corresponding spectral subspaces with respect to $Q_{1,n(T)}$.
Moreover we decompose
\begin{equation} \label{zeta.dec}
\zeta = \tsum _{j \not=1}\zeta_j, \quad \zeta_j = \bar{z}_ju_j^- +
z_j \bar u_j^+, \quad [\eta] = e^{i\theta}\eta_+ +
e^{-i\theta}\eta_{-}.
\end{equation}

Define
\begin{equation} \label{rho.def}
\rho(t)  := \frac{1}{n} (\Delta t + \ga_0 t)^{-1/2} , \quad \Delta t
:= (n\rho_0)^{-2}, \quad \rho(0)=\rho_0,
\end{equation}
where $\ga_0$ is given in \eqref{gamma0.def}, and let
\begin{equation}\label{tc.dec}
t_c : = %\max\big(n^{-3},\
\sup_{0 < T \leq T_e}  \{T : |z_0(t)| \leq
\e_4 n^{-1}\rho(t)^2,\ \ 0\leq t \leq T \},%\big),
\end{equation}
where $\e_4>0$ is a small constant to be chosen in
\eqref{L-error.est}, and $z_0$ is the coefficient of $\zeta_0$ in
\eqref{zeta.dec} with respect to $Q_{1,n(T)}$. If there does not
exist any $T$ satisfying the right side of \eqref{tc.dec}, we let
$t_c=0$.

Be definition $t_c \le T_e$ could be finite or infinite and is
independent of the choice of $Q$ in \eqref{psi.best.dec}. If it is
finite, it is the first time that $z_0$ becomes large enough, and
will not be destroyed by other components in the future. The
subscript $_c$ means ``change'' (of behavior). The function
$\rho(t)$ is an upper bound for higher bound states for $0 \le t \le
t_c$.

If $t_c=0$, we may skip most of this section and go
directly to Lemma \ref{OG.t-ch} and section \ref{S:5.2}.
%====================================================================
%\subsection{Converging to an excited state} \label{S:5.1}
%====================================================================

We will bound $\eta$ in $L^p$ and $L^2_{loc}$, with fixed $p$
satisfying
\begin{equation}
\frac {27}5 < p < 6, \quad \si =\si(p)= \frac{3p-9}{2p}, \quad
\frac 23 < \si < \frac 34.
\end{equation}

From now on let $0\le T<t_c$ and $\psi$ be decomposed as in
\eqref{psi.best.dec} with respect to $Q_{1,n(T)}$. We start with
the following lemma.
%================================
%===================================
\begin{lemma} [Initial estimates]\label{OG.eta}
Fix $ \frac {27}5 < p < 6$ with $\sigma(p) = \frac{3p-9}{2p}$. We
have
\begin{equation}
\sum_{k \not = 1}|z_k(0)|^2 \le \tfrac98 \rho_0^2,\quad
\norm{e^{\bL t} \eta_\pm(0)}_{L^p} \wei{t}^{\sigma(p)} +
\norm{e^{\bL t} \eta_\pm(0)}_{L^2\loc} \wei{t}^{7/6}\le C_2 \rho_0
\end{equation}
for $t \ge 0$, for  some $C_2>0$ uniformly in $n=n(T)$, $0\le
T<T_e$.
\end{lemma}

\myproof Let $ \psi':= e^{-i\theta(0)} \psi(0) -Q$. From
\eqref{psi.best.dec} at $t=0$, we have
\begin{equation}
a(0)R + \zeta(0) + \eta(0) = \psi'= e^{-i\theta(0)} \bke{\tsum
_{j=0}^K x_j(0)\phi_j + \xi(0)}  -Q.
\end{equation}
For $k \not=1$, applying the projection $P_k$ on this equation, we
get
\begin{equation}
|z_k(0)| \leq |2 c_k|[|(u_k^+, \psi')| + |(u_k^-,
\overline{\psi'})|] \leq (1+o(1))[|x_k(0)| + n^3] .
\end{equation}
Thus $\sum_{k \not = 1}|z_k(0)|^2 \leq \frac98 \rho_0^2$ by
\eqref{eq6.1}. Moreover, since $\psi'$ is localized and
$\norm{\psi'}_{H^1\cap L^1} \lec \rho_0$, using Lemma
\ref{Pipmest}, we get the estimates of $\eta_\pm(0)$ for $t>1$ by
Lemma \ref{th:decay} and for $0\le t\le 1$ by Lemma
\ref{H1conserve}. \myendproof

Recall $\eta^{(3)}$ and $z_H$ are defined in \eqref{wt-eta3.def}
and \eqref{zLzH.def}.
%\begin{definition}
We now define
\begin{equation} \label{MT.def}
M_T := \sup _{0 \le t \le T} \max \bket{
\begin{aligned}
& \rho(t)^{-1} z_H(t),\quad 2D^{-1}\rho^{-2}(t)|a(t)|,\\
& \bkt{n^{2\si-1} \rho(t)^{2\si-2\al}
+ 2C_2\rho_0\wei{t}^{-\si(p)}}^{-1}\norm{\eta(t)}_{L^p},\\
& \bkt{n^{-\al/2} \rho^{3} + n^{4/5} \rho^{7/3} + 2C_2
\rho_0\wei{t}^{-7/6}}^{-1} \|\eta^{(3)} (t) \|_{L^2 \loc}
\end{aligned}
}.
\end{equation}
%\end{definition}
Above $ \al >0$ is a small constant to be chosen. We can choose
$\al =0.01$.

Clearly $M_0 \le 3/2$ if $n$ is sufficiently small. By continuity
we have $ M_T \le 2$ for $T>0$ sufficiently small. Our main result
in this section is the following proposition, which implies $M_T
\le 3/2$ for all $T<t_c$ by a continuity argument.

\begin{proposition} \label{th:S-4}
Suppose that for some $T \in [0 , t_{c})$, $M_T$ is well-defined
and $M_T \leq 2$.  Then we have $M_T\le 3/2$ and $n(T)/n \in
(\frac 34, \frac 54)$.
\end{proposition}

The proof of Proposition \ref{th:S-4} is decomposed to Lemmas
\ref{wtXX.r}--\ref{bsz.est}.

Note that $T < t_c$ and $M_T\le 2$ imply
\begin{equation} \label{MT2}
\begin{split}
& |z_0(t)| \leq \e_4 n^{-1} \rho^2(t), \quad z_H(t) \le 2 \rho(t)
, \quad
|a(t)|\le D \rho(t)^2, \\
& \norm{\eta(t)}_{L^p} \le 2 n^{2\si-1} \rho(t)^{2\si-2\al}
+ 4C_2\rho_0\wei{t}^{-\si}, \\
& \| \eta^{(3)}(t) \|_{L^2 \loc} \le 2n^{-\al/2} \rho^{3} +
2n^{4/5} \rho^{7/3} + 4C_2 \rho_0\wei{t}^{-7/6}.
\end{split}
\end{equation}
Since $[\eta]=\eta^{(2)}+\eta^{(3)}$ and $\norm{\eta^{(2)}}_{L^2
\loc}\lec n \rho^2 $ by its definition, we get
\begin{equation} \label{zall-eta.est}
 \norm{\eta(t)}_{L^2 \loc} \les n\rho(t)^2 +\rho_0\wei{t}^{-7/6}.
\end{equation}
It is sometimes convenient to use
\begin{equation} \label{rho0rho}
\rho_0 \bka{t}^{-1/2} \lec \rho(t) \lec n^{-1}\bka{t}^{-1/2},
\quad \norm{\eta}_{L^p} +\norm{\eta(t)}_{L^2 \loc} \lec \rho.
\end{equation}

%===============================================================================================================
\begin{lemma} \label{wtXX.r}
Recall $X$, $\td X$, $F$ and $F_1$ are defined in \eqref{XwtX.def},
\eqref{F.def}, and \eqref{F.dec}, with $\frac {27}5<p<6$. Assume
$M_T \leq 2$, then we have
\begin{equation} \label{X-wtX.est}
\begin{split}
\wt{X} &\les n\rho^4 +  \rho_0 \rho(t)^2 \wei{t}^{-7/6}
+ n\rho_0^2 \wei{t}^{-7/3}, \\
X &\les n^2\rho^3  + n\rho_0\rho(t)\wei{t}^{-7/6} + n\rho_0^2
\wei{t}^{-7/3},
\end{split}
\end{equation}
and, with $o(1)$ denoting small positive constants which go to $0$
as $n+\norm{\psi_0}_{H^1} \to 0$,
\begin{equation} \label{FN1.est}
\begin{split}
\norm{F}_{L^{p'}} &\les n\rho^2 + o(1)\rho_0^{2} \wei{t}^{-1.4}, \\
\norm{F-F_1}_{L^{\frac{9}{8}} \cap L^{\frac 32}} &\les \rho^3
 + n^{0.64} \rho^{2.54} + \rho_0^{7/4} \bka{t}^{-5/4}.
\end{split}
\end{equation}
\end{lemma}
%===============================================================================================================
\myproof By H\"{o}lder's inequality for $p\ge 9/2$, and
$\norm{\eta}_{L^2\cap L^p} \ll 1$,
\begin{equation} \label{eta.m3}
\begin{split}
& \norm{\eta^3}_{L^1_{loc}}  \leq
\norm{\eta}_{L^2_{loc}}^{\frac{2p-6}{p-2}}
\norm{\eta}_{L^p}^{\frac{p}{p-2}}, \quad
\norm{\eta^3}_{L^1}  \leq
%\norm{\eta}_{L^2}^{\frac{2p-6}{p-2}}
o(1)\norm{\eta}_{L^p}^{\frac{p}{p-2}} ,\\
& \norm{\eta^3}_{L^{p'}} \leq
o(1)%\norm{\eta}_{L^2}^{\frac{2(p-4)}{p-2}}
\norm{\eta}_{L^p}^{\frac{p+2}{p-2}}, \quad \norm{\eta^3}_{L^{9/8}
\cap L^{3/2}} \leq o(1)%\norm{\eta}_{L^2}^{\frac{2(2p-9)}{2(p-2)}}
\norm{\eta}_{L^p}^{\frac{11p}{9(p-2)}}.
\end{split}
\end{equation}
From  \eqref{XwtX.def} with $\be=\rho$ and $n$ replaced by
$n(T)\sim n$,
\begin{equation}
\wt{X} \lec \rho^2\norm{\eta}_{L^2_{loc}} +X_1, \quad X \les n
\rho\norm{\eta}_{L^2_{loc}} + X_1, \quad   X_1=
n\norm{\eta}_{L^2_{loc}}^2 + \norm{\eta^3}_{L^1_{loc}}.
\end{equation}
Using \eqref{MT2}$_2$, \eqref{zall-eta.est}, and \eqref{eta.m3}$_1$,
one gets for $\frac{27}5<p<6$ that
\begin{equation}
\label{X1.est} X_1 \lec n^2 \rho^4  + \rho_0 \rho^2
\bka{t}^{-7/6}+ n \rho_0^2 \bka{t}^{-7/3}.
\end{equation}
One gets \eqref{X-wtX.est} from the above two equations.

To bound $F= \kappa  Q(2|h_\si|^2 + h_\si^2) + \kappa  |h_\si|^2
h_\si$ in $L^{p'}$ with $h_\si = aR + \zeta + \eta$, since
$\norm{aR} \lec n^{-1}\rho^2$, $\norm{\zeta}_{L^q} \lec \rho$ for
$q \ge 2$, and $\norm{\eta}_{L^p} \le \rho$, by \eqref{eta.m3}$_2$
and \eqref{MT2} we get
\begin{equation} \label{Ftc.est}
\norm{F}_{L^{p'}} \lec n \rho^2 +  o(1)\norm{\eta}_{L^{p}}^{\frac
{p+2}{p-2}} \lec n \rho^2 + o(1)\rho_0^{2} \bka{t}^{-1.4510}.
\end{equation}
Similarly, to bound $F-F_1$ with $F_1 =\kappa  Q(2|\zeta|^2 +
\zeta^2)$, by \eqref{eta.m3} we have
\begin{equation}
\norm{F-F_1}_{L^{\frac{9}{8}} \cap L^{\frac 32}} \lec  \rho^3 +  n
\rho \norm{\eta}_{L^2_{loc}}
+o(1)\norm{\eta}_{L^{p}}^{\frac{11p}{9(p-2)}}.
\end{equation}
By \eqref{MT2}, $\rho \le n^{-1}\bka{t}^{-1/2}$, and
$\frac{27}5<p<6$, it is bounded by
\begin{equation}
\begin{split}
&\lec \rho^3 + n\rho[n \rho^{2} + \rho_0\wei{t}^{-7/6}] + [
n^{0.6471} \rho^{2.5494} + \rho_0^{1.8333} \bka{t}^{-1.2941}] \\
& \lec  \rho^3 + n^{0.64} \rho^{2.54} + \rho_0^{7/4} \bka{t}^{-5/4}.
\end{split}\end{equation} \myendproof

\begin{lemma}[Dispersion estimates] \label{disp.eta.est}
Assume $M_T \leq 2$, then for all $0 \leq t \leq T$, we have
\begin{equation}\begin{split}
\norm{\eta(t)}_{L^p} &\leq \frac{3}{2} n^{2\si-1}
\rho(t)^{2\si-2\al}
+3C_2\rho_0\wei{t}^{-\si}, \\
\|\eta^{(3)}(t) \|_{L^2_{loc}} &\leq \frac{3}{2} [
 n^{-\al} \rho^{3} + n^{4/5} \rho^{7/3}] + 3C_2\rho_0\wei{t}^{-7/6} .
\end{split}
\end{equation}
\end{lemma}
%=======================================
\myproof We first prove the $L^p$-bound. Since $[\eta] =
e^{i\theta}\eta_+ + e^{-i\theta}\eta_{-}$,  it suffices to
estimate $\norm{\eta_\pm}_{L^p}$. By \eqref{eta-pm.eqn} with
$t_0=0$, and by Lemmas \ref{th:decay} and \ref{Pipmest},
\begin{equation} \label{L5.etapm.est}
\norm{\eta_\pm}_{L^p} \lec \norm{e^{t\bL}\eta_\pm(0)}_{L^p} +
\int_0^t\al_p(t-s)[\norm{F_{L\pm}}_{L^{p'}} +
\norm{F}_{L^{p'}}](s) ds.
\end{equation}
By Lemma \ref{OG.eta},
\begin{equation} \label{L52.ini}
%\norm{e^{t\bL} \eta_\pm(0)}_{L^2\loc} \leq C_2
%\rho_0\wei{t}^{-7/6}, \quad
\norm{e^{t\bL} \eta_\pm(0)}_{L^p} \leq C_2 \rho_0\wei{t}^{-\si}.
\end{equation}
By \eqref{basic.est}, Lemma \ref{wtXX.r}, and \eqref{rho0rho},
\begin{equation} \label{F.theta.est}
|\dot{\theta}| = |F_\theta| \les \rho^2 + n^{-1} X  \les \rho(t)^2
+ \rho_0\rho(t)\wei{t}^{-7/6}  + \rho_0^2 \wei{t}^{-7/3} \les
\rho(t)^2.
\end{equation}
By \eqref{FL}, \eqref{rho0rho}, and Lemma \ref{Kpmest},
\begin{equation} \label{FL1.est}
\norm{F_{L\pm}}_{L^{p'}} \les |F_\theta|(\norm{\eta}_{L^p} +
n^{-1}|a| + |z|) \les \rho^2\cdot \rho = \rho^3.
\end{equation}
By Lemma \ref{wtXX.r},
%\begin{equation} \label{FN1.est}
$\norm{F}_{L^{p'}} \lec n\rho^2  + \rho_0^{2}\wei{t}^{-7/5}$.
%\end{equation}
Thus the integral in \eqref{L5.etapm.est} is bounded by
\begin{equation} \label{L5eta-p1.est}
\les  \int_{0}^t \alpha_p(t-s) [n\rho^2(s)  +
\rho_0^{2}\wei{s}^{-7/5}] ds \les \rho_0^{2\al} n^{2\si-1}
\rho(t)^{2\si -2\al} + \rho_0^{2}\wei{t}^{-\si}.
\end{equation}
Here we have used \eqref{rho.def}, $n\rho^2(s) \sim n^{-1}(\De
t+s)^{-1}$, and $\forall 0<\al < \si<1$
\begin{equation}
\label{int-est1} \int_{0}^t |t-s|^{-\si}(\Delta t + s)^{-1} ds
\les (\Delta t)^{-\al }(\Delta t + t)^{-\si +\al}.
\end{equation}
Combining \eqref{L52.ini} and \eqref{L5eta-p1.est}, we get the first
estimate of Lemma \ref{disp.eta.est}.

We next prove the second estimate. Recall that $\eta_{\pm}^{(3)} =
\sum_{j=1}^4\eta_{\pm,j}^{(3)}$, where $\eta_{\pm,j}^{(3)}$ are
defined in \eqref{wt-eta14} and \eqref{wt-eta23.def} with $t_0=0$.
By Lemmas \ref{OG.eta} and \ref{th:sdecay}, we get
\begin{equation} \label{eta(3).12.est}
\norm{\eta_{\pm,1}^{(3)}}_{L^2_{loc}} \leq C_2\rho_0\wei{t}^{-7/6},
\quad \norm{\eta_{\pm,2}^{(3)}}_{L^2_{loc}}  \leq Cn
\rho_0^2\wei{t}^{-3/2}.
\end{equation}

For $\eta_{\pm,3}$, by Lemma \ref{au.est.lemma}, \eqref{rho0rho},
and \eqref{X-wtX.est},
\begin{equation}%\label{fkl.est}
\max |p_k| \lec n\rho^2+\hat X_p+X \lec n \rho^2.
\end{equation}
By \eqref{fkl-th.est}, \eqref{F.theta.est} and the above,
\begin{equation}%\label{fkl.est}
\||\dot f_{kl}| +|\dot{\theta}f_{kl}|\|_{L^2_r} \lec n |\dot
\th|\rho^2 + n \rho \max |\dot p_k| \lec  n \rho^2 \rho^2 + n \rho
(n\rho^2) \lec n^2 \rho^3.
\end{equation}
It follows from Lemma \ref{th:sdecay} that
\begin{equation} \label{eta(3).3.est}
\norm{\eta_{\pm,3}^{(3)}}_{L^2\loc} \leq C\int_0^t\wei{t-s}^{-3/2}
n^2 \rho^3(s) ds \leq Cn^2\rho^3(t).
\end{equation}
Here we have used, for $a,b>1$ and $S \ge 1$,
\begin{equation} \label{int-est2}
\int_{0}^t\bka{t-s}^{-a} (S+s)^{-b}ds \lec  S^{1-b} (S+t)^{-a} +
(S+t)^{-b},
\end{equation}
which is bounded by $(S+t)^{-b}$ if $a \ge b$.

For $\eta_{\pm,4}$, by Lemma \ref{th:decay}, we have
\begin{equation} \label{eta3-4}
\norm{\eta_{\pm,4}^{(3)}}_{L^2\loc} \leq C\int_0^t \al_\infty(t-s)
[\norm{F_{L\pm}}_{L^{9/8}\cap L^{3/2}} + \norm{F-F_1}_{L^{9/8}\cap
L^{3/2}}](s) ds,
\end{equation}
where $\al_\I(t) = t^{-1/2} \bka{t}^{-2/3}$.  It follows from
\eqref{int-est2} that
\begin{equation}%\label{int-est2}
\int_{0}^t\al_\infty(t-s) \rho(s)^{r} ds \lec \rho(t)^{r}+
n^{1/3}\rho_0^{r-2} \rho(t)^{7/3}  , \quad r >2.
\end{equation}
As for \eqref{FL1.est}, we have $\norm{F_{L\pm}}_{L^{9/8}\cap
L^{3/2}} \lec \rho^3$. By Lemma \ref{wtXX.r}, $\norm{F-F_1}_{L^{9/8}
\cap L^{3/2}}   \les \rho^3 + n^{0.64} \rho^{2.54} + \rho_0^{7/4}
\bka{t}^{-5/4}$. Thus
\begin{equation}\label{FN2.term}
\begin{split}
\norm{\eta_{\pm,4}^{(3)}}_{L^2\loc} & \lec (\rho^3 +
n^{1/3}\rho_0\rho^{7/3}) + (n^{0.64} \rho^{2.54}+ n^{0.97}
\rho_0^{0.54}\rho^{7/3})+ \rho_0^{7/4} \bka{t}^{-5/4}
\\
&\lec  \rho^3 + o(1) n^{4/5} \rho^{7/3} + \rho_0^{7/4}\wei{t}^{-5/4}
.
\end{split}
\end{equation}

Summing \eqref{eta(3).12.est}, \eqref{eta(3).3.est} and
\eqref{FN2.term}, we get the bound of $ \|\eta_{\pm}^{(3)} \|_{L^2
\loc} $ in the lemma. \myendproof

%================================================================================================
\begin{lemma}[Bound states estimates] \label{bsz.est}
Assume $M_T\le 2$, then for all $0 \leq t \leq T$, we have
\begin{equation}z_H(t) \leq \frac{3}{2} \rho(t),  \quad
|a(t)| \leq \frac{3}{4} D \rho(t)^2, \quad |n(t) - n| \le \frac 14
n.
\end{equation}
\end{lemma}
%================================================
\myproof For $1 < k \leq K$, from Lemma \ref{L-NF}, we have a perturbation $q_k$ of $p_k$ such that
\begin{equation} \label{NF-ZH}
\dot q_k = \sum_{l \not =1} D_{kl} |q_l|^2 q_k + Y_k q_k + g_k,
\end{equation}
where
\begin{equation}
\begin{split}
& |q_k- p_k| \les Cn\rho^2, \quad |\Re(Y_k)| \leq Cn^2z_L^2 \leq C\rho^4(t),\\
& |g_k| \les  n \rho^4 + n^3\rho\norm{\eta}_{L^2\loc}  + n \rho \|
\eta^{(3)} \|_{L^2_{loc}} + \wt{X} + n\rho \hat{X}.
\end{split}
\end{equation}
From \eqref{hat.X.def} and $\norm{\eta}_{L^p} \le \rho$, we have
$\hat{X} \lec   \rho^3$.  Thus, from \eqref{MT2},
\eqref{zall-eta.est} and Lemma \ref{wtXX.r}, we get
\begin{equation}\begin{split}
|g_k| & \les o(1) n^2\rho^3 + n\rho_0\rho\wei{t}^{-7/6} +
n\rho_0^2\wei{t}^{-7/3}.
\end{split}\end{equation}
Since $\rho_0 = n^{1+\de}$ and $0< \delta < \frac{3}{2}$, it follows that
\begin{equation}\label{gH.a.est}
\int _0^{n^{-3} \mmin T}|g_k|(t) dt \le C n \rho_0;\quad |g_k|(t)
\leq o(1) n^2\rho^3(t), \quad \forall t \geq n^{-3}.
\end{equation}
Now, from \eqref{NF-ZH}, we get
\begin{equation}\label{NF-ZH2}
\frac{d}{dt}|q_k| = \sum_{l\not = 1} \Re (D_{kl}) |q_l|^2 |q_k| + (\Re
Y_k) |q_k| + \Re (\frac {\bar q_k}{|q_k|} g_k).
\end{equation}
for all $ 0\leq t \leq n^{-3}$, by integrating this equation on
$(0,t)$, we see that $|q_k(t) - q_k(0)| \ll \rho_0$ . Using $z_H =
(\tsum_{k >1} |p_k|^2)^{1/2}$, $z_H(0) \le \sqrt{9/8}\rho_0$ and
$|q_k - p_k|\lec n \rho^2$, we get
\begin{equation}\label{zk.in.layer}
z_H(t) \leq 1.1 \rho_0, \quad \forall\ 0 \leq t \leq n^{-3}.
\end{equation}
Now, let $f_H = (|q_{2}|^2 + \cdots +|q_{K}|^2)^{1/2}$, from
\eqref{NF-ZH2} and \eqref{q-p.est}, in particular $D_{k0} |q_0|^2
\lec n^2(n^{-1}\rho^2)^2=\rho^4$, we get
\begin{equation}\begin{split}
\dot{f}_H & \leq -\frac{\gamma_0 n^2}2 f_H^3 + C[f_H\rho^4 +
\tsum_{k=2}^K|g_k|].
\end{split}\end{equation}
By \eqref{MT2} and \eqref{gH.a.est}, we get
\begin{equation} \label{fH.eqn}
\dot{f}_H \leq -\frac{\gamma_0 n^2}2 f_H^{3} + o(1) n^2 \rho(t)^3,
\quad n^{-3} \leq  t \leq t_c.
\end{equation}
Let $g(t):=\frac{7}{5} \rho(t)$. We have $f_H(n^{-3}) < g(n^{-3})$
and $\dot g =- \frac {\ga_0 n^2}{2}\frac{25}{49} g^3$, thus $\dot
f_H(t) < \dot g(t)$ if $f_H(t)=g(t)$. By comparison principle,
\begin{equation} %\label{fH.eqn}
f_H(t) \leq g(t)=\frac{7}{5} \rho(t), \quad (n^{-3} \leq t \leq T),
\end{equation}
which together with \eqref{zk.in.layer} give the first estimate of
the Lemma.

For the second estimate, recall that $a = a^{(2)} + b$ with
$|a^{(2)}|\leq Cn^2 \rho^2(t)$. From Lemma \ref{b.NF}, there is a
perturbation $\wt{b}$ such that
\begin{equation} \label{bper.eqn}
\frac{d}{dt} \wt{b} = b_0 + \hat b_0 + \sum_{1 < l, k \leq K} B_{kl}
|z_l|^2 |z_k|^2 + g_b ,
\end{equation}
where $g_b$ and $B_{kl}$ are defined in Lemma \ref{b.NF} and $\hat
b_0 = B_{00}|z_0|^4 + 2\sum_{1 < k \leq K} B_{k0} |z_0|^2 |z_k|^2$.
We have $|b -\wt{b}| \leq C n^2\rho^2$ and $|b_0|+|\hat b_0| \les
n^4|z_0|^2 \lec \e_4^2 n^2 \rho^4$. By Lemma \ref{b.NF},
\eqref{MT2}, (in particular $|z_0|\le \e_4 n^{-1}\rho^2$ and this is
where we choose $\e_4$), \eqref{zall-eta.est}, Lemma \ref{wtXX.r},
\eqref{X1.est} and $\hat{X} \lec  n^4 \rho^3+\norm{\eta}_{L^p}^3$,
\begin{equation} \label{L-error.est}
\begin{split}
|g_b| & \les n^3\rho^4 + n\rho^5 + \e_4 n\rho^2 \| \eta
\|_{L^2_{loc}} +n\rho^2 \| \eta^{(3)} \|_{L^2_{loc}} + nX_1
+ n\rho \hat{X} \\
& \lec o(1) n^2\rho(t)^4 + \wt{g}_b, \qquad \wt{g}_b = n^{2}\rho_0^2
\bka{t}^{-7/3}+n\rho_0 \rho^2 \bka{t}^{-7/6}.
\end{split}
\end{equation}
Then, for $t \geq \De t = n^{-2}\rho_0^{-2}$, we have $\rho(t) \sim
n^{-1} t^{-1/2}$ and
\begin{equation}
\int_t^{T} |\wt g_b|(s)ds \lec \int_t^{\infty}[n^4 s^{-7/3} +
s^{-7/6-1}] ds \lec n^4 t^{-4/3} + t^{-7/6}   \lec n^2\rho(t)^2.
\end{equation}
For $0 \leq t \leq \De t $, we have $\rho(t) \sim \rho_0$ and
\begin{equation}
\int_t^{T} |\wt g_b|(s)ds  \leq  \bke{ \int_t^{\De t }+\int_{\De
t}^{\I}} |\wt g_b|(s)ds  \lec \int_t^{\De t }n^{2}\rho_0^2
\bka{s}^{-7/6}ds + n^2 \rho_0^2 \lec n^2 \rho_0^2.
\end{equation}
Using $\int_t^\I n^2\rho^4 ds \lec \rho(t)^2$, we get
have
\begin{equation} \label{gb.int}
 \int_t^{T} |b_0 + g_b|(s)ds \leq o(1)\rho(t)^2, \ \forall \ t \in [0, T).
\end{equation}
Integrating \eqref{bper.eqn} on $(t ,T)$ and using $\max_{kl}
(|B_{kl}|)/(K^{-1} \gamma_{0} n^2) \leq \frac{D}{2}$, we get
\begin{equation}
|\wt{b}(t)|  \leq |\wt{b}(T)|+ \frac{D}{2} \rho^2(t) + o(1)\rho^2(t)
\leq |\wt{b}(T)|+ \frac 59 D \rho^2(t).
\end{equation}
Now, since $a(T) =0$, we get
\begin{equation}
|\wt{b}(T)| = |a(T) - b(T)| + |b(T) - \wt{b}(T)| \leq |a^{(2)}(T)| +
Cn^2\rho(T)^2 \lec n^2 \rho(t)^2.
\end{equation}
Thus we have $|\wt{b}(t)| \leq |\wt{b}(T)| +  |\wt{b}(t)- \wt{b}(T)|
\leq \frac 58 D\rho(t)^2$ and
\begin{equation}
|a(t)| \leq |a^{(2)}(t)| + |\wt{b}(t)| + |\wt{b}(t) -b(t)| \leq
\frac 34 D\rho(t)^2.
\end{equation}
Finally, Lemma \ref{Linearized.dec} shows $|n(T) - n(t)| \lec n^{-1}
|a(t)|+ n^3 \ll n$ and the last claim of the Lemma. \myendproof
%===============================================================================================================

The proof of  Lemma \ref{disp.eta.est} and Lemma \ref{bsz.est}
complete the proof of Proposition \ref{th:S-4}.

We now distinguish the two cases that $t_c = \I$ and $t_c<\I$.

Suppose $t_c = \I$. By Lemma \ref{Linearized.dec} (iii) we have for
any $t < T < \I$
\begin{equation}
|n(t)^2 - n(T)^2| \lec   |a_{n(T)}(t)| \lec \rho^2(t),
\end{equation}
which shows that $n(t)$ converges to some $ n_\I\sim n$ as $t \to
\I$. Furthermore $n(t) \sim n(0) \sim n_\I$ and $|n(t) - n_\I|\lec
n^{-1}\rho^2(t)$. Together with the estimate $M_T\le 3/2$ we have
shown the main theorem in the case the solution converges to an
excited state.

In the case $t_c < \I$, by continuity we also have $M_{t_c} \le
3/2$. we will show that the solution escapes from the  first excited
state family in the next section. We prepare it with the
following lemma, whose proof is the same as that for $\eta_\pm(t)$
in Lemma \ref{disp.eta.est} with the nonlinear
terms set to zero for $t_c < s< t$.
%================================================================
\begin{lemma} \label{OG.t-ch}
Suppose $t_c <\infty$. Let $\Delta t = n^{-2}\rho_0^{-2}$ and
$\eta_\pm(t) =e^{\mp i \th(t)} P_\pm [\eta(t)]$ where $\eta(t)$ is
as in \eqref{psi.best.dec} with respect to $Q_{1,n(t_c)}$. Then for
all $t \geq t_c$, we have
\begin{equation}\label{eq:4.64}
\norm{e^{\bL (t-t_{c})} \eta_\pm (t_{c})}_{L^p} \leq \frac 14
\La_1(t) ,\quad \norm{e^{\bL(t-t_c)} \eta_\pm(t_c)}_{L^2_{\loc}}
\leq\frac 14 \La_2(t),
\end{equation}
where for $C_2$ from Lemma \ref{OG.eta}, some $C_3>0$ and $\rho_c
= \rho(t_c)$,
\begin{equation}\label{eq:4.65}
\begin{split}
\La_1(t) &= C_3[C_2 \rho_0\wei{t}^{-\sigma(p)}
%+ n^{-1}(\De t)^{-\alpha}[\Delta t + t]^{-\sigma(p) + \alpha}],
+ n^{2\si-1}\rho_0^{2 \al} \rho(t)^{2\si-2\al}],
\\
\La_2(t) &=C_3[C_2 \rho_0\wei{t}^{-7/6} +
n\rho_c^2\wei{t-t_c}^{-7/6}  + \rho^{3}(t) +
n^{4/5}\rho^{7/3}(t)].
\end{split}
\end{equation}
Moreover, with $\si_2 := \min(\de, \frac 32-\de,
\frac{2+5\de}{15})>0$ and $\tcp := t_c+n^{-3}$,
\begin{align}
\nonumber &\La_1(t) + \La_2(t) \lec \rho_c, \quad & (\forall t>t_c),
\\
\label{eq:4.66} &\La_1\lec \rho_0\bka{t}^{-\si} +
n^{1/3}\rho_c^{4/3},\quad
\La_2\lec \rho_0\bka{t}^{-7/6} + n\rho_c^2, & (t_c<t<\tcp),%
\\
\nonumber &\La_1(t) \lec n^{1/3}\rho_c^{4/3}, \quad \La_2(t) \lec
n^{1+\si_2}\rho_c^2, \quad   & (t
> \tcp).
\end{align}
\end{lemma}

\myproof From \eqref{wt-eta.eqn}, we have
\begin{equation} \label{t.ch-L}
e^{\bL (t-t_c)} \eta_\pm(t_c) = e^{\bL t} \eta_\pm(0) + \int_0^{t_c}
e^{\bL(t-s)} P_\pm \{F_{L\pm} + e^{\mp i\theta}J[F]\}(s) ds.
\end{equation}
We also decompose $\eta_\pm(t_c) = \eta^{(2)}_\pm (t_c)+
\eta^{(3)}_{\pm}(t_c)$ with a similar formula for $e^{\bL (t-t_c)}
\eta^{(3)}_\pm(t_c)$. We can bound $e^{\bL (t-t_c)} \eta_\pm(t_c)$
in $L^p$ and $e^{\bL (t-t_c)} \eta^{(3)}_\pm(t_c)$ in $L^2_{loc}$
using the same proof for Lemma \ref{disp.eta.est} with the integrand
set to zero for $t_c < s < t$. We also have
\begin{equation}
\norm{e^{\bL(t-t_c)}\eta^{(2)}_{\pm}(t_c)}_{L^2\loc} \lec
\wei{t-t_c}^{-3/2}n\rho_c^2
\end{equation}
using the explicit definition of $\eta^{(2)}_{\pm}$ in
\eqref{wt-eta23.def} and Lemma \ref{th:sdecay}. The above shows
\eqref{eq:4.64}.

We now show \eqref{eq:4.66}. Its first part is because $\rho_0
\bka{t}^{-1/2} \le \rho_c$  for all $t\ge t_c$, which follows from
\eqref{rho0rho}.

Its second part follows from $\rho(t) \sim \rho_c < \rho_0$.

For the third part with $t>\tcp$,  it suffices to show
\begin{equation}
\label{eq:4.72}\rho_0 \bka{t}^{-\si} \lec
n^{2\si-1}\rho_c^{2\si-2\al}, \quad \rho_0 \bka{t}^{-7/6} \lec
n^{1+\si_2}\rho_c^2.
\end{equation}

If $t_c< \De t$, then $\rho \sim \rho_c \sim \rho_0$. Writing all
factors as powers of $n$ using $\bka{t}^{-1} \le n^3$,
\eqref{eq:4.72} is reduced to $1+\de + 3 \si
> 2\si -1 + (2 \si - 2\al)(1+\de)$
and $1+\de + 7/2> 1+\si_2 + 2(1+\de)$. Both are valid using
$2/3<\si<3/4$, $0<\de<3/2$ and $\si_2<3/2-\de$.

If $t_c> \De t$, then $\rho_c \sim n^{-1}t_c^{-1/2}$, and
\eqref{eq:4.72} is reduced to $n^{1+\de} \bka{t}^{-\si} \lec
n^{-1+2\al} t_c^{-\si+\al}$ and $n^{1+\de} \bka{t}^{-7/6} \lec
n^{-1+\si_2}t_c^{-1}$, both are correct. \myendproof

%==================================================================

%==================================================================

\section{Escaping from an excited state}\label{S:5.2}
In this section we study the dynamics near an excited state when
$t>t_c$ assuming $t_{c} < \I$. We want to show that the solution
will escape from the $\rho_0$-neighborhood of the excited state.
Recall $\rho_0 = n^{1+\de}$ with $0<\de<3/2$. (We need $\de \ll 1$
in next section but not here.)

Fix $Q= Q_{1,n(t_c)}$ and decompose $\psi(t)$ for $t_c\le t<T_e$ as
in \eqref{psi.best.dec} and \eqref{zeta.dec} with respect to this
fixed $Q$. At $t=t_c$ we have Lemma \ref{OG.t-ch} and, by definition
of $t_c$ and $M_{t_c} \le 3/2$,
\begin{equation}
|z_0(t_c)| \ge \e_4 n^{-1} \rho_c^2, \quad z_H(t_c) \le
\frac 32 \rho_c, \quad |a(t_c)| \le
\frac 34 D \rho_c^2, \quad \rho_c: = \rho(t_c).
\end{equation}

Let
\begin{equation}
\label{ga.def} \gamma(t) := |q_0(t)| + n^{5}|q_0(t)|^{1/2} +
\rho_c,
\end{equation}
where $q_0(t)$ is the perturbation of $p_0(t)$ defined in Lemma
\ref{L-NF}. It will be shown to be an upper bound for bound
states.%
\footnote{The term $n^{5}|q_0|^{1/2}$ is included in $\ga$ so that
$z_H \lec \ga$. Explicitly: The bound of $\norm{\eta}_{L^p}$
includes $n^{11}|q_0|$, see \eqref{eq:5.30b}. By \eqref{eq:5.20},
the bound of $\norm{\eta^3}_{L^{9/8}\cap L^{3/2}}$ and hence
$\|\eta^{(3)}\|_{L^2_{loc}}$ contains $n^{18}z_L^{m}$ where $m \to
11/6$ as $p \to 6$. To bound $z_H$ by $\ga$, we need $
\norm{\eta}_{L^2_{loc}}\lec n \ga^2$ for \eqref{eq:5.48} and $\ga=
|q_0|+\rho_c$ is insufficient.}
We have defined $\ga(t)$ in
terms of $|q_0|$ instead of $|z_0|$ so that it is non-decreasing in
$t$ (for $t>\tcp:=t_c+n^{-3}$).

Define
\begin{equation}
\label{to.def} t_{o} : = \sup \left \{ t \geq t_{c} :  z_L(s) <
2n^{1+\de}, \ \forall \ s \in [t_{c}, t)\right \}.
\end{equation}
The time $t_{o}$ is the time that $z_L$ becomes powerful enough in
orthogonal coordinates. The subscript $_o$ means ``out'' (of the
neighborhood). It follows from Proposition \ref{tc-pro} below that
$t_o < T_e$ and hence the decompositions \eqref{psi.best.dec} and
\eqref{zeta.dec} are valid at least slightly beyond $t_o$.

Recall %In this section, we fix
\begin{equation}
\frac {27}5 < p < 6, \quad \si =\si(p)= \frac{3p-9}{2p}, \quad \frac
23 < \si < \frac 34.
\end{equation}

The main result of the section is the following proposition.
%====================================================================
\begin{proposition} \label{tc-pro}
There exist  constants $C_3, D_1>0$, uniform in $n$, (with $C_3$
greater than that in Lemma \ref{OG.t-ch}), such that for all
$t_{c} \leq t \leq t_{o}$, we have
\begin{equation} \label{t.ch-to-t.out}
\begin{split}
& |q_0(t)-q_0(s)| \le \frac 1{10} \e_4 n^{-1}\rho_c^2,\quad
(t_c \le s \le t \le \tcp:=t_c+n^{-3}) ,\\
& \frac{|q_0(t)|}{|q_0(s)|} \in [e^{\frac{1}{2}(\Re \LLLL_0)
(t-s)},
e^{\frac{3}{2}(\Re \LLLL_0) (t-s)}],\quad (\tcp \le s < t),\\
& z_H(t) \leq  \sqrt{\frac {6D}{\ga_0}} \gamma(t), \quad |a(t)| \leq D_1\gamma^2, \\
& \norm{\eta(t)}_{L^p} \leq n^{\si_1}\gamma(t)^2
+\frac 12 \La_1(t), \quad \si_1 =4\si-3-\al  ,\\
& \norm{\eta^{(3)}(t)}_{L^2\loc} \leq C_3 n^5\gamma(t)^2 + C_3
\gamma(t)^3 +\frac 12 \La_2(t) ,
\end{split}
\end{equation}
where $\al>0$ is so small that $-\frac 13+2\al <\si_1=
\frac{3(p-6)}p -\al <0$, and $\La_1(t)$ and $\La_2(t)$ are defined
in \eqref{eq:4.64}. In particular, $t_0 \le T_e$ and for some
constants $c_1$ and $c_2$,
\begin{equation} \label{to-tc}
\begin{split}
t_c + c_1 n^{-4} \log\frac{2\rho_0}{z_L(t_c)} \le t_{o} \le t_c +
c_2 n^{-4} \log\frac{2\rho_0}{z_L(t_c)} .
\end{split}
\end{equation}
\end{proposition}

The main term in the integrand of $\eta$ is of order $nz^2$. In the
first term of its ${L^p}$-bound we lose some powers of $n$ due to
integration over a time interval of order $n^{-4}$. On the other
hand, the first term $\ga^3$ of $\norm{\eta(t)}_{L^2\loc} $ estimate
is optimal and comes from recent time terms of order $z^3$ in the
integrand.

\medskip
%==============================================================================
\myproof The lemma clearly holds true for $t= t_{c}$. By a
continuity argument, it suffices to prove the lemma with additional
weaker assumptions:
\begin{equation} \label{t.ch-to-t.out-weak}
\begin{split}
& |q_0(t)-q_0(s)| \le \frac 1{2} \e_4 n^{-1}\rho_c^2,\quad
(t_c \le s \le t \le \tcp) ,\\
& \frac{|q_0(t)|}{|q_0(s)|} \in [e^{\frac{1}{4}(\Re \LLLL_0)
(t-s)}, e^{2(\Re \LLLL_0) (t-s)}],\quad (\tcp \le s < t),
%\quad |z_0(t)| = (1+o(1))|q_0(t)|,
\\
& z_H(t) \leq 2 \sqrt{\frac {6D}{\ga_0}}\gamma(t),\quad  |a(t)| \leq
2D_1\gamma^2,
\\
& \norm{\eta(t)}_{L^p} \leq 2n^{\si_1}\gamma(t)^2
+2\La_1(t), \\
& \norm{\eta^{(3)}(t)}_{L^2\loc} \leq 2C_3 n^5\gamma(t)^2 +
2C_3\ga(t)^3 +2\La_2(t) .
\end{split}
\end{equation}
At least for $t$ near $t_c$, the assumptions of Lemma \ref{L-NF}
are satisfied and hence $|z_0| \le |q_0|+|p_0-q_0| \le \ga + C
n\ga^2=(1+o(1))\ga$. Together with \eqref{t.ch-to-t.out-weak} and
$[\eta] = \eta^{(2)} + \eta^{(3)}$, the assumptions of Lemmas
\ref{Basic.lemma}--\ref{L-NF} are valid until $t=t_o$ with $\be =
(1+o(1))\ga(t)$, and
\begin{equation} \label{sim-L25}\begin{split}
|z_0(t)| &\le (1+o(1)) \ga(t) , \\
\norm{\eta(t)}_{L^2_{loc}} &\le C n\gamma^2(t)+  \La_2(t),\\
\norm{\eta(t)}_{L^2_{loc}\cap L^p} &\le \ga(t).
\end{split}\end{equation}
Here we have used \eqref{eq:4.66}.

It is convenient to have an upper bound of $\ga$ in terms of
$|q_0|$. Clearly
\begin{equation}
\ga^2(t) \sim |q_0|^2 +n^{10} |q_0|+\rho_c^2 \lec \e_4^{-1}n|q_0(t)|
+ \e_4^{-1} n |z_0(t_c)|.
\end{equation}
Since $|z_0(t_c)|\le |q_0(t_c)|+C n\ga(t_c)^2 \le |q_0(t)|+C
n\ga(t)^2$, we get
\begin{equation} \label{ga.est}
\ga^2(t) \lec \e_4^{-1}n|q_0(t)|.
\end{equation}
Thus we get an improved $z_0$ estimate,
\begin{equation} \label{z0.est}
|z_0|\le |q_0| + C n \ga^2 \le (1+o(1)) |q_0|.
\end{equation}
We can also derive from \eqref{t.ch-to-t.out-weak} and $|z_0(t_c)|
\ge \e_4 n^{-1} \rho_c^2$ that, for any $t_c \le s < t < t_{o}$,
\begin{equation} \label{eq:5.11}
|q_0(s)|\le \frac 65 |q_0(t)| e^{-\frac{1}{4}(\Re \LLLL_0) (t-s)}.
\end{equation}

We now give error estimates. For ${X_1}= n\norm{\eta}_{L^2_{loc}}^2
+ \norm{\eta^3}_{L^1_{loc}}$, using \eqref{t.ch-to-t.out-weak},
\eqref{sim-L25}, and H\"older inequality, we have
\begin{equation}% \label{XXYY}
\begin{split}
{X_1}  & \lec n(n^2 \ga^4 + \La_2^2)   +(n \ga^{2} + \La_2)^A (
n^{\si_1 } \ga^{2} +\La_1)^B ,
\end{split}
\end{equation}
with $A=\frac{2p-6}{p-2}$ and $B=\frac{p}{p-2}$. We claim that
\begin{equation} \label{de1.est}
X_1(t) \lec \left \{
\begin{aligned}
 &n\ga^2,  &&(\forall t>t_c),
\\
& n\rho_0^2 \bka{t}^{-7/6} + n^{2.8} \ga^4
%+n^{16.4}\rho_c^{2.8}z_L^{3/2}
, \quad &&( t_c < t <\tcp),
\\
&n^{2.8} \ga^4,  &&(t >\tcp) .
\end{aligned}\right .
\end{equation}
The first estimate is because $\La_1+\La_2 \lec \rho_c$. The last
estimate is, using \eqref{eq:4.66}$_3$ and $1.4 <A<1.5<B<1.6$ with
$A+B=3$,
\begin{equation}% \label{XXYY}
\begin{split}
X_1(t) & \lec  n^3 \ga^4 + (n\ga^2)^A (n^{1/3}\ga^{4/3})^B = n^3
\ga^4 + (n\ga)^{2A/3} n\ga^4 \lec n^{2.8} \ga^4 .
\end{split}
\end{equation}
When $t_c < t <\tcp$, using $\rho \sim \rho_c < \rho_0$,
\eqref{eq:4.66}$_2$, $\si_1>-1/3$, and the previous estimate,
\begin{equation}\begin{split}
X_1(t)
%&\lec  n^3 \ga  ^4 + n\La_2^2 +(n \ga ^{2} + \La_2)^A (
%n^{\si_1 } \ga ^{2} + \La_1)^B
%\\
&\lec  n^3 \ga  ^4 + n\rho_0^2 \bka{t}^{-7/3} +(\rho_0\bka{t}^{-7/6}
+ n\ga ^2)^A (\rho_0\bka{t}^{-\si} +n^{1/3}\ga ^{4/3})^B
\\
& \lec n\rho_0^2 \bka{t}^{-7/6} + n^{2.8} \ga ^4 .
%n^{16.4}\ga ^{2.8} z_L^{3/2} .
\end{split}\end{equation}

%and $a+b/4>1.772$.

For $\wt X$ and $X$  defined in \eqref{XwtX.def}, we have
\begin{equation} \label{XXYY}
\begin{split}%
\wt X &\leq \ga^2\norm{\eta}_{L^2_{loc}} + X_1 \leq
n\gamma^4+\ga^2 \La_2+X_1(t), \\
 X &\leq n \ga\norm{\eta}_{L^2_{loc}} +
X_1 \leq n^2\gamma^3+n\ga \La_2 + X_1(t).
\end{split}
\end{equation}

For $\hat X_p$ defined in \eqref{hat.X.def} we have
\begin{equation} %\label{XXYYhat}
\begin{split}
\hat{X}_p &= n^4z_L\norm{\eta}^2_{L^p} + n^6z_L^2\norm{\eta}_{L^p}
+ n^{6(6-p)/p}\norm{\eta}_{L^p}^3 %
\\
& \lec n^4z_L(n^{2\si_1}\ga^4 +  \La_1^2) + n^6 z_L^2(n^{\si_1}\ga^2
+ \La_1) + n^{6(6-p)/p}(n^{3\si_1}\ga^6 + \La_1^3).
\end{split}
\end{equation}
Using Young's inequality on $n^4z_L \La_1^2+n^6 z_L^2 \La_1$, and
$6(6-p)/p+3 \si_1 = \si_1 -2 \al> -1/2$,  we get
\begin{equation} \label{XXYYhat}
\begin{split}
\hat{X}_p &\lec n^{3/2} \ga^4 + n^{8.5} z_L^3 + n^{6(6-p)/p} \La_1^3
.
\end{split}
\end{equation}

From \eqref{eq:all}, \eqref{z.eqn}, Lemmas \ref{Basic.lemma},
\ref{au.est.lemma} and \eqref{t.ch-to-t.out-weak}, \eqref{XXYY} and
\eqref{eq:4.66}$_2$, we get
\begin{equation} \label{dot-theta.p}
\begin{split}
|\dot{\theta}| & \les \beta^2 + n^{-1}X \les
\ga^2+n^{-1}(n^2\ga^3+n\ga\La_2+X_1)\lec \ga^2
, \\
|\dot{p}_k| & \les n^4z_L + n\beta^2 +\hat X_p + X\les n^4z_L + n
\gamma^2+X_1 \les n^4z_L + n \gamma^2.
\end{split} \end{equation}

We now estimate the main terms. By H\"older inequality,
\begin{equation}
\label{eq:5.20}
\begin{split}
& \norm{\eta^3}_{L^{p'}} \leq \norm{\eta}_{L^2}^{\frac{2(p-4)}{p-2}}
\norm{\eta}_{L^p}^{\frac{p+2}{p-2}}, \quad \norm{\eta^3}_{L^{9/8}
\cap L^{3/2}} \leq \norm{\eta}_{L^2}^{\frac{2(2p-9)}{2(p-2)}}
\norm{\eta}_{L^p}^{\frac{11p}{9(p-2)}}.
\end{split}
\end{equation}
Using $36/7< p < 6$ and $-\frac12 < \si_1 =4\si-3 -\al =
3-\frac{18}p-\al<0$,
\begin{equation}
(n^{4\si -3-\al}\ga^2)^{\frac{p+2}{p-2}} \le (n^{4\si
-3-\al}\ga^2)^{\frac{11p}{9(p-2)}} \le o(1) \gamma^3,
\end{equation}
for $\al>0$ sufficiently small.  By Lemma \ref{Basic.lemma} and
$\norm{\eta}_{L^{2}} \le o(1)$, we get
\begin{equation} \label{NLterm.tc}
\begin{split}
\norm{F}_{L^{p'}} &\lec n\gamma^2+ X + n \norm{\eta}_{ L^p}^2 +
\norm{\eta^3}_{L^{p'}}\lec n\ga^2+\de_2 ,
\\
\norm{F-F_1}_{L^{9/8}\cap L^{3/2}} &\lec \gamma^3 + X + n
\norm{\eta}_{ L^p}^2 + \norm{\eta^3}_{L^{9/8}\cap L^{3/2}}\lec
\ga^3+\de_2 ,
\\
\de_2(t):&= n \ga(t) \La_2(t)+n \La_1^2(t).%+n^{18}z_L^{1.8} .
%\\
%& \lec  \ga(t) \La_2(t)+n \La_1^2(t) +n^a\ga^{2a}\La_1^b.
\end{split}
\end{equation}
In deriving the above estimates most terms in $X_1$ are controlled
by $\de_2$ except
\begin{equation}
n^A\ga^{2A}\La_1^B \le (n^{A-B/2} \ga^{2A}) (n^{B/2} \La_1^B) \lec
(n^{A-B/2} \ga^{2A})^{2/(2-B)} +(n^{B/2} \La_1^B)^{2/B} \lec \ga^3 +
n \La_1^2.
\end{equation}

Estimates \eqref{t.ch-to-t.out} now follows from Lemmas \ref{dis-tc}
and \ref{bou-tc} below.

In particular, taking $s=\tcp$ and $t=t_o$,
\eqref{t.ch-to-t.out}$_2$ together with $\Re \LLLL_0 \sim n^{-4}$
and $|z_0|=(1+o(1))|q_0|$ imply \eqref{to-tc}. \myendproof

% ==============================================================
\begin{lemma}[Dispersion estimates]\label{dis-tc}
For all $ t_c \leq t \leq t_o$, we have
\begin{equation}\label{eq:th5-3}
\begin{split}
& \norm{\eta(t)}_{L^p} \leq [n^{\si_1} \ga^2 +\La_1](t), \quad
\norm{\eta^{(3)}(t)}_{L^2\loc} \leq [C_3 n^5\gamma^2 + C_3\gamma^3
+\La_2](t) .
\end{split}
\end{equation}
\end{lemma}

Note that $\La_j(t)$ may compete with the main terms for $t$ near
$t_c$ but decay rapidly.

\medskip

\myproof  We first estimate $\norm{\eta(t)}_{L^p}$. It suffices to
estimate $\eta_\pm$ with
\begin{equation}
\eta_\pm(t) = e^{\bL (t-t_{c})}\eta_\pm(t_{c}) + \int_{t_{c}}^t
e^{\bL(t-s)} P_\pm \{ F_{L\pm} + e^{\mp i\theta}J[F]\} ds.
\end{equation}
By Lemma \ref{th:decay}, we have
\begin{equation} \label{after-tc-L5.est}
\norm{\eta_\pm(t)}_{L^p} \lec \norm{e^{\bL
(t-t_{c})}\eta_\pm(t_{c})}_{L^p} + \int_{t_{c}}^t \alpha_p(t-s)
\{\norm{F_{L\pm}}_{p'} + \norm{F}_{p'}\}(s) ds.
\end{equation}
By Lemma \ref{OG.t-ch}, we have $\norm{e^{\bL (t-t_{c})}
\eta_\pm(t_{c})}_{L^p} \le \frac 14 \La_1(t)$. By \eqref{FL} and
\eqref{dot-theta.p}, we get
\begin{equation}
\label{FLpm.est} \norm{F_{L\pm}}_{L^{p'}\cap L^{9/8} \cap L^{3/2}}
\les |\dot{\theta}| [\norm{\eta}_{L^p} +n^{-1}|a| + |z|]  \les
\ga^2\cdot \ga.
\end{equation}
From this, \eqref{NLterm.tc}, \eqref{after-tc-L5.est}, and $X_1 \ll
n \rho_c^2$, we get
\begin{equation} \label{L5.tc.1}
\begin{split}
& \int_{t_{c}}^t \al_p(t-s)[\norm{F_{L\pm}}_{L^{p'}} +
\norm{F}_{L^{p'}}](s) ds \les \int_{t_{c}}^t \al_p(t-s)(n
\gamma(s)^2 +\de_2(s))ds.
\end{split}
\end{equation}

Recall $\ga^2 \sim |q_0|^2  + n^{10}|q_0|+\rho_c^2$. By
\eqref{t.ch-to-t.out-weak}, $\Re \LLLL_0 \sim n^4$ and $\int^t
|t-s|^{-\si} e^{-a(t-s)} ds\lec a^{\si-1}$,
\begin{equation} %\label{tch-tout-i}
\begin{split}
\int_{t_{c}}^t \al_p(t-s)n|q_0|^2(s) ds &\leq \int_{t_{c}}^t
\al_p(t-s)
n|q_0|(t)^2 e^{-\frac 14 \Re \LLLL_0 (t-s)} ds\\
& \leq C n^{4\si-3} |q_0|^2 (t) .%\le C n^{4\si-3} \ga^2(t).
\end{split}
\end{equation}
The integral of $nn^{10}|q_0|$, part of $\de_2$, is bounded in the
same way by $Cn^{4(\si-1)+11} |q_0|(t)$.

For $\rho_c^2$, we have
\begin{equation} \label{L5.tc.2}
\int_{t_{c}}^t \al_p(t-s)n\rho_c^2 ds  \lec n \rho_c^2
\bka{t-t_c}^{1-\si} = n^{4\si-3-\al/2}\cdot \rho_c^2 n^{\al/2}
T^{1-\si}
\end{equation}
where $\al>0$ is to be chosen and $T= n^{4} \bka{t-t_c}$. Let $A=
\frac 18 n^{-4}\Re\LLLL_0$ which is of order 1. If $AT \le 10\log
\frac 1n$, then $ n^{\al/2} T^{1-\si} = o(1)$ if $n$ is
sufficiently small. If $AT \ge 10\log \frac 1n$, then by
\eqref{eq:5.11}
\begin{equation}
\label{eq:5.30b} \rho_c^2  T^{1-\si} \le C  n  |q_0(t_c)|T^{1-\si}
\le C n|q_0(t)| e^{-2AT} T^{1-\si}.
\end{equation}
Since $e^{-AT} \le n^{10}$ and  $e^{-AT} T^{1-\si}\le C$, it is
bounded by $Cn^{11}|q_0(t)|$.

Using \eqref{eq:4.66}, the error term $\de_2(t)=n \ga(t) \La_2(t)+n
\La_1^2(t)$ is bounded by $n^{7/3} \rho_c^2$ when $t>\tcp$ and by $
n^{7/3} \rho_c^2+ n \rho_0^2\bka{t}^{-7/6}$ when $t<\tcp$.   The
term $n^{7/3} \rho_c^2$ is smaller than the main term $n \ga^2$ in
\eqref{L5.tc.1} and can be absorbed, while
\begin{equation}\label{eq:5.30}
\int_{t_c}^{\tcp} n\rho_0^2 \bka{t}^{-7/6} dt \lec n \rho_c^2
\end{equation}
which can be checked using $\rho_c \sim \rho_0$ for $t_c<\De t$ and
$\rho_c \sim n^{-1}t_c^{-1/2}$ for $t_c>\De t$.

Thus the integral in \eqref{after-tc-L5.est} is bounded by $n^{\si_1
}\ga^2$ with $\si_1 =4\si-3-\al$, and we have shown the first
estimate of \eqref{eq:th5-3} for $\norm{\eta}_{L^p}$.

\medskip

Next, we  estimate $\|\eta^{(3)} \|_{L^2_{loc}}$. Decompose
$\eta_{\pm}^{(3)} = \sum_{j=1}^4\eta_{\pm,j}^{(3)}$, where
$\eta_{\pm,j}^{(3)}$ are defined explicitly in \eqref{wt-eta14}
and \eqref{wt-eta23.def} with $t_0 = t_c$. From Lemmas
\ref{th:sdecay} and \ref{OG.t-ch}, we get
\begin{equation} \label{eta12.tc.est}
\begin{split}
\norm{\eta_{\pm,1}^{(3)}}_{L^2\loc} &\leq \frac 14 \La_2(t), \quad
\norm{\eta_{\pm,2}^{(3)}}_{L^2\loc} \leq \frac 14C_3
n\rho_c^2\wei{t-t_c}^{-3/2} \le \frac 14 \La_2(t).
\end{split}
\end{equation}
By \eqref{fkl-th.est} and \eqref{dot-theta.p}, we have
\begin{equation}
\begin{split}
\norm{|\dot f_{kl}| +|\dot{\theta}f_{kl}|}_{L^2_r} &\lec n |\dot
\th|\ga^2 + n \ga |\dot p| \lec n(\ga^2)\ga^2 + n\ga(n^4
\ga+n\ga^2)
\\
&\lec n^5\ga^2 + n^2 \ga^3 .
\end{split}\end{equation}
By Lemma \ref{th:sdecay} again and $\ga(s) \lec \ga(t)$ for $s<t$,
we obtain
\begin{equation} \label{eta3.tc.est}
\norm{\eta_{\pm,3}^{(3)}}_{L^2\loc} \lec\int_{t_c}^t
\wei{t-s}^{-3/2} [n^5\ga^2 + n^2 \ga^3 ](s)ds \lec [n^5\ga^2 + n^2
\ga^3 ](t).
\end{equation}

Finally, $\norm{\eta_{\pm,4}^{(3)}}_{L^2\loc}$ is  bounded by
$\int_{t_c}^t \alpha_\infty(t-s)I_4(s)  ds$ by Lemma
\ref{th:decay}, with
\begin{equation}%\label{FN2.term}
I_4 =\norm{F_{L\pm}}_{L^{9/8}\cap L^{3/2}} + \norm{F
-F_1}_{L^{9/8}\cap L^{3/2}}\lec \ga^{3} + \de_2
\end{equation}
by \eqref{FLpm.est} and \eqref{NLterm.tc}$_2$.  Using $\de_2(t)= n
\ga(t) \La_2(t)+n \La_1^2(t)$ and the explicit form of $\La_j$ in
\eqref{eq:4.65} together with the integral bound \eqref{int-est2},
we get
\begin{equation} \label{eta4.tc.est}
\begin{split}
\norm{\eta_{\pm,4}^{(3)}}_{L^2\loc} &\lec\int_{t_c}^t
\alpha_\infty(t-s) [\ga^{3} + \de_2](s) ds  \\ &\lec \ga^{3}(t) +
n\rho_0^2\bka{t}^{-7/6}+ n^{5/3}\rho^{7/3} \lec \ga^3(t) +
o(1)\La_2(t).
\end{split}
\end{equation}
Summing the above estimates, we get the second estimate of
\eqref{eq:th5-3} for $\norm{\eta^{(3)}}_{L^2\loc}$. \myendproof

%%%%%%%%%%%%%%%%%%%%%%%%%%%%%%%%%%%%%%%%%%%%%%%%%%%%%%%%%%%%%%%%%%

\begin{lemma}[Bound states estimates] \label{bou-tc}
There is a uniform in $n$ constant $D_1>0$ such that for
all $t_{c} \leq t \leq t_{o}$, we have
\begin{equation} %\label{t.ch-to-t.out}
\begin{split}
& |q_0(t)-q_0(t_c)| \le \frac 1{10} \e_4 n^{-1}\rho_c^2,\quad
(t_c \le t \le \tcp) ,\\
& \frac{|q_0(t)|}{|q_0(s)|} \in [e^{\frac{1}{2}(\Re \LLLL_0)
(t-s)},
e^{\frac{3}{2}(\Re \LLLL_0) (t-s)}],\quad (\tcp \le s < t),\\
& z_H(t) \leq \sqrt{ \frac {6D}{\ga_0}}\gamma(t), \quad |a(t)| \leq
D_1\gamma(t)^2.
\end{split}
\end{equation}
\end{lemma}

\myproof
First we estimate $q_0(t)$. From Lemma \ref{L-NF}, we have
\begin{equation} \label{qL.normal}
\dot q_0(t) = (\Re \LLLL_0) q_0 + \wt Y_0q_0 + g_0, \quad |q_0 -
p_0| \les n\gamma^2, \quad |\Re(\wt Y_0)| \leq Cn^2\gamma^2 \ll
n^4.
\end{equation}
Here $\wt Y_0  = Y_0 + \sum_{l \not = 1} D_{0l}|q_l|^2$. Moreover,
from \eqref{gk.est}, \eqref{XXYY} and \eqref{ga.est}, we have
\begin{equation}\begin{split}
|g_0| &  \leq C[n^5\gamma^2 % n^4 \ga^3
+ n\gamma^4 +
n^3\gamma\norm{\eta}_{L^2\loc} + n\gamma\norm{\eta^{(3)}}_{L^2\loc}
+ \hat{X}_p + \wt{X}] \le
%n^5\gamma^2 + o(1) n^2\gamma^3
o(1)n^4|q_0| + \de_3,
\end{split}\end{equation}
where $\de_3 = C(n^{6(6-p)/p}\La_1^3 + \ga^2 \La_2 + X_1)$. If
$t<\tcp$, by \eqref{eq:4.66}$_2$, \eqref{de1.est}$_2$ and
\eqref{eq:5.30},
\begin{equation}\begin{split}
\de_3(t) &\lec n\rho_0^2\bka{t}^{-7/6}+ n\rho_c^{4} + n \ga^2
\rho_c^2 + n^{2.8}\ga^4,
\\
|q_0(t)-q_0(t_c)|&\le \int_{t_c}^{\tcp} Cn^4|q_0| +  \de_3(s) ds \le
o(1)(|q_0(t_c)|+\e_4 n^{-1}\rho_c^2),
\end{split}\end{equation}
This shows the $q_0(t)$-estimate for $t<\tcp$. Suppose now $\tcp<t$.
By \eqref{eq:4.66}$_3$, \eqref{de1.est}$_3$, and \eqref{ga.est},
\begin{equation}
\de_3(t)\lec n^{6(6-p)/p}(n^{1/3} \rho_c^{4/3})^3 + \ga^2
n^{1+\si_2} \rho_c^2+ n^{2.8}\ga^4 \ll n^4 |q_0|.
\end{equation}
%using {\color{red}$\si>2/3$ (i.e. $p>27/5$)}.
Since $\Re \LLLL_0>0$ is of order $n^4$, Eq.~\eqref{qL.normal}
gives
\begin{equation}
0<\frac{1}{2} (\Re \LLLL_0) |q_0| \leq \frac{d}{dt}|q_0| \leq
\frac{3}{2}(\Re \LLLL_0) |q_0|,
\end{equation}
which implies the estimate of $|q_0(t)|$ for $t>\tcp$.

Next, we estimate $z_H(t)$. For any $k>1$, by Lemma \ref{L-NF}, we have
\begin{equation} \label{qH.n}
\frac{d}{dt}q_k = \sum_{l > 1} D_{kl} |q_l|^2 q_k + Y_kq_k + g_k,
\quad |q_k -p_k| \leq Cn\gamma^2.
\end{equation}
Moreover, we have
\begin{equation}
|D_{kl}| \leq Dn^2, \quad |\Re(Y_k)| \leq Dn^2|z_0|^2, \quad
\Re(D_{kl}) \leq -\frac {\gamma_0} 2 n^2, \quad \forall l >1.
\end{equation}
So, we have
\begin{equation}
\frac{d}{dt}(|q_k|) \leq -\frac{\ga_0n^2}2 \sum_{l {>1}} |q_l|^2
|q_k| + 2Dn^2 |q_0|^2|q_k| + |g_k|.
\end{equation}
Let $f(t) = (\sum_{l>1}|q_l|^2)^{1/2}$. We have $f(t_c)\lec \rho_c$
and
\begin{equation}
\dot f(t) \leq -\frac{\gamma_0n^2}{2}f^3 + 2D n^2|q_0|^2f(t) +
\sum_{k >1}|g_k|.
\end{equation}
On the other hand, from \eqref{gk.est}, we have
\begin{equation}
\label{eq:5.48} |g_k|  \leq C[n\gamma^4 + n^4 \ga^3+
n^3\gamma\norm{\eta}_{L^2\loc} +
n\gamma\norm{\eta^{(3)}}_{L^2\loc} + n\gamma\hat{X}_p + \wt{X}]
\le o(1) n^2\gamma^3 + \de_4,
\end{equation}
where $\de_4 = C (n\ga n^{6(6-p)/p}\La_1^3 + n\ga \La_2 + X_1)$.  If
$t \le \tcp$, by \eqref{eq:4.66}$_2$, \eqref{de1.est}$_2$ and
\eqref{eq:5.30}
\begin{equation}
\begin{split}
&\de_4(t) \lec n\rho_0^2\bka{t}^{-7/6}+ n^{2}\rho_c^{2} \ga + n^{2.8}\ga^4
\\
&|f(t) -f(t_c)| \le   \int_{t_c}^{\tcp} Cn^2 \rho_c^3 +\de_4(s) ds
\le C n \rho_c^2 \ll \rho_c.
\end{split}
\end{equation}
Thus $f(t) \lec \rho_c$ for $t < \tcp$. When $\tcp<t$, since
$\de_4(t)\le n^2\ga^5 + n^{2+\si_2}\ga^3 + n^{2.8}\ga^4\ll n^2
\ga^3$, for $\wt \ga = (\frac {16D}{3\ga_0})^{1/2}\ga$,
\begin{equation}
\dot f(t)  \leq \frac{\gamma_0n^2}{4}[\wt \gamma^3 -f^3], \quad
(t>\tcp).
\end{equation}
Since $\ga(t)$ is nondecreasing and $f(\tcp)< \wt \gamma(\tcp)$, by
comparison we get
\begin{equation}
f(t) \leq \wt \gamma(t), \quad \forall\ t> \tcp .
\end{equation}
Thus $z_H(t) \le f(t)+|f(t)-z_H(t)| \le \wt \gamma(t)+ C n
\ga^2(t)<\sqrt{\frac {6D}{\ga_0}}\ga(t)$.

\medskip

%============================================
Finally, we estimate $a(t)$. By \eqref{a.dec} and Lemma \ref{b.NF}, $a
= a^{(2)} + (b-\wt b) + \wt b$, where
\begin{equation}
|a^{(2)}| \les n^2\ga^2, \quad |\wt{b} - b| \leq Cn\gamma[\gamma^2
+n\norm{\eta}_{L^2\loc}]\leq Cn^2\gamma^2,
\end{equation}
and
\begin{equation}
\frac{d}{dt} \wt{b} = b_0 +\sum_{k,l \not = 1} B_{kl} |z_k|^2|z_l|^2
+ g_b.
\end{equation}
Using $a(t_c)=0$,
\begin{equation}
\begin{split}
|a(t)-0| & \le |a^{(2)}(t)|+ |a^{(2)}(t_c)|  + |(b-\wt b)(t)|
+|(b-\wt b)(t_c)|+ |\wt b(t)-\wt b(t_c)| \\
& \le C n^2 \ga^2(t) + \int_{t_c}^t|\frac{d}{dt} \wt{b} | .
\end{split}
\end{equation}
From \eqref{b0.def}, $b_0(t) = b_{00}|z_0(t)|^2$ with $b_{00}=
2\Im \kappa  c_0(Q^2,\bar u_0^+ u_0^-)$ and $|b_{00}|n^{-4}\le
C_4$ for some explicit $C_4=O(1)$. We also have $|B_{kl}|
|z_k|^2|z_l|^2 \lec n^2 \ga^4$ and
\begin{equation} \label{b0.gb.est}
\begin{split}
|g_b| & \le C[n^3\gamma^4 + n^2\beta^2|z_0|^2 + n\beta^5 + n^2|z_0|
\norm{\eta}_{L^2\loc}
+ nX_1 %n^2\norm{\eta}_{L^2\loc}^2 + n\norm{\eta^3}_{L^1\loc}
+ n\gamma^2\norm{\eta^{(3)}}_{L^2\loc} + n\gamma \hat{X}_p]
\\
& \le o(1) n^4|z_0|^2 +  C n^2\gamma^4 + \de_5,
\end{split}
\end{equation}
where $\de_5 = nX_1 +(n^2z_L+n\ga^2) \La_2 +  n \ga n^{6(6/p-1)}
\La_1^3$. Thus
\begin{equation}\label{eq:5.52}
|a(t)| \le C n^2 \ga^2(t) + \int_{t_c}^t(C_4+o(1)) n^4|q_0(s)|^2 +
C n^2 \ga^4 (s) +\de_5(s) ds.
\end{equation}
By \eqref{t.ch-to-t.out-weak},
\begin{equation}
\label{eq:4.115} \int_{t_c}^t(C_4+o(1)) n^4|q_0(s)|^2 ds  \le
\frac 65 C_4 n^4 |q_0(t)|^2\int_{t_c}^t e^{-\frac 14 \Re
\LLLL_0(t-s)} ds
 \le \frac{24 C_4 n^4 }{5\Re \LLLL_0}|q_0(t)|^2.
\end{equation}
Moreover, by the definition of $\ga$,
\begin{equation}
\int_{t_c}^tC n^2 \ga^4 (s)ds  \lec \int_{t_c}^t [n^2 |q_0|^4+
n^{22}|q_0|^2] (s)ds + n^2 \rho_c^4(t-t_c).
\end{equation}
The integral is bounded by $n^{-2}|q_0|^4+ n^{18}|q_0|^2= o(1)
|q_0|^2$ similarly as in \eqref{eq:4.115}, while the last term is
bounded by $n^2 \rho_c^4 C n^{-4} \log \frac {|z_0|(t)}{\e_4 n
\rho_c^2} = o(1) \rho_c^2$. Thus this term is $o(1) \ga^2$.

For the error term $\int_{t_c}^t \de_5(s) ds$, if $t \le \tcp$, by
\eqref{eq:4.66}$_2$ and \eqref{de1.est}$_2$ we have
\begin{equation}
\begin{split} \de_5(s) &\le
n^2\rho_0^2\bka{t}^{-7/6}+n^{3.8}\ga^4 + (n^2|q_0(t_c)|+n\ga^2)
(\rho_0
\bka{t}^{-7/6}+n\rho_c^2)\\
&\qquad + n \ga(t_c) (\rho_0^3 \bka{t}^{-3\si} + n\rho_c^{4})
\\
& \le n^2\rho_0^2\bka{t}^{-7/6}+o(1)n^{4}\ga^2 .
\end{split}
\end{equation}
Thus, using \eqref{eq:5.30}, we have $\int_{t_c}^t\de_5(s) ds \le
o(1) n\ga(t_c)^2$. If $t>\tcp$, by \eqref{eq:4.66}$_3$ and
\eqref{de1.est}$_3$ we have $ \de_5(s) \le n^{3.8} \ga^4 + n^2 \ga
n^{1+\si_2}\rho_c^2 + n \ga n \rho_c^4 = o(1) (n^2 \ga^4 + n^4
\ga^2)$, which is dominated by other terms in \eqref{eq:5.52}.

In conclusion, we have shown
\begin{equation}
\label{D1.def} |a(t)| \le D_1\gamma^2(t), \quad D_1 :=\frac{5 C_4
n^4 }{\Re \LLLL_0}=O(1).
\end{equation}

This completes the proof of the Lemma \ref{bou-tc}. \myendproof

\medskip

The above finishes the proof of Proposition \ref{tc-pro}.

\bigskip

%==================================================================
We now prove the following out-going estimate of $\eta$ at $t_o$.

\begin{lemma} \label{eta*.OG}
For some $C_5>0$, for all $t \geq t_{o}$, we have
\begin{equation}\begin{split}
\norm{e^{(t-t_o) \bL } \eta_\pm(t_o)}_{L^p} & \leq  \wt{\La}_1(t):
=\La_1(t) + C_5 n^{-2}\rho_0 (n^{-4}+t-t_o)^{-\si},
\\
\norm{e^{(t-t_o) \bL} \eta_\pm(t_o)}_{L^2_{loc}} & \leq
\wt{\La}_2(t):= \La_2(t) + C_5 n \rho_0^2 \bka{t-t_o}^{-7/6} \\
&\quad + C_5\rho_0^3 \bka{t-t_o}^{-1/6} n^{-4} (t-t_o+n^{-4})^{-1}
\\ &\quad
 + C_5 n^{-3} (n^{7/3}\rho_c +
\rho_c^2) (t-t_o+n^{-4})^{-7/6} .
\end{split}\end{equation}
\end{lemma}

\myproof %Denote $\tau = t-t_o$.
For all $t \geq t_{o}$, we have
\begin{equation}
e^{\bL (t-t_o)} \eta_\pm(t_o) = e^{\bL (t-t_{c})} \eta_\pm(t_{c}) +
\int_{t_{c}}^{t_{o}} e^{\bL(t -s)} P_\pm\{F_{L\pm} +
Je^{i\theta}[F]\} ds.
\end{equation}
We first bound it in $L^p$. By Lemma \ref{OG.t-ch}, the first term
is bounded in $L^p$ by $\La_1(t)$. The second term is bounded  in
$L^p$ as in \eqref{L5.tc.1} by
\begin{equation}\begin{split}
&\lec \int_{t_{c}}^{t_{o}}\al_p(t-s) [\norm{F_{L\pm}}_{L^{p'}} +
\norm{F}_{L^{p'}}] ds \lec \int_{t_{c}}^{t_{o}}\al_p(t-s) [
n\ga^2(s) + \de_2(s)] ds.
\end{split}\end{equation}
Note $n\ga^2+\de_2 \sim  n|q_0|^2  + n^{11}|q_0|+ n\rho_c^2+\de_2$.
By \eqref{t.ch-to-t.out-weak},
\begin{equation} \label{tch-tout-i}
\int_{t_{c}}^{t_o} \al_p(t-s)n|q_0|^2 (s) ds \leq \frac
65\int_{t_{c}}^{t_o}  \al_p(t-s) n\rho_0^2 e^{-\frac 14 \Re
\LLLL_0 (t_o-s)}ds.
\end{equation}
Using
\begin{equation}
\int^{t_o} |t-s|^{-\si} e^{-(t_o-s)/T} ds \lec \int_{t_o-T}^{t_o}
|t-s|^{-\si} e^{-(t_o-s)/T} ds \lec \int_{t_o-T}^{t_o} |t-s|^{-\si}
 ds \lec T(t-t_o+T)^{-\si}
\end{equation}
with $T=4/ \Re \LLLL_0 \sim n^{-4}$, \eqref{tch-tout-i} is bounded
by $C n^{-3} \rho_0^2 (t-t_o+n^{-4})^{-\si}$.

Similarly $\int_{t_{c}}^{t_o} \al_p(t-s)n^{11}|q_0| (s) ds$ is
bounded by $n^{11}\rho_0n^{-4}(t-t_o+n^{-4})^{-\si}$.

Let $t_k$ denote the first time in $[t_c, t_o)$ so that
$|q_0(t)|=\rho_c$.  When $t>t_k$, the integrand $\rho_c^2$ is
dominated by $|q_0|^2$ and can be absorbed. By
\eqref{t.ch-to-t.out-weak}, $t_o - t_k \gec n^{-4}\log \frac
{2\rho_0}{\rho_c}$. We have
\begin{equation} \label{eq:5.62}
\int_{t_{c}}^{t_k} \al_p(t-s)n\rho_c^2 ds  \lec n \rho_c^2 |t_k-t_c|
|t-t_k|^{-\si}.
\end{equation}
Using
\begin{equation} \label{eq:5.63}
 \rho_c^2 \lec \e_4^{-1} n q_0(t_c) \lec \frac 65 \e_4^{-1} n
\rho_c e^{-\frac 14 \Re \LLLL_0(t_k-t_c)},
\end{equation}
and $n^4|t_k-t_c|e^{-\frac 14 \Re \LLLL_0(t_k-t_c)}\le C$, the
integral in \eqref{eq:5.62} is bounded by $C
n^{-2}\rho_c|t-t_k|^{-\si}$.

Using \eqref{eq:4.66}, the error term $\de_2(t)$ is bounded by
$n^{7/3} \rho_c^2$ when $t>\tcp$ and by $ n^{7/3} \rho_c^2+ n
\rho_0^2\bka{t}^{-7/6}$ when $t<\tcp$.   The term $n^{7/3} \rho_c^2$
is much smaller than the main terms and can be absorbed, while by
\eqref{eq:5.30},
\begin{equation}
\int_{t_c}^{\tcp}\al_p(t-s) n\rho_0^2 \bka{s}^{-7/6} ds \lec n
\rho_c^2 |t-t_c|^{-\si}.
\end{equation}

Summing the above estimates gives the first estimate of Lemma
\ref{eta*.OG}.
\medskip

For the second estimate, we have $\eta_\pm(t_o) =
\eta_\pm^{(2)}(t_o) + \eta_\pm^{(3)}(t_o)$. By \eqref{wt-eta14},
\eqref{wt-eta23.def} and \eqref{wt-eta3.def} with $t_0$ replaced by
$t_c$, we have for $\tau = t- t_o\ge 0$
\begin{equation}
e^{\bL \tau}\eta_\pm(t_o) = e^{\bL \tau} \eta^{(2)}_\pm(t_o) +
\sum_{j=1}^4 e^{\bL \tau} \eta_{\pm,j}^{(3)}(t_o),
\end{equation}
with
\begin{equation}
e^{\bL \tau} \eta^{(3)}_{\pm,1}(t_o) = e^{(t-t_c)\bL}
\eta_\pm(t_c), \quad e^{\bL \tau} \eta^{(3)}_{\pm,2}(t_o) =
-e^{(t-t_c)\bL}\eta^{(2)}_\pm(t_c),
\end{equation}
\begin{equation*}
e^{\bL \tau} \eta^{(3)}_{\pm,3}(t_o)  = - \int_{t_c}^{t_o}
 e^{(t-s) \bL} e^{\mp
i\theta(s)} \Pi_\pm  \tsum_{k,l \in \Omega_m} \big(\Re R_{
kl}e^{-i\om_{kl}s} \dot f_{kl}  \mp i\Re R_{ kl}e^{-i\om_{kl}s}
\dot{\theta}f_{kl} \big) (s)ds,
\end{equation*}
\begin{equation}
e^{\bL \tau} \eta^{(3)}_{\pm,4}(t_o)  =
\int_{t_c}^{t_o}e^{(t-s)\bL} P_\pm \{F_{L\pm} + J e^{\mp
i\theta}[F-F_1]\} ds.
\end{equation}

From the explicit definition of $\eta^{(2)}_{\pm}(t_o)$ in
\eqref{wt-eta23.def} and Lemma \ref{th:sdecay} we obtain
\begin{equation} \label{eta(2).to.est}
\norm{e^{\bL \tau} \eta^{(2)}_\pm(t_o)} \leq Cn \rho_0^2
\wei{t-t_o}^{-3/2}.
%= Cn \rho_0^2\wei{t-t_o + \Delta t}^{-3/2}.
\end{equation}
By  Lemma \ref{OG.t-ch},
\begin{equation} \label{eta12.to.est}
%\begin{split}
\norm{e^{\bL \tau} \eta^{(3)}_{\pm,1}(t_o)}_{L^2\loc}  \leq \frac
14 \La_2(t),  \quad %\\
\norm{e^{\bL \tau} \eta^{(3)}_{\pm,2}(t_o)}_{L^2\loc} \leq
C_3n\rho_c^2\wei{t-t_c}^{-3/2}.
%\end{split}
\end{equation}
As in \eqref{eta3.tc.est} and \eqref{eta4.tc.est}, we obtain
\begin{equation} %\label{eta34.to.est}
\begin{split}
\norm{e^{\bL \tau} (\eta^{(3)}_{\pm,3}+\eta^{(3)}_{\pm,4})
(t_o)}_{L^2\loc} & \lec \int_{t_c}^{t_o} \al_\I(t-s)[n^5\gamma^2 +
\gamma^3+\de_2](s)ds \le I_1 + I_2 + I_3,
\end{split}
\end{equation}
where  $I_j$ are integrals over the same time interval with the
following integrands  \begin{equation} (n^5|q_0|^2 + n^{15}|q_0| +
|q_0|^3 + n^{15}|q_0|^{3/2} ), \quad  (n^{7/3} \rho_c^2 +
\rho_c^3)1_{[t_c,t_k]} , \quad n\rho_0^2 \bka{s}^{-7/6}
1_{[t_c,\tcp]}.
\end{equation}
Then
\begin{equation} %\label{eta34.to.est}
\begin{split}
I_1(t) &\lec \int_{t_c}^{t_o} \al_\I(t-s) \rho_0^3 e^{-\frac 14\Re
\LLLL_0(t_o-s)}ds \lec \rho_0^3
\int_{t_o-n^{-4}}^{t_o} \bka{t-s}^{-7/6} ds \\
&\le \rho_0^3 \bka{t-t_o}^{-1/6}n^{-4}(t-t_o+n^{-4})^{-1}.
\end{split}
\end{equation}
With constant $\e =n^{7/3} \rho_c^2 + \rho_c^3$, using
\eqref{eq:5.63} and $n^4(t_k-t_c) e^{-\Re \frac 14
 \LLLL_0(t_k-t_c)} \le C$,
\begin{equation} %\label{eta34.to.est}
\begin{split}
I_2(t) &\lec \int_{t_c}^{t_k} \al_\I(t-s) \e ds \le \e
(t-t_k)^{-1/6} (t_k-t_c) (t-t_c)^{-1} \\
& \le  \e (t-t_k)^{-1/6}
 (t-t_c)^{-1} n^{-4} n^4(t_k-t_c) n \rho_c^{-1} e^{-\Re \frac 14
 \LLLL_0(t_k-t_c)}
\\
&\le \e_4^{-1}n^{-3} (n^{7/3}\rho_c+ \rho_c^2) (t-t_k)^{-1/6}
 (t-t_c)^{-1}.
\end{split}
\end{equation}
Finally, $I_3(t) \lec \int_{t_c}^\tcp \al_\I(t-s) n\rho_0^2
\bka{s}^{-7/6}ds \le (t-t_c)^{-7/6}n\rho_0^2$. Summing the estimates
we get the second part of the Lemma.
 \myendproof

%==========================================================================

%=============================================================
\section{Dynamics away from bound states}
In this section, we study the dynamics of the solution $\psi(t)$
for $ t_o\le t \le t_i$, where $t_o$ is the time it leaves
$2\rho_0$ neighborhood of first excited states, and $t_i$ is the
time it enters the $\rho_0$-neighborhood of ground states, to be
defined in \eqref{ti.def}. In this time interval we use orthogonal
coordinates and decompose
\begin{equation}
\psi(t) = \sum_{j=0}^K x_j(t) \phi_j + \xi(t), \quad \xi(t)\in
\bE_c^{H_0}, \quad (t \ge t_o).
\end{equation}

We first estimate $x_j(t_o)$ and $\xi(t_o)$ in Lemma \ref{to.ch},
for which we recall some definitions. Recall that $\Delta t =
n^{-2}\rho_0^{-2} = n^{-2(2 +\de)}$, $0 < \al \ll 1$ is fixed and $0
< \de \leq \frac{1}{10}$. Moreover, $\frac{27}{5} < p <6$ is fixed,
$\frac{2}{3} < \sigma = \frac{3(p-3)}{2p} < \frac{3}{4}$, and
$\sigma' := \frac{3(p-2)}{2p}> \sigma$.  Recall from Lemma
\ref{eta*.OG} that $\wt{\La}_2 = \wt{\La}_{2,1} + \wt{\La}_{2,2}$
with
\begin{equation}\label{La2.dec}
\begin{split}
\wt{\La}_{2,1}(t) & := \La_2(t) + C_5 n \rho_0^2 \bka{t-t_o}^{-7/6}  + C_5 n^{-3} (n^{7/3}\rho_c + \rho_c^2) (t-t_o+n^{-4})^{-7/6},\\
\wt{\La}_{2,1}(t) & := C_5\rho_0^3 \bka{t-t_o}^{-1/6} n^{-4}
(t-t_o+n^{-4})^{-1}.
\end{split}
\end{equation}
We also define
\begin{equation} \label{La34.def}
\begin{split}
\Lambda_3(t) & := 3\wt{\La}_2(t) + C_6n^3(1+t-t_o)^{-3/2},  \quad \Lambda_4(t) := \sum_{j=1}^3 \La_{4,j}(t), \\
%\Lambda_4(t) & := 3C_6[n^{-1}(\Delta t)^{-\al}[\Delta t +t]^{-\sigma +\al} + \rho_0(1+t-t_o)^{-\sigma}  + n^{-2}\rho_0 (n^{-4}+t-t_o)^{-\si}].
\end{split}
\end{equation}
where $C_6$ is some uniform constant defined in \eqref{xi*1.est}
and
\begin{equation} \label{La4.dec}
\begin{split}
& \La_{4,1}: = C_6n^{-1 + (4+2\de)\al}(\Delta t +t)^{-\sigma
+\al}, \quad
\La_{4,2} := C_6\rho_0(1+t-t_o)^{-\sigma}, \\
& \La_{4,3} := C_6n^{-1+\de} (n^{-4}+t-t_o)^{-\si}.
%, \quad \sigma'= \frac{3}{2} - \frac{3}{p}.
\end{split}
\end{equation}
Note that $\La_{4,1}$ is the second term in $\La_1$ and comes from
the out-going estimate at $t_c$; $\La_{4,3}$ is from the out-going
estimate at $t_o$ and $\La_{4,2}$ is from \eqref{xi*1.est}.  Also
note that
\begin{equation} \label{La34p.est}
\La_3(t) \leq 3C_6n^3, \quad \La_4(t) \leq 2C_6n^{\frac{5p-18}{p}
+\de} + C_6\rho_0 \wei{t-t_o}^{-\sigma}, \quad \frac{5p-18}{p} >
\frac{5}{3}.
\end{equation}
%================================================================================
\begin{lemma} \label{to.ch}At $t = t_{o}$ we have
\begin{equation} \label{bes.tout}
(1.9)n^{1+\delta} \leq |x_0| \leq (2.1) n^{1+\delta}, \quad (\sum_{k
> 1} |x_k|^2)^{\frac12} \leq 6 \sqrt {\frac D {\ga_0}}\rho_0, \quad (0.9) n \leq
|x_1| \leq (1.1)n.
\end{equation}
Moreover, we have for all $t \ge t_o$
\begin{equation} \label{disp.tout}
%\begin{split}
\norm{e^{-i(t-t_o) H_0}\xi(t_o)}_{L^2\loc}  \leq  \Lambda_3(t), \quad \norm{e^{-i(t-t_o)H_0} \xi(t_o)}_{L^p} \leq \Lambda_4(t).
%\end{split}
\end{equation}
\end{lemma}
%=========================
\myproof For all $0 \leq t \leq t_{o}$, we have
\begin{equation} \label{t-out-e}
 \psi = [Q + a(t) R + \zeta + \eta]^{-iEt + i\theta} = \sum_{j =0}^K x_j \phi_j + \xi.
\end{equation}
Here $Q = Q_{1, n(t_c)}$. Recall $n(t_c) = n+O(n^{1+2\de})$ by
substituting \eqref{t-out-e} with $t=0$ into $n=|(\phi_1, \psi_0)|$.
For $j \not= 1$, taking the inner product of \eqref{t-out-e} at
$t=t_o$ with $\phi_j$ we get
\begin{equation}
\begin{split}
|x_j(t_{o})| & %\leq O(n^3) + n|a| + (1+O(n^2))|z_j(t_{o})| + \norm{\eta}_{L^2_{loc}} + O(n^2)\sum_{l \not= 1, j}|z_l(t_{o})|\\
= O(n^3) + (1+O(n^2))|z_j(t_{o})|, \quad (j \not = 1).
\end{split}
\end{equation}
We also have
\begin{equation}
|x_1(t_{o})| = (\phi_1, Q) + O(n^3) = n(t_c) +  O(n^3) = n +
O(n^{1+2\de}). \end{equation} Since $|z_0(t_o)|=(1+o(1))2 \rho_0$
and $z_H(t_o) \le \sqrt{ {6D}/{\ga_0}} (1+o(1))|z_0(t_o)|$, we have
\eqref{bes.tout}.

Next, we shall prove \eqref{disp.tout}. Denote $\theta_*:= iEt_o -
i\theta(t_{o})$ and
\begin{equation}
x^* = e^{\theta_*}\sum_{j=0}^K x_j(t_{o})\phi_j, \quad \xi^* =
e^{\theta_*}\xi(t_o),\quad \eta^* = \eta(t_o).
\end{equation}
From
\eqref{t-out-e}, we get
\begin{equation} \xi^* = P_c^{H_0}\{ Q + a(t_o) R + \zeta(t_o) + \eta^* - x^* \}. \end{equation}
We write $\xi^* = \xi^*_1 + \xi^*_2 + \xi^*_3$ where
\begin{equation}\begin{split}
\xi^*_1 & := P_c^{H_0}\bigg \{Q + a(t_{o}) R + \sum_{j \not
=1}z_j(t_{o}) \bar{u}_j^+
+ \sum_{j > 1}\bar{z}_j(t_{o})u_j^-  - x^* \bigg \}, \\
\xi^*_2 & := P_c^{H_0}[\bar{z}_0(t_{o}) u_0^-], \quad
 \xi^*_3 : = P_c^{H_0} \eta^*.
\end{split}\end{equation}
From the explicit formulae of $Q, R, u_j^+$, we see that $\xi_1^*$
is localized and $\norm{\xi_1^*} \les n^3 + n|a(t_{o})| + \max_{j
\not = 1} |z_j| n^2 \les n^3$. Therefore, for all $t \geq t_{o}$,
$\tau = t-t_{o}$, by Lemma \ref{uv-est} and Lemma \ref{th2-12}, we
have a uniform constant  $C_6 > \max\{C_3, C_5\}$ such that
\begin{equation} \label{xi*1.est}
\begin{split}
& \norm{e^{-i\tau H_0} \xi_1^*}_{L^2_{loc}} \leq \frac 12 C_6n^3
(1 +\tau)^{-3/2}, \qquad  \norm{e^{-i\tau H_0} \xi_1^*}_{L^p} \leq
\frac 12
C_6n^3 (1 +\tau)^{-\si'}, \\
& \norm{e^{-i \tau H_0} \xi_2^*}_{L^2_{loc}} \leq \frac
12C_6n^{3+\de}(1 + \tau)^{-3/2},\quad  \norm{e^{-i\tau H_0}
\xi_2^*}_{L^p} \leq \frac 12  C_6\rho_0(1 +\tau )^{-\sigma}.
\end{split}
\end{equation}
Here for $\tau <1$ we have used $\norm{e^{-i \tau H_0}
\xi_2^*}_{L^p} \lec \norm{\xi_2^*}_{H^1} \lec |z_0(t_o)|$. Next, we
estimate $e^{-i\tau H_0} \xi_3^*$ in $L^2_{loc}$ and ${L^p}$. Note
$[e^{-i\tau (H_0-E)} \xi_3^*] = e^{\tau J(H_0-E)}[\xi_3^*]$. Recall
that
\begin{equation}
\bL = J(H_0-E)  - W, \quad [\eta^*] = e^{i\theta(t_o)}\eta_{+}^* +
e^{-i\theta(t_o)} \eta_{-}^*, \quad \eta_\pm^* = \eta_{\pm}(t_o),
\end{equation}
for some localized potential $W$ of
order $n^2$. By Duhamel's principle, we have
\begin{equation} \label{xi3.for}
e^{\tau J(H_0-E)}[\xi_3^*] = P_c^{H_0} e^{\tau \bL} [\eta^*] +
\int_{0}^{\tau} e^{J(H_0-E)(\tau -s)} P_c^{H_0} W e^{\bL s} [\eta^*]
ds.
\end{equation}
From Lemma \ref{eta*.OG}, we get
\begin{equation} \label{L5xi3}
\begin{split}
\norm{e^{\tau J(H_0-E)}[\xi_3^*]}_{L^p} & \leq \sum_{\pm }\norm{e^{\tau \bL} \eta_{\pm}^*}_{L^p} + C n^2\sum_{\pm }\int_{0}^{\tau}|\tau -s|^{-\sigma'}\norm{e^{\bL s} \eta_{\pm}^*}_{L^2_{loc}}ds\\
& \leq 2\wt{\La}_1(t) + C n^2\sum_{\pm }\int_{0}^{\tau}|\tau -s|^{-\sigma'} \wt{\La}_2(s + t_o)ds.
%& \leq 3\wt{\La}_1(t).
\end{split}
\end{equation}
Using the fact that
\begin{equation}
\int_0^t(t-s)^{-\beta_1}(\ep^{-1} + s)^{-\beta_2} \leq
C\ep^{\beta_2-1}(\ep^{-1} + t)^{-\beta_1}, \quad 0 < \beta_1 < 1 <
\beta_2,
\end{equation}
we have
\begin{equation}
\begin{split}
n^2\int_0^{\tau}|\tau -s|^{-\sigma'} \wt{\La}_2(s + t_o)ds  &
\le C n^\de\rho_0\wei{t-t_o}^{-\sigma'} + C \rho_0 (\De t + t)^{-\sigma'}
\end{split}
\end{equation}
which is $o(1)\rho_0\wei{t-t_o}^{-\sigma}$. From this and
\eqref{L5xi3}, we get
\begin{equation} \label{Lpxi3}
\norm{e^{\tau J(H_0-E)}[\xi_3^*]}_{L^p} \leq 2\wt{\La}_{1} + o(1)\rho_0\wei{t-t_o}^{-\sigma}.
\end{equation}
Similarly from \eqref{xi3.for} with $\td \al_\I(t)
=\min\{t^{-3/2},\,t^{-9/10}\}$,
\begin{equation}\label{L2xi3}
\begin{split}
\norm{e^{\tau J(H_0-E)}[\xi_3^*]}_{L^2\loc} & \leq \sum_{\pm
}\norm{e^{\tau \bL} \eta_{\pm}^*}_{L^2\loc} +  C n^2\sum_{\pm
}\int_{0}^{\tau}  \td \al_\I(\tau-s)\norm{e^{\bL s}
\eta_{\pm}^*}_{L^2_{loc}} ds
\\
& \leq 2\wt{\La}_2(t)
+ C n^2\sum_{\pm }\int_{0}^{\tau} \td \al_\I(\tau-s) \wt{\La}_2(s + t_o) d s
 \leq 3\wt{\La}_2(t).
\end{split}
\end{equation}
So, \eqref{disp.tout} follows from \eqref{xi*1.est}, \eqref{Lpxi3}, and \eqref{L2xi3}. This completes the proof of Lemma \ref{to.ch}.
\myendproof\\
%=========================================================

For $j \in\{0,1,\cdots, K\}$, let $f_j := |\mu_j(t)|^2$, where
$\mu_j$ is the perturbation of $x_j$ defined in Lemma
\ref{normal-form-u}. Since $\frac{d}{dt}|\mu|^2 = 2\Re
\bar{\mu}\dot{\mu}$ and $c_l^j$ are all purely imaginary, from
\eqref{muj.eqn} we have
\begin{equation} \label{fj.eqn}
\dot{f}_j = \sum_{a,b=0}^K 2(\Re d_{ab}^j) f_a f_b f_j + 2\Re\bar{\mu}_j g_j.
\end{equation}
Let
\begin{equation} \label{fh.def}
f = \sum_{l =1}^{K} f_l,  \quad h = \sum_{l =1}^{K} 2^{-l}f_l, \quad
\gamma \  := \min\{ \gamma_{ab}^0, \ \text{for}\ a, b \geq 1 \} > 0.
\end{equation}
Then, from \eqref{fj.eqn}, Lemma \ref{normal-form-u} and as in \cite[(4.58)]{Tsai}, we have
\begin{equation} \label{ffhh.def}
\begin{split}
\frac{d}{dt}(f_0 + f)(t) & \leq 2(K+1) \max_{l}|\bar{\mu}_l g_l|,\quad \frac{d}{dt}(f_0 + h)(t) \geq -2(K+1) \max_{l}|\bar{\mu}_l g_l|.\\
\end{split}
\end{equation}
Moreover, we have the following lemma.
\begin{lemma} Assume as in Lemma \ref{normal-form-u}. We have
\label{th:7-1}
\begin{equation} \label{hj+.lower}
\begin{split}
\dot f_0  \geq 2\gamma f^2 f_0 +  2 \Re \bar{\mu}_0 g_0,\quad \dot
f  \leq -4\gamma f_0 f^2 + \sum_{l=1}^{K} 2\Re \bar{\mu}_l g_l.
%\frac{d}{dt}(f_j^+) & \geq -4 \gamma_j^+ (f_j^+)^2 f_j^+ + 2\gamma_j (f_j^-)^2 f_j^+ + \sum_{l=0}^{j-1}  \Re \bar{\mu}_l g_l
\end{split}
\end{equation}
\end{lemma}
%================================================
\myproof  From \eqref{fj.eqn} and Lemma \ref{normal-form-u} in
particular \eqref{Redabj}, we have
\begin{equation}
\begin{split}
& \dot{f}_0 - 2\Re \bar \mu_0 g_0  = \sum_{a,b=0}^K
2\Re(d_{ab}^0)f_af_b f_0 = \sum_{a,b=0}^K
[2(2-\delta_a^b)\gamma_{ab}^0 - 4(2-\delta_0^b)\gamma_{0b}^a
]f_af_bf_0.
\end{split}
\end{equation}
Note that  $\gamma^a_{0b} =0$ for any $a$ and  $b$. Thus
\begin{equation}
\begin{split}
\dot{f}_0 - 2\Re \bar \mu_0 g_0 &  = \sum_{a, b=1}^K
2(2-\delta_a^b)\gamma_{ab}^0f_af_bf_0 \geq 2\gamma f^2 f_0.
\end{split}
\end{equation}
This proves the first part of \eqref{hj+.lower}. For the second
part,
\begin{align}
\nn  \dot{f} &- 2 \sum_{l=1}^{K}\Re \bar{\mu}_l g_l =
\sum_{l=1}^{K}\sum_{a,b=0}^{K} 2[(2-\delta_a^b)\gamma_{ab}^l -
2(2-\delta_{l}^b)\gamma_{lb}^a ] f_a f_b f_l , \\ \label{f.ien} &
= \sum_{b=0}^{K}\sum_{a,l=1}^{K} 2[(2-\delta_a^b)\gamma_{ab}^l -
2(2-\delta_{l}^b)\gamma_{lb}^a ] f_a f_b f_l +  \sum_{l=1}^{K}
\sum_{b=0}^{K}  -4(2-\delta_{l}^b)\gamma_{lb}^0 f_0 f_b f_l.
%
%\\ \nn
%& \leq -4  \sum_{l=1}^{K} \sum_{b=1}^{K}  (2-\delta_{l}^b)\gamma_{lb}^0 f_0 f_b f_l -4  \sum_{l=1}^{K}\gamma_{l0}^0 f_0^2 f_l\\
%\nn  & \leq -4 \sum_{b,l=1}^{K} (2-\delta_{l}^b)\gamma_{lb}^0 f_0
%f_b f_l \leq -4\gamma f_0 (f)^2.
\end{align}
By switching $a$ and $l$ in the terms with factor $\gamma_{ab}^l$,
the summands in the first sum become $-2(2-\delta_l^b)\gamma_{lb}^a
f_a f_b f_l \le 0$. The summands of the second sum are also
nonpositive. Keeping only terms with $b>0$ in the second sum, we get
\begin{equation}
\begin{split}
\nn  \dot{f} &- 2 \sum_{l=1}^{K}\Re \bar{\mu}_l g_l \le -4
\sum_{b,l=1}^{K} (2-\delta_{l}^b)\gamma_{lb}^0 f_0 f_b f_l \leq
-4\gamma f_0 f^2.
\end{split}
\end{equation}
This proves the second part of \eqref{hj+.lower}. \myendproof

\medskip

The following proposition estimates the solution in a time interval
containing $[t_o,t_i]$.

\begin{proposition} \label{outside-in}
Let $\de_6(t):= \rho_0^2\wei{t-t_o}^{-\frac{6}{p}}$. For all $t \in
[t_o, t_o + \frac 6\ga n^{-2(2+\de)}]$, we have
\begin{equation} \label{oin.est}
\begin{split}
& \frac{n}{5} \leq \max_{j} |x_j| \le (\tsum_{j=0
}^K|x_j(t)|^2)^{\frac 12} \leq 2n,\\
& \norm{\xi(t)}_{L^2_{loc}} \le  n^{3-\alpha}  + \de_6(t),\quad
\norm{\xi(t)}_{L^p} \le n^{3-\alpha} |t-t_{o}|^{\frac{6-p}{2p}} +
\frac{3}{2}\Lambda_{4}(t).
\end{split}
\end{equation}
%Moreover, at $t = t_i$, we obtain:
%\begin{equation} |x_0(t)| \in [\frac{1}{2}n, \frac{3}{2}n], \quad \max_{j \geq 1}|x_j| \leq n^{1+\delta}. \end{equation}
\end{proposition}

\myproof Since \eqref{oin.est} holds at $t =t_o$, we then prove it
by using the continuity argument. So, we can assume the following
weaker estimates: For $ t_{o} \leq t \leq t_o +\frac 6\ga
n^{-2(2+\de)}$,
\begin{equation} \label{weak-est.u}
\begin{split}
& \frac{n}{10} \leq \max_{j} |x_j| \le (\tsum_{j=0
}^K|x_j(t)|^2)^{\frac 12} \le 3 n, \\
&\norm{\xi(t)}_{L^2_{loc}} \le 2[ n^{3-\al}  + \de_6(t)] \le n^2,\\
& \norm{\xi(t)}_{L^p} \le 2 n^{3-\al} |t-t_{o}|^{\frac{6-p}{2p}} +
3\Lambda_{4}(t) \le n^{2.7} + 3\Lambda_{4}(t).
\end{split}
\end{equation}
In particular $\norm{\xi(t)}_{L^2_{loc}}+\norm{\xi(t)}_{L^p}\ll n$.
The proof of Proposition \ref{outside-in} follows from Lemma
\ref{xi.in.est} and Lemma \ref{bst.est} below. \myendproof
%============================================================================================================
%============================================================================================================
\begin{lemma}\label{xi.in.est} For all $t \in [t_o, t_o + \frac 6\ga n^{-2(2+\de)}]$, we have
\begin{equation} \norm{\xi(t)}_{L^2_{loc}} \le  n^{3-\alpha}  + \de_6(t), \quad \norm{\xi(t)}_{L^p} \le n^{3-\alpha} |t-t_{o}|^{\frac{6-p}{2p}} +   \frac{3}{2}\Lambda_{4}(t). \end{equation}
\end{lemma}
%============================================================================================================
%============================================================================================================
\myproof For all $t -t_{o} \leq C n^{-2(2+ \delta)}$, by \eqref{weak-est.u}, we have
\begin{equation}\label{L5xi.est}
\norm{\xi(t)}_{L^p} \les n^{3-\al -\frac{2(2+\de)(6-p)}{2p}} + \Lambda_{4}(t) \leq C[n^{\frac{(5+\de)p - 6(2+\de)}{p} - \al} + 3\La_4(t)].% \quad \Lambda(t) \les n^{2.9}.
\end{equation}
We have
\begin{equation} \xi(t) = e^{-iH_0(t-t_o)} \xi(t_{o}) + \int_{t_{o}}^t e^{-iH_0(t-s)} \Pc i^{-1} G(s) ds. \end{equation}
So, we have %From the linear estimates of the semi group $e^{-iH_0t}$ and the assumption of Proposition \ref{outside-in}, we have
\begin{equation} \label{L4.xi}
 \norm{\xi(t)}_{L^p} \leq \Lambda_4(t) + C \int_{t_{o}}^t |t -s|^{-\frac{3(p-2)}{2p}} \norm{G(s)}_{L^{p'}} ds.
\end{equation}
Note that $\norm{G}_{L^{p'}} \les \norm{G_3}_{L^{p'}} + \norm{G -
G_3 - \kappa  \xi^2\bar{\xi}}_{L^{p'}} + \norm{\kappa
\xi^2\bar{\xi}}_{L^{p'}}$ and $\norm{G_3}_{L^{p'}} \les n^3$. On
the other hand, from Lemma \ref{nonlin.est}, \eqref{weak-est.u}
and \eqref{L5xi.est}, we get
\begin{equation} \label{localG.term}
\begin{split}
& \norm{G - G_3 - \kappa  \xi^2\bar{\xi}}_{L^1 \cap L^{p'}} \les
n^2 \norm{\xi}_{L^2_{loc}} \les [n^{5-\alpha} + n^2\de_6(t)].
\end{split}
\end{equation}
On the other hand, using H\"{o}lder's inequality, we get
\begin{equation} \label{Lp.xicube}
\norm{\kappa  \xi^2\bar{\xi}}_{L^{p'}} \leq
\norm{\xi}_{L^2}^{\frac{2(p-4)}{p-2}}
\norm{\xi}_{L^p}^{\frac{p+2}{p-2}}, \quad \norm{|\xi|^2\xi}_{L^1}
\leq
\norm{\xi}_{L^2}^{\frac{2(p-3)}{p-2}}\norm{\xi}_{L^p}^{\frac{p}{p-2}}.
 \end{equation}
From this, \eqref{L5xi.est} and since $0 < \de \le \frac{1}{10}$, we
get
\begin{equation} \label{nlocalG.term}
\begin{split}
& \norm{\kappa  \xi^2\bar{\xi}}_{L^{p'}}\leq
\norm{\xi}_{L^2}^{\frac{2(p-4)}{p-2}}
\norm{\xi}_{L^p}^{\frac{p+2}{p-2}} \leq o(1)[n^{5 -2\al} +
\La_4(t)^{\frac{p +2}{p-2}}].
\end{split}
\end{equation}
By \eqref{La34p.est}, \eqref{localG.term}, and \eqref{nlocalG.term}, we have
\begin{equation} \label{G.Lpd}
\norm{G(s)}_{L^{p'}} \leq C[n^3 + o(1)\wt{\de}_2(t)], \quad \wt{\de}_2(t):= [\rho_0\wei{t-t_o}^{-\sigma}]^{\frac{p+2}{p-2}}.
\end{equation}
Therefore, using $\si \frac{p+2}{p-2}>1$,
\begin{equation}
\begin{split}
\norm{\xi(t)}_{L^p} & \leq  \Lambda_{4}(t) + C \int_{t_{o}}^t |t - s|^{-(\frac{3}{2} - \frac{3}{p})}[n^3 + o(1)\wt{\de}_2(s)] ds \\
& \leq Cn^3 |t -t_{o}|^{\frac{6-p}{2p}} +  \frac{3}{2}\Lambda_{4}(t).
\end{split}
\end{equation}
So, we have proved the estimate of $\norm{\xi(t)}_{L^p}$.

We now estimate $\norm{\xi(t)}_{L^2_{loc}}$. By \eqref{xi2},
\eqref{xi.all}, \eqref{weak-est.u} and Lemma \ref{to.ch}, we have
\begin{equation}\label{L2.xi3.12}
\norm{\xi_1^{(3)}(t)}_{L^2_{loc}} \leq  \Lambda_{3}(t), \quad
\norm{\xi_2^{(3)}(t)}_{L^2_{loc}} \les n^3(1+t -t_o)^{-3/2}.
\end{equation}
By \eqref{xi.all} and the estimate of $\max_{j}|\dot{u}_j|$ in Lemma
\ref{nonlin.est}, we get
\begin{equation} \label{L2.xi3.3}
\norm{\xi_3^{(3)}(t)}_{L^2_{loc}} \les
\int_{t_{o}}^{t}|1+t-s|^{-3/2} n^5 ds \les n^5.
\end{equation}
For $\xi_4^{(3)}(t)$, bounding its integrand by
either $L^\infty$ or $L^p$-norm and using \eqref{localG.term}, we
have
\begin{equation}  \label{L2.xi3.4}
\begin{split}
\norm{\xi_4^{(3)}(t)}_{L^2_{loc}} & \les
\int_{t_{o}}^t\min\{|t-s|^{-3/2}, |t-s|^{-\frac{3(p-2)}{2p}} \}
\norm{G-G_3 - \kappa \xi^2\bar{\xi}}_{L^1\cap L^{p'}} ds
\\
& \les \int_{t_{o}}^t\min\{|t-s|^{-3/2}, |t-s|^{-\frac{3(p-2)}{2p}}\}
[n^{5-\al} + n^2\de_6(s)]ds\\
& \les n^{5-\al} + n^2\de_6(t).
\end{split}
\end{equation}
For $\xi_5^{(3)}(t)$, bounding its integrand in either
$L^{\frac{2p}{p-4}}$ or $L^p$, we have
\begin{equation}
\norm{\xi_5^{(3)}(t)}_{L^2_{loc}} \leq C
\int_{t_{o}}^t\min\{|t-s|^{-\frac{6}{p}}, |t-s|^{-\frac{3(p-2)}{2p}}
\}\norm{|\xi|^2\xi}_{L^{\frac {2p}{p+4}} \cap L^{p'}} ds.
\end{equation}
By \eqref{Lp.xicube}, $\frac {p+2}{p-2}>2$ and $ 2 < \frac{6p}{p+4}
< p$ because $\frac{27}{5} < p <6$,
\begin{equation}
\label{l-choice} \norm{|\xi|^2\xi}_{L^{\frac{2p}{p+4}}\cap L^{p'}}
\leq C \norm{\xi}_{L^p\cap L^2}\norm{\xi}_{L^p}^{2}\leq
o(1)\norm{\xi}_{L^p}^{2}.
\end{equation}
Therefore, by \eqref{weak-est.u},
\begin{equation} \label{L2loc.xi35}
\begin{split}
\norm{\xi_5^{(3)}(t)}_{L^2_{loc}} & \leq o(1)\int_{t_{o}}^t
\min\{|t-s|^{-\frac{6}{p}}, |t-s|^{-\frac{3(p-2)}{2p}} \}
[n^{5.4} + \La_4^{2}(s )] ds\\
& \leq o(1)[n^{5.4} + \La_4(t)^{2} + \de_7(t)],
\end{split}
\end{equation}
where
\begin{equation} \label{delta7.def}
 \de_7(t) :=  \rho_0^2\wei{t-t_o}^{-\frac{6}{p}}
 + n^{\frac{-2+2\de }{3}}(n^{-4} +t -t_o)^{-\frac{6}{p}},
\end{equation}
and we have used $\frac{2}{3} < \si < \frac{3}{4}$, \eqref{La4.dec},
and \eqref{int-est2} with $a=6/p< b = 2 \si -2\al$, (or $b=2\si$).

Collecting all of the estimates of $\xi_j^{(3)}$ with $j =1,2,3,4$,
we have
\begin{equation} \label{xi3.est}
\begin{split}
& \norm{\xi^{(3)}(t)}_{L^2_{loc}} \leq \Lambda_{3}(t) + Cn^{5}
+o(1)[\La_4(t)^{2} +\de_7(t)].
\end{split}
\end{equation}
By \eqref{La34.def}, we have $\La_3(t) \les n^3$ and $\La_4(t)^{2}
+\de_7(t) \leq n^3 + \de_6(t)$. Thus
\begin{equation}
\norm{\xi(t)}_{L^2_{loc}} \leq \norm{\xi^{(2)}(t)}_{L^2_{loc}} +
\norm{\xi^{(3)}(t)}_{L^2_{loc}} \leq Cn^3 +  o(1)\de_6(t) .
\end{equation}
This completes the proof of the lemma.
\myendproof
%==================================
%=====================================================================================
\begin{lemma}\label{err.g}
For $t \in [t_o, t_o + \frac 6\ga n^{-2(2+\delta)}]$, the error
terms $g_j(t)$ in \eqref{muj.eqn} satisfy
\begin{equation} \label{gj-out-in.est}
|g_j(t)| \leq o(1)n^{6.7 +\de} + Cn^2g(t),
\end{equation}
where
\begin{equation}
g(t) := \La_3(t) + o(1)[n^{1+3\de}\wei{t-t_o}^{-\frac{p\sigma}{p-2}}
+ \La_4^2(t) + \de_7(t)]
\end{equation}
satisfies
\begin{equation} \label{OG.cond}
\begin{split}
& \int_{t_o}^\I g(s) ds \leq o(1)n^{-\frac{2}{3}} ; \quad g(t) \leq
o(1)n\rho_0^2, \ \forall\ t \geq t_o + n^{-3}.
\end{split}
\end{equation}
\end{lemma}
%===================================================================
\myproof Recall \eqref{error.est.mu},
\begin{equation} |g_j(t)| \les n^7 + n^2 \norm{\xi^{(3)}}_{L^2_{loc}} + n \norm{\xi}_{L^2_{loc}}^2 + \norm{\xi}_{L^2_{loc}}^{\frac{2(p-3)}{p-2}}\norm{\xi}_{L^p}^{\frac{p}{p-2}}.\end{equation}
From \eqref{weak-est.u} and \eqref{xi3.est}, we get
\begin{equation}
\begin{split}
n^2 \norm{\xi^{(3)}}_{L^2_{loc}} &\leq n^2\La_3 + C n^{7} +
o(1)n^2[\La_4^{2} + \de_7],
\\
n\norm{\xi}_{L^2_{loc}}^2 &\leq C[n^{7-2\al} + n\de_6(t)^2],
\end{split}
\end{equation}
and, using $[n^{2.7} +\La_{4,1} +\La_{4,3}]^{\frac{p}{p-2}}\le o(1)
n^{\frac{5+3\de}{2}}$,
\begin{equation}
\begin{split}
 \norm{\xi}_{L^2_{loc}}^{\frac{2(p-3)}{p-2}}
 \norm{\xi}_{L^p}^{\frac{p}{p-2}}
 & \les [n^{3-\al} + \de_6(t)]^{\frac{2(p-3)}{p-2}}
 [n^{2.7} + \La_4]^{\frac{p}{p-2}} \\
\ & \leq o(1)[n^{\frac{2(p-3)(3-\al)}{p-2}}
+ \de_6^{\frac{2(p-3)}{p-2}}]
[n^{\frac{5+3\de}{2}} + \rho_0^{3/2}\wei{t-t_o}^{-\frac{p\sigma}{p-2}}]\\
& \leq o(1)[ n^{6.7 +\de}  +
\rho_0^3\wei{t-t_o}^{-\frac{p\sigma}{p-2}}].
\end{split}
\end{equation}
Summing the estimates we get \eqref{gj-out-in.est}. The estimates
\eqref{OG.cond} follow from direct checking. \myendproof

%=====================================================================================
\begin{lemma} \label{bst.est}
For all $t \in [t_o, t_o + \frac 6\ga n^{-2(2+\delta)}]$, we have
\begin{equation}
\frac{1}{5}n  \leq \max_{j}|x_j(t)| \leq (\tsum_{j=0
}^K|x_j(t)|^2)^{\frac 12} \leq 2 n.
\end{equation}
\end{lemma}
%===================================================================================================
\myproof
From the first equation of \eqref{ffhh.def}, \eqref{OG.cond} and
$\de \leq \frac{1}{10}$, we get
\begin{equation} \label{u.max}
\begin{split}
(f_0 + f)(t) & \leq (f_0 + f)(t_o) + C n\max_{j } \int_{t_o}^{t} |g_j(s)| ds\\
&  \leq (f_0 + f)(t_o) + C[o(1)n^{7.7 +\de}(t-t_o) +  n^3\int_{t_o}^t g(s)ds] \\
& \leq (f_0 + f)(t_o) + o(1)\rho_0^2 \leq [1+o(1)](f_0 + f)(t_o).
\end{split}
\end{equation}
By \eqref{error.est.mu}, \eqref{weak-est.u}, we have $[1-o(1)]\sum_j
{|x_j|^2}\le f_0  +f$. By Lemma \ref{to.ch}, we get $(f_0 + f)(t_o)
\leq 2n^2$. It follows from \eqref{u.max} that $(\sum_{j=0
}^K|x_j(t)|^2)^{\frac 12} \leq 2 n$.

Similarly, by integrating the second equation of \eqref{ffhh.def},
we obtain
\begin{equation} (f_0 + h)(t) \geq [1-o(1)](f_0 + h)(t_o). \end{equation}
By \eqref{error.est.mu}, \eqref{weak-est.u} and the definition of $f_0, h$, we get
\begin{equation}
(f_0 + h)(t) \leq [\tsum_{k=0}^K 2^{-k} +o(1)] \max_j |x_j(t)|^2.
\end{equation}
Therefore,
\begin{equation}
2\max |x_j(t)|^2 \geq [1-o(1)](f_0 + h)(t_o) \geq [1-o(1)]
\frac{1}{2}|x_1(t_o)|^2.
\end{equation}
Hence $\max_{j} |x_j(t)|^2 \geq \frac{n^2}{25}$ for all $t \in [t_o,
t_o +\frac 6\ga n^{-2(2+\de)}]$.
%This completes the proof of Lemma \ref{bst.est}.
\myendproof
%==============================================================================================================
\begin{proposition} \label{ti.lemma}
There exists $t_i$ such that $t_o + \frac{\de}{10 \td \ga}
n^{-4}\log\frac{1}{n} \leq t_i \leq t_o +
\frac{7}{\gamma}n^{-4-2\delta}$ and
\begin{equation}\label{ti.lemma.eq1}
\frac{n}{5} \leq |x_0(t_i)| \leq 2n, \quad (0.9)\rho_0 \leq
(\sum_{j=1}^{K}|x_j(t_i)|^2 )^{1/2} \leq (1.1)\rho_0.
\end{equation}
Above $\wt{\gamma} = \max\{1,(d_{ab}^l)_-: \ \forall a,b,l
=0,\ldots, K\}$ and $d_{ab}^l=O(1)$ are  given in \eqref{Redabj}.
\end{proposition}

\myproof By Lemma \ref{bst.est}, we already have $|x_0| \leq 2n$.
The proof is divided into four steps.

\medskip
%===============================================================================================================
\noindent\textbf{Step 1:} Let $t_1 := t_{o} + n^{-3}$. For $t_{o}
\leq t \leq t_1$, for any $j$, by \eqref{fj.eqn},
\eqref{weak-est.u}, \eqref{gj-out-in.est}, and \eqref{OG.cond}, we
get
\begin{equation} \label{in.to} \begin{split}
|f_j(t) - f_j(t_{o})| & \les \int_{t_{o}}^{t_1} [ n^6 + n|g_j(s)|]
ds  \les n^3 + n \int_{t_{o}}^{t_1} [n^2g(s)]ds \leq o(1)\rho_0^2.
\end{split}\end{equation}
In particular, for $j =0,1$, we get
\begin{equation} [1-o(1)]f_j(t_o) \leq f_j(t) \leq [1 + o(1)]f_j(t_o), \quad \forall \ t \in [t_o, t_1]. \end{equation}
By \eqref{error.est.mu} and the definitions of $f_j$, we get
\begin{equation} [1-o(1)]|x_j(t_o)| \leq |x_j(t)| \leq [1+o(1)]|x_j(t_o)|, \ \forall \ t \in [t_o, t_1],\ j = 0,1. \end{equation}
Together with \eqref{bes.tout}, for $t \in [t_o, t_1]$, we have
\begin{equation} \label{inlayer.est}
1.8\rho_0 \leq |x_0(t)| \leq 2.2 \rho_0, \quad 0.8n \leq |x_1(t)|
\leq 1.2 n.
\end{equation}
On the other hand, for $j >1$, from \eqref{in.to}, we obtain $f_j(t)
\leq f_j(t_o) + o(1) \rho_0^2$ for $t \in [t_o, t_1]$. So, by
\eqref{error.est.mu}, \eqref{bes.tout}, and the definition of $f_j$,
we get
\begin{equation} \label{xHin}
|x_j(t)| \leq [1+o(1)]f_j(t)^{1/2} \leq 7 \sqrt {\frac D
{\ga_0}}\rho_0, \quad \forall \ t \in [t_o, t_1], \ \forall \ j >1.
\end{equation}
\medskip \noindent
\textbf{Step 2:} Let us define
\begin{equation} t_2 : = \sup \{t \geq t_1 : f_0(s) < \frac{n^2}{10}, \ \forall \ s \ \in [t_1,t] \}.\end{equation}
By \eqref{inlayer.est}, $t_2 < t_1$. We shall prove that
\begin{equation} \label{t1.est}
t_1 < t_2 \leq t_2': = t_1 +  a^{-1}\log \frac{n^2}{5f_0(t_1)},
\quad a := 2\gamma [\frac{n^2}{50}]^2.
\end{equation}
For all $t_1 \leq  t \leq t_2$, $f_0(t) < \frac{n^2}{10}$. Note
$h(t_1) \ge f_1(t_1)/2 \ge (1+o(1))(0.8n)^2/2 \ge (0.3)n^2$. From
\eqref{ffhh.def} and Lemma \ref{err.g}, we get
\begin{equation} \label{t1.h} \begin{split}
h(t) & \geq (f_0 + h)(t_1) - f_0(t)
- 2(K+1) \int_{t_1} ^{t} \max_{j} |\mu_j||g_j|(s) ds
\\
\ & \geq (0.3)n^2 - \frac{n^{2}}{10} - Cn\int_{t_1}^t[n^{6.7 + \de}
+ n^2 g(s)]ds \geq \frac{n^2}{100}.
\end{split}
\end{equation}
By \eqref{hj+.lower}, \eqref{weak-est.u} and \eqref{t1.h}, we have, for $t \in [t_1, t_2']$,
\begin{equation} \label{hk+.moi}
\dot f_0 \geq 2\gamma f^2 f_0 - 2|\mu_0||g_0| \geq 2\gamma_0 (2h)^2 f_0 - 4n|g_0|  \geq 2\gamma [\frac{n^2}{50}]^2 f_0 - 4n|g_0|.
\end{equation}
Note the coefficient of $f_0$ is $a$. Thus
\begin{equation}
f_0(t) \geq e^{a(t-t_1)}[f_0(t_1) -
4n\int_{t_1}^te^{-a(s-t_1)}g_0(s) ds].
\end{equation}
On the
other hand, from \eqref{gj-out-in.est}, we have
\begin{equation}
\begin{split}
 n\int_{t_1}^te^{-a(s-t_1)}g_0(s) ds & \leq n\int_{t_1}^t[n^{6.7 + \de} + n^2 g(s)]ds \\
 & \leq n^{7.7 + \de}(t-t_1) + n^3\int_{t_1}^t g(s) ds \leq o(1)\rho_0^2 \leq o(1) f_0(t_1).
\end{split}
\end{equation}
Therefore,
\begin{equation} \label{f0.t1}
f_0(t)  \geq \frac{1}{2}e^{a(t-t_1)}f_0(t_1), \quad \forall t \in
[t_1, t_2].
\end{equation}
This shows $t_2 \le t_2'$ is finite, and $f_0(t_2) =
\frac{n^2}{10}$.

%================================================================================================================
\medskip

\noindent\textbf{Step 3:} Define
\begin{equation} \label{ti.def}
t_{i} := \sup \{t \geq t_2 : f (s) > \rho_0^{2}, \ \forall \ s \in
[t_2, t)\}.
\end{equation}
From \eqref{t1.h}, we get $t_i
> t_2$. We shall prove in Steps 3 and 4 that
\begin{equation}
t_2 + \frac{\de}{10 \td \ga} n^{-4}\log\frac{1}{n} \leq t_{i} \leq
t_3 := t_2 + \frac{6}{\gamma} n^{-4 - 2\de}.
\end{equation}
%Suppose by contradiction that $t_i > t_3$. So, $t_2 < t_3 < t_i$.
By
definition of $t_i$, we get
\begin{equation} \label{ass.ti}
f(t) > \rho_0^2, \quad \forall\ t \in [t_2, t_i).
\end{equation}
From Lemma \ref{th:7-1} and \eqref{ass.ti}, we have
\begin{equation} \frac{d}{dt}(f_0(t))
\geq 2\gamma \rho_0^4 f_0(t)- 4n|g_0|, \ \forall \ t \in [t_2, t_i).
\end{equation}
From this and as in \eqref{f0.t1}, we also obtain
\begin{equation} \label{f0.ti}
 f_0(t) \geq \frac{1}{2}e^{2\gamma \rho_0^4(t-t_2)}f_0(t_2)
 \geq \frac{n^2}{20}, \ \forall \ t \in [t_2, t_i).
 \end{equation}
From this, \eqref{hj+.lower}, and Lemma \ref{err.g}, for $ t \in
[t_1, t_i)$,
\begin{equation}
\begin{split}
\frac{d}{dt}(f (t)) & \leq -4\gamma f_0(t) f(t)^2 + C n \max_{k>0}
|g_k| \\
%& \leq -4\gamma f_0(t) f(t)^2 + Cn[n^{6.7 + \de} + n^2 g(t)] \\
\ & \leq -\frac{\gamma n^{2}}{5} f (t)^2 + C n[n^{6.7+\de} + n^2g(s)].
\end{split}\end{equation}
From this and \eqref{ass.ti}, (and $\de \leq \frac{1}{10}$), we get
\begin{equation} \label{eq6.79}
\frac{n^2\gamma}{6} - Cn^{3}\rho_0^{-4} g(t) < \frac{n^2\gamma}{5}
-\frac{Cn[n^{6.7+ \de} + n^2g(s)]}{f^2} \leq - \frac{\dot{f}}{f^2},
\quad \forall t \in [t_2, t_i).
\end{equation}
Note that by
\eqref{OG.cond}, \eqref{t1.h}, Proposition \ref{outside-in} and $\de
\leq \frac{1}{10}$, we have $\forall \ t \geq t_2$
\begin{equation}
n^{-1-4\delta} \int_{t_2}^t g(s) ds
\leq o(1)n^{-1-4\de}n^{-2(1 - \de)/3} =o(1)n^{-\frac{5(1+2\de)}{3}}
\leq o(1) f(t_2)^{-1}.
\end{equation}
Integrating \eqref{eq6.79} in $[t_2,t]$, we get
\begin{equation} \label{ub.f}
f (t) < [f (t_2)^{-1}/2 + \frac{n^2\gamma}{6} (t - t_2)]^{-1}, \quad
\forall \ t \in [t_2, t_i].
\end{equation}
In particular, $\rho_0^2 < f(t) < [\frac{n^2\gamma}{6} (t -
t_2)]^{-1}$, which shows $t_i \le t_3$, and $f(t_i) = \rho_0^2$.
From this, \eqref{error.est.mu} and \eqref{f0.ti}, we get the
estimates \eqref{ti.lemma.eq1}. Since
\begin{equation} %\label{ub.f}
t_i -t_o\le (t_i - t_2) + (t_2-t_1)+(t_1-t_o)\le \frac 6\ga
n^{-4-2\de} +C n^{-4}\log \frac 1n +n^{-3}
\end{equation}
by \eqref{t1.est} and \eqref{inlayer.est}, we get the upper bound of
$t_i -t_o$ in Prop.~\ref{ti.lemma}.

\medskip

\noindent{\bf Step 4:} It remains to show that $t_i \ge t_2 +
\frac{\de}{10 \td \ga} n^{-4}\log\frac{1}{n}$. Recall $g(t) \leq
o(1)n\rho_0^2$ for all $t \geq t_1= t_o + n^{-3}$ from Lemma
\ref{err.g}. By \eqref{fj.eqn} and Prop.~\ref{outside-in},
\begin{equation}
\dot{f} (t) \geq -9\wt{\gamma} n^4 f (t) - Cn[n^{6.7+\de} +n^2g(t)]
\geq -10\wt{\gamma} n^4 f(t), \ \forall \ t \in [t_1, t_{i}],
\end{equation}
where $\wt{\gamma} = \max\{1,(d_{ab}^l)_-: \ \forall a,b,l
=0,\ldots, K\}$. This implies that
\begin{equation}
t_{i}- t_2 \geq \frac{n^{-4}}{10\wt{\gamma}} \log \frac{f (t_2)}{f
(t_{i})} \geq \frac{\de}{10 \td \ga} n^{-4}\log\frac{1}{n}.
\end{equation}
For the second inequality we have used $f(t_2)\ge h(t_2) \ge n^2/50$
by \eqref{t1.h}. This completes the proof of Proposition
\ref{ti.lemma}. \myendproof

\medskip

At $t=t_i$ the solution enters $\rho_0$-neighborhood of ground
states and we change to linearized coordinates. For that purpose we
prepare outgoing estimates at $t=t_i$.

%============================================================
\begin{lemma} \label{OGxi.est} Let $t_i$ be as in Proposition
\ref{ti.lemma}. For any $t > t_i$, we have
\begin{equation}
\begin{split}
& \norm{e^{-iH_0(t-t_i)}\xi(t_i)}_{L^2\loc} \leq\frac{1}{2}[
\La_{L,1}(t) +\La_{L,2}(t)],
\\
& \norm{e^{-iH_0(t-t_i)}\xi(t_i)}_{L^p} \leq
\frac{1}{2}[\La_{G,1}(t) + \La_{G,2}(t)],
\end{split}
\end{equation}
where for some constant $C_7 \ge C_6$ and $\sigma'
=\frac{3(p-2)}{2p}$,
\begin{equation} \label{LaLG1.def}\begin{split}
& \La_{L,1}(t): = 2C_7[n^{-1+2\de}\wei{t-t_o}^{-7/6}
+ \rho(t)^3 + n^{4/5}\rho(t)^{7/3}],\\
& \La_{G,1}(t): = 2C_7[n^{-1+\de}\wei{t-t_o}^{-\sigma} + n^{-1
+2(2+\de)\al}(\Delta t +t)^{-\sigma +\al}],
\\
&\La_{L,2}(t): = \frac{2n^{\frac{5p -18 +
p\de}{p-2}}(t_i-t_o)}{t-t_o}\wei{t-t_i}^{-1/2}, \\
&\La_{G,2}(t): = 2C_7n^{3}(t_i-t_o)(t-t_o)^{-\sigma'}.
\end{split}
\end{equation}
\end{lemma}

\medskip
%=================================================
\myproof Decompose $ e^{-i(t-t_i) H_0} \xi(t_{i})=\chi(t) + J(t)$,
where
\begin{equation}
\chi(t):= e^{-i(t-t_{o})H_0} \xi(t_{o}) ,\quad J(t) :=
\int_{t_{o}}^{t_{i}} e^{-i(t-s)H_0} \Pc G(s) ds.
\end{equation}
Denote $T= t_{i} - t_{o}$. By Lemma \ref{to.ch} and using
$n^{-4}\log \frac 1n \lec T \lec n^{-2(2+\delta)}$, we have
\begin{equation} \label{ti.OG}
\norm{\chi(t)}_{L^p} \leq \La_4(t)\le \frac 12 \La_{G,1}(t), \quad
\norm{\chi(t)}_{L^2\loc} \leq \La_3(t)\le \frac 12 \La_{L,1}(t),
\end{equation}
for some $C_7$. By \eqref{G.Lpd}, we have
\begin{equation}
\norm{G(s)}_{L^{p'}} \leq C[n^3 + o(1)\wt{\de}_2(s)] \quad \forall \
s \in [t_{o}, t_{i}], \quad \wt{\de}_2(s) =
[\rho_0\wei{s-t_o}^{-\sigma}]^{\frac{p+2}{p-2}}.
\end{equation}
 So, we have (using $\frac{p+2}{p-2}>2$)
\begin{equation}
\begin{split}
\norm{J(t)}_{L^{p}} & \leq C \int_{t_{o}}^{t_{i}} |t -s|^{-\sigma'} \norm{G(s)}_{L^{p'}} ds \leq C \int_{t_{o}}^{t_{i}} |t -s|^{-\sigma'}[n^3  +o(1)\wt{\de}_2(s)]ds \\
&  \leq C n^3T( t-t_o)^{-\sigma'} + \rho_0^2(t-t_o)^{-\sigma'} \leq
\frac{1}{2}\La_{G,2}(t).
\end{split}
\end{equation}

It remains to estimate $\norm{J(t)}_{L^2\loc}$. By
\eqref{localG.term} and \eqref{Lp.xicube},
\begin{equation}
\norm{G(s)}_{L^1\cap L^{p'}} \leq Cn^{3} + C n^2
\norm{\xi(s)}_{L^2_{loc}} +o(1)\norm{\xi(s)}_{L^p}^{\frac{p}{p-2}} .
\end{equation}
By \eqref{weak-est.u} and \eqref{La34p.est},
\begin{equation}
\norm{G(s)}_{L^1\cap L^{p'}} \leq  o(1)[n^{\frac{5p -18 +
p\de}{p-2}} + \rho_0^{3/2}
\wei{s-t_o}^{-\frac{p\sigma}{p-2}}].
\end{equation}
Thus
\begin{equation}
\begin{split}
\norm{J(t)}_{L^2\loc} & \leq C\int_{t_o}^{t_i} \min\{(t-s)^{-3/2}, (t-s)^{-\sigma'}\} \norm{G(s)}_{L^1\cap L^{p'}} ds \\
& \leq o(1) \int_{t_o}^{t_i} \min\{(t-s)^{-3/2}, (t-s)^{-\sigma'}\}[n^{\frac{5p -18 + p\de}{p-2}} +  \rho_0^{3/2} \wei{s-t_o}^{-\frac{p\sigma}{p-2}}] ds\\
& \leq o(1) n^{\frac{5p -18 +
p\de}{p-2}}\frac{T}{t-t_o}\wei{t-t_i}^{-1/2} + o(1)
\rho_0^{3/2}(t-t_o)^{-1}\wei{t-t_i}^{-1/2} ,
\end{split}
\end{equation}
which is bounded by $\frac{1}{2}\La_{L,2}(t)$. This  completes the
proof of the lemma. \myendproof
%% ===================================================================================================================

%% ===================================================================================================================
\section{Converging to a ground state}
%==============================================================
In this section we study the solution when it is already inside a
neighborhood of the ground states. It is similar to the estimates in
\cite{BP2,TY1,C2,Tsai}, however,  it requires a proof because the
dispersive component has much worse estimates. As in Section 4, for
fixed $T \ge t_i$ we shall decompose $\psi(t)$ as (see
\eqref{psi.dec})
\begin{equation} \label{gs.dec}
\psi(t) = [Q_{0, n(T)} + a(t) R_{0, n(T)} +\zeta(t) +
\eta(t)]e^{-iEt + i\theta(t)}, \quad t \in [t_i, T].
\end{equation}
We have $a(T) = 0$, and
\begin{equation}
\zeta = \sum_{j=1}^K\zeta_j, \quad \zeta_j = \bar{z}_j u_j^{-} +
z_ju_j^+, \quad [\eta] =\begin{bmatrix}\Re \eta \\ \Im \eta
\end{bmatrix} = e^{i\theta}\eta_+ + e^{-i\theta}\eta_-.
\end{equation}
Denote $z_H(t) = (\sum_{j=1}^{K}|z_j(t)|^2)^{1/2}$. From Lemma
\ref{Linearized.dec} and Proposition \ref{ti.lemma}, \eqref{gs.dec}
is valid at least  for $T > t_i$ sufficiently close to $t_i$. We
prove in this section that this is true with suitable estimates for
all $T \geq t_i$ and, moreover, $n(T)$ converges to some $n_+ \sim
n$ as $T \rightarrow \infty$.
%=======================================
\begin{lemma} \label{OG.ti.ln}
There exists $C_8>0$ such that if $T> t_i$ and $n(T)/n(t_i) \in
(\frac 12, \frac 32)$, then
\begin{equation}\label{th71-eq1}
\frac45\rho_0 \leq z_H(t_i) \leq \frac 65\rho_0,
\end{equation}
and, for $t \ge t_i$,
\begin{equation}
\begin{split}
& \norm{e^{\bL(t-t_i)}\eta_\pm(t_i)}_{L^2\loc} \leq \La_L(t): =
\La_{L,1}(t) + \La_{L,2}(t) + \La_{L,3},
\\
& \norm{e^{\bL(t-t_i)}\eta_\pm(t_i)}_{L^p} \leq \La_G(t): =
\La_{G,1}(t) + \La_{G,2}(t) + \La_{G,3},
\end{split}
\end{equation}
where $\La_{L,1}$, $\La_{L,2}$, $\La_{G,1}$ and $\La_{G,2}$ are
defined in Lemma \ref{OGxi.est},
$\La_{L,3}(t)=C_8n^3\wei{t-t_i}^{-3/2}$, and
$\La_{G,3}(t)=C_8n^3\wei{t-t_i}^{-\si'}$.
\end{lemma}

\myproof At $t=t_i$, with $Q=Q_{n(T)}$ and $\Th= E_{n(T)}t_i -
\theta(t_i)$ we have
\begin{equation} \label{gs.ti}
Q + a(t_i) R +\zeta(t_i) + \eta(t_i) = e^{i\Th}[\sum_{j=0}^{K}
x_j(t_i)\phi_j +\xi(t_i)].
\end{equation}
For each $j \geq 1$, applying the projection $P_j$ (see
Prop.~\ref{L-spectral}, (iii)) to \eqref{gs.ti}, we get $z_j(t_i) =
e^{i\Th} x_j(t_i) + O(n^3)$. By Proposition \ref{ti.lemma} we get
\eqref{th71-eq1}. Denote
\begin{equation}
\eta_{1}  :=\sum_{j=0}^{K} e^{i\Th} x_j(t_i)\phi_j - Q - a(t_i) R
-\zeta(t_i), \quad  \eta_{2}  : = e^{i\Th} \xi(t_i).
\end{equation}
Then $[\eta] = P_c^{\bL}[\eta_1] + P_c^{\bL}[\eta_2]$. Since
$\eta_1$ is localized and of order $O( n^3)$, we get
\begin{equation}
\norm{e^{\bL(t-t_i)}P_\pm[\eta_1]}_{L^2\loc} \leq
C_8n^3\wei{t-t_i}^{-3/2}, \quad
\norm{e^{\bL(t-t_i)}P_\pm[\eta_1]}_{L^p} \leq
C_8n^3\wei{t-t_i}^{-\sigma'}.
\end{equation}
On the other hand, we have $\bL = J(H_0 -E) + W$ with $W = O(n^2)$
which is localized. By Duhamel's formula,
\begin{equation}
e^{\bL(t-t_i)}P_{\pm}[\eta_2] = P_{\pm} e^{J(H_0-E)(t-t_i)}[\eta_2]
+ \int_{t_i}^t P_{\pm} e^{\bL(t-s)} W e^{J(H_0-E)s}[\eta_2] ds.
\end{equation}
Thus, using Lemma \ref{OGxi.est},
$\norm{e^{\bL(t-t_i)}P_{\pm}[\eta_2]}_{L^p}$ is bounded by
\begin{flalign*}
 & \leq \norm{e^{J(H_0-E)(t-t_i)}[\eta_2]}_{L^p} + Cn^2\int_{t_i}^t |t-s|^{-\sigma'}\norm{e^{J(H_0-E)s}[\eta_2]}_{L^2\loc}ds\\
& \leq\frac{1}{2}[\La_{G,1}(t) + \La_{G,2}(t)] + Cn^2\int_{t_i}^t
|t-s|^{-\sigma'}[\La_{L,1} +\La_{L,2}](s ) ds,
\end{flalign*}
which is bounded by $\La_{G,1}(t) + \La_{G,2}(t)$. Similarly,
$\norm{e^{\bL(t-t_i)}P_{\pm}[\eta_2]}_{L^2\loc}$ is bounded by
\begin{flalign*}
& \leq \norm{e^{J(H_0-E)(t-t_i)}[\eta_2]}_{L^2\loc} + Cn^2\int_{t_i}^t \min\{ |t-s|^{-3/2}, |t-s|^{-\sigma'}\}\norm{e^{J(H_0-E)s}[\eta_2]}_{L^2\loc}ds\\
& \leq \frac{1}{2}[\La_{L,1}(t) + \La_{L,2}(t)] + Cn^2\int_{t_i}^t
\min\{ |t-s|^{-3/2}, |t-s|^{-\sigma'}\}[\La_{L,1} + \La_{L,2}](s)
ds,
\end{flalign*}
which is bounded by $\La_{L,1}(t) + \La_{L,2}(t)$. Summing the
estimates we get the Lemma. \myendproof
%======================================================

%====================================================================================
\medskip

Denote
\begin{equation} \begin{split}
\hat{\rho}(t) &= \rho(t-t_i) =[\rho_0^{-2} +\gamma_0 n^2
(t-t_i)]^{-1/2},
\\
\delta_8(t)&= n^{-\frac{2}{3}(1-\de)}(t-t_o)^{-6/p} +
n^{6}\wei{t-t_i}^{-6/p} \leq o(1)n\hro(t)^2,
\end{split}
\end{equation}
and
\begin{equation}
M_{T}^* : = \sup_{t_i \leq t \leq T}\max \left \{
\begin{array} {ll}
& \hro(t)^{-1}|z_H(t)|, \quad  [2D \hro(t)]^{-1}|a(t)|, \\
& [\La_G(t) + n^{7/9}\hro(t)^{5/3}]^{-1} \norm{\eta}_{L^p}, \\
& [\La_L(t) + \La_G^{2}(t) + n^{-\al}\hro(t)^{3} +\delta_8(t)]^{-1} \norm{\eta^{(3)}}_{L^2\loc}
\end{array} \right \}.
\end{equation}
%==============================================

\begin{proposition} \label{gs.con}
Suppose  for $T \ge t_i$ we have $n(T)/n(t_i) \in (\frac 12, \frac
32)$ and $M^*_T \leq 3$. Then we have  $M_{T}^* \leq \frac52$ and
$n(T)/n(t_i) \in (\frac{3}{4}, \frac{5}{4})$.
\end{proposition}
%=====================================================

This Proposition implies Theorem \ref{mainthm} in the case $k=0$,
see e.g.~\cite{BP2,TY1,C2,Tsai}.

\medskip

{\it Proof of Proposition \ref{gs.con}}.\quad The condition $M^*_T
\leq 3$ means, $\forall \ t_i \leq t <T$,
\begin{equation} \label{MT*.w}
\begin{split}
& z_H(t) \leq 3\hro(t), \quad |a(t)| \leq 6D \hro(t), \quad \norm{\eta}_{L^p} \leq 3[\La_G(t) + n^{7/9}\hro(t)^{5/3}],\\
& \norm{\eta^{(3)}}_{L^2\loc} \leq 3[\La_L(t) + \La_G^{2} +
n^{-\al}\hro(t)^{3} +\delta_8(t)].
\end{split}
\end{equation}
Let $\de_9(t) := \La_{L,2}(t) + \La_{L,3}(t)$. Note
\begin{equation}
\La_{L,1}(t) +\La_{G}^2(t) +\de_8(t) \leq o(1)n\hro(t)^2, \quad
\sum_{j=1}^3[\La_{L,j}(t)+ \La_{G,j}(t) ]\leq o(1)\hro(t).
\end{equation}
Thus
\begin{equation} \label{eta3.sim.est}
\norm{\eta^{(3)}}_{L^2\loc} \leq o(1)n\hro(t)^2 + 3\de_9(t), \quad
\norm{\eta}_{L^2\loc} \leq Cn\hro(t)^2 +3\de_9(t).
\end{equation}
We also have
\begin{equation} \label{etaL2Lp}
\norm{\eta}_{L^2\loc} \leq o(1)\hro(t), \quad \norm{\eta}_{L^p} \leq o(1)\hro(t), \quad t \in [t_i, T].
\end{equation}
Recall that $X$ and $\wt{X}$ are defined in \eqref{XwtX.def}. From
\eqref{MT*.w}, \eqref{eta3.sim.est} and
\begin{equation} \label{ti.nonloc}
\norm{\eta^3}_{L^1\loc} \leq
\norm{\eta}_{L^2\loc}^{\frac{2(p-3)}{p-2}}
\norm{\eta}_{L^p}^{\frac{p}{p-2}} \leq
\norm{\eta}_{L^2\loc}\norm{\eta}_{L^p}^2,
\end{equation}
we have
\begin{equation} \label{wtX.ti.est}
\wt{X} \les n\hro(t)^4 + [\hro^2 + n\de_9]\de_9(t), \quad X \les
n^2\hro(t)^3 + n\hro\de_9.
\end{equation}
On the other hand, from \eqref{basic.est}, \eqref{Z.est}, \eqref{eta3.sim.est} and \eqref{wtX.ti.est}, we also obtain
\begin{equation} \label{p.theta.dot}
|\dot{\theta}| = |F_\theta| \les \hro^2 + n^{-1}X \les \hro(t)^2,
\quad |\dot{p}_k| = |Z_k| \les n\hro(t)^2, \quad k \geq 1.
\end{equation}
The proof of this proposition is divided into two lemmas:
\myendproof
%========================================
\begin{lemma} Suppose for some $T > t_i$ that $M_T^* \leq 3$. Then, $t \in [t_i, T]$, we have
\begin{flalign*}
& \norm{\eta(t)}_{L^p} \leq \frac52[\La_G(t) + n^{7/9}\hro(t)^{5/3}],\\
& \norm{\eta^{(3)}(t)}_{L^2\loc} \leq \frac52[\La_L(t) + \La_G^{2} +
n^{-\al}\hro(t)^{3} +\delta_8(t)].
\end{flalign*}
\end{lemma}
%================================================
\myproof Recall $[\eta] = e^{-i\theta}\eta_- + e^{i\theta}\eta_+$.
By \eqref{wt-eta.eqn},
\begin{equation}
\eta_\pm(t) = e^{\bL(t-t_i)}\eta_{\pm}(t_i) + \int_{t_i}^t
e^{\bL(t-s)}P_\pm\{F_{L\pm} +e^{\mp i \th}J[F]\}(s) ds.
\end{equation} From \eqref{FL}, \eqref{basic.est},
\eqref{wtX.ti.est}, \eqref{p.theta.dot}, \eqref{eq:5.20} and as in
\eqref{FL1.est} and \eqref{FN1.est}, we get
\begin{equation}
\norm{F_{L\pm}}_{L^{p'}} + \norm{F}_{L^{p'}} \leq
Cn\hro(t)^2+o(1)\norm{\eta}_{L^p}^{\frac {p+2}{p-2}}.\end{equation}
By Lemma \ref{OG.ti.ln}, $\frac {p+2}{p-2}>2$ and
$\norm{\eta}_{L^p}^2 \le o(1) n \hro^2$,
\begin{flalign*}
\norm{\eta_{\pm}}_{L^p} & \leq \La_G(t)
+ \int_{t_i}^{t}(t-s)^{-\sigma'}n\hro(s)^2 ds \\
& \leq \La_G(t) + C n^{-1}(\Delta t)^{-\al}[\Delta t + t-t_i]^{-\sigma' +\al} \leq \La_{G}(t) + o(1) n^{7/9}\hro(t)^{5/3}.
\end{flalign*}
Above we have used \eqref{int-est1} with $\al=\si'-5/6>7/18$. So, we
get
\begin{equation}
\norm{\eta}_{L^p} \leq \norm{\eta_-}_{L^p} + \norm{\eta_+}_{L^p}\leq
2[\La_{G}(t) + n^{7/9}\hro(t)^{5/3}].
\end{equation}

To estimate $\eta^{(3)}$ in $L^2_{loc}$,  write $\eta^{(3)} =
e^{-i\theta}\eta_-^{(3)} + e^{i\theta}\eta_+^{(3)}$ with
$\eta_\pm^{(3)} =\sum_{j=1}^4 \eta_{\pm,j}^{(3)}$ as in
\eqref{wt-eta.eqn2} and \eqref{wt-eta23.def}. Using \eqref{MT*.w},
\eqref{wtX.ti.est}, \eqref{p.theta.dot} and the argument of
\eqref{eta(3).12.est} and \eqref{eta(3).3.est}, we obtain
\begin{equation} \label{eta3.123}
\sum_{j=1}^{3} \norm{\eta_{\pm,j}^{(3)}(t)}_{L^2\loc} \leq \La_L(t) + Cn[\rho(t)^3 + \rho_0^2\wei{t-t_i}^{-3/2}].
\end{equation}
To estimate $\eta^{(3)}_{\pm,4}$, we write $\eta^{(3)}_{\pm,4} = \eta^{(3)}_{\pm,4,1} + \eta^{(3)}_{\pm,4,2}$ with
\begin{flalign*}
\eta^{(3)}_{\pm,4,1}(t) & = \int_{t_i}^t e^{\bL(t-s)}P_{\pm}\{F_{L\pm} +e^{\mp i\theta} J[F-F_1 -\kappa |\eta|^2\eta]\}(s) ds, \\
\eta^{(3)}_{\pm,4,2}(t) &= \int_{t_i}^t
e^{\bL(t-s)}P_{\pm}\{e^{\mp i\theta} J[\kappa  |\eta|^2\eta]\}(s)
ds\}.
\end{flalign*}
Note that $\{F_{L\pm} +e^{\mp i\theta} J[F-F_1 -|\eta|^2\eta]\}$
is localized and
\begin{equation} \norm{F_{L\pm}}_{L^1\cap L^{p'}} + \norm{F-F_1 -|\eta|^2\eta}_{L^1\cap L^{p'}} \leq C[\hro(t)^3 + n\La_G^2(t)]. \end{equation}
So, we get
\begin{equation} \label{eta3.4.1}
\begin{split}
\norm{\eta^{(3)}_{\pm,4,1}(t)}_{L^2\loc} & \leq  C\int_{t_i}^{t} \min\{|t-s|^{-3/2}, |t-s|^{-\sigma'}\}[\hro(s)^3 +  n\La_G^2(s)]ds\\
& \leq C[\hro(t)^3 +  n\La_G^2(t)].
\end{split}
\end{equation}
As in \eqref{l-choice}, we have
\begin{equation}\norm{|\eta|^2\eta}_{L^{\frac{2p}{p+4}} \cap L^{p'}} \leq o(1)\norm{\eta}_{L^p}^2 \leq o(1)[\La_{G}^2 + n^{\frac{14}{9}} \hro(t)^{\frac{10}{3}}].  \end{equation}
So, it follows as in \eqref{L2loc.xi35} that
\begin{equation} \label{eta3.4.2}
\norm{\eta^{(3)}_{\pm,4,2}(t)}_{L^2\loc} \leq o(1)[\La^2_G(t) + \hro(t)^3 + \delta_8(t)].
\end{equation}
Collecting \eqref{eta3.123}, \eqref{eta3.4.1}, and \eqref{eta3.4.2}, we obtain the second estimate of the lemma.
\myendproof

%==================================
\begin{lemma} \label{BZ}
Suppose for some $T > t_i$ that $M_T^* \leq 3$. Then, for $t \in
[t_i, T]$, we have
\begin{equation} z_H(t) \leq 2 \hro(t), \quad |a(t)| \leq 4 D\hro(t)^2, \quad n(T)/n(t_i) \in (\frac{3}{4}, \frac{5}{4}), \quad \forall k \geq 1. \end{equation}
\end{lemma}
\myproof From \eqref{etaL2Lp} and \eqref{MT*.w}, we can apply
Lemma \ref{L-NF} with $\beta = \hro(t)$. So, for each $k \geq 1$, we can find $q_k$ such that
\begin{equation}
\dot{q}_k = \sum_{l \geq 1}D_{kl}|q_l|^2 q_k + Y_k q_k + g_k, \quad
|q_k - p_k| \leq Cn\hro(t)^2.
\end{equation}
Moreover, $\Re(Y_k) =0$, $\Re D_{kl}\leq -\frac{\gamma_0n^2}{2}$ and
\begin{equation} |g_k| \leq C[n\hro(t)^4 + n^3\hro(t)\norm{\eta}_{L^2\loc} + n\hro(t)\norm{\eta^{(3)}}_{L^2\loc} + \wt{X}]. \end{equation}
So, from \eqref{eta3.sim.est} and \eqref{wtX.ti.est}, we get
\begin{equation}
|g_k| \leq o(1)n^2\hro(t)^3 + Cn\hro(t)\de_9(t).
\end{equation}
Note that
\begin{equation} \label{gk.ti.est}
\int_{t_i}^{t_i+n^{-3}} \de_{9}(t) ds \leq o(1) ; \quad
\de_{9}(t)\leq o(1)n\hro(t)^2 \quad \forall t > t_i + n^{-3}.
\end{equation}
By the argument in the proof of Lemma \ref{bsz.est}, we obtain
\begin{equation} |z_k(t)| \leq z_H(t) \leq  2\hro(t), \quad \forall t \in [t_i, T]. \end{equation}
On the other hand, from Lemma \ref{b.NF}, we have
\begin{flalign*}
& \dot{\wt{b}} = \sum_{k, l \geq 1}B_{kl}|z_k|^2|z_l|^2 + g_b, \quad |b -\wt{b}| \leq Cn\hro^2, \\
& |g_b| \leq C[n^3\hro^4(t) + n^2\norm{\eta}_{L^2\loc} + n\norm{\eta^3}_{L^1\loc} + n\hro^2\norm{\eta^{(3)}}_{L^2\loc}].
\end{flalign*}
Again, by using \eqref{eta3.sim.est} and \eqref{ti.nonloc}, we get
\begin{equation}
|g_b(s)| \leq o(1)n^2\hro(t)^4 + Cn\de_9(t)[n\de_9 +\hro(t)^2].
\end{equation}
If $ T \geq t > t_i +n^{-3}$, we get $\de_9(t) \leq o(1)n\hro(t)^2$. So,
\begin{equation}
\int_{t_i +n^{-3}}^{t}|g_b(s)|ds \leq o(1)\int_{t_i+n^{-3}}^t
n^2\hro(s)^4 ds \leq o(1)\hro(t)^2.
\end{equation}
For $t_i \leq t < \min(n^{-3}, T)$, we have $|g_b(s)| \leq o(1)
[n^2\rho_0^4 + n^2\rho_0 \de_9(s)]$ and, using \eqref{gk.ti.est},
\begin{equation}
\int_{t_i}^{t}|g_b(s)|ds \leq o(1)\int_{t_i}^{t_i+n^{-3}}
[n^2\rho_0^4 + n^2\rho_0 \de_9(s)] ds \leq o(1)\rho_0^2 \sim
o(1)\hro(t)^2.
\end{equation}
So,
we get
\begin{equation}
\int_{t_i}^{t}|g_b|(s) ds \leq o(1)\hro(t)^2, \quad \forall \ t_i \leq t \leq T.
\end{equation}
From this and as in the proof of Lemma \ref{bsz.est}, we also get
\begin{equation}
|a(t)| \leq 4D\rho(t)^2, \quad n(T)/n(t_i) \in (\frac{3}{4}, \frac{5}{4}).
\end{equation}
This completes the proof of Lemma \ref{BZ}. \myendproof
%=============================================================================================================

%====================================================================

\section*{Acknowledgments}
We thank S.~Gustafson for his constant interest in this work. Part
of this work was conducted while the first author visited the
University of British Columbia, by the support of the 21st century
COE program and the Kyoto University Foundation, and while the
second author was a postdoctoral fellow at  the University of
British Columbia. The research of Nakanishi was partly supported
by the JSPS grant no.~15740086. The research of Tsai was partly
supported by the NSERC grant no. 261356-08.

%====================================================================

%\bigskip

\noindent{\it Kenji Nakanishi}\quad Department of Mathematics,
Kyoto University, Kyoto 606-8502, Japan;
n-kenji@math.kyoto-u.ac.jp

\medskip

\noindent{\it Tuoc Van Phan}\quad Department of Mathematics,
University of Tennessee, Knoxville, TN 37996, USA;
phan@math.utk.edu

\medskip

\noindent{\it Tai-Peng Tsai} \quad Department of Mathematics,
University of British Columbia, Vancouver, BC V6T 1Z2, Canada;
ttsai@math.ubc.ca

\end{document}